  \documentstyle[12pt]{article}
    \title{{\bf  A theory of tensor products for module categories for
a vertex operator algebra, IV}}
    \author{Yi-Zhi Huang}
    \date{}
    
\begin{document}
    \bibliographystyle{alpha}
    \maketitle

    \input amssym.def
    \input amssym
    \newtheorem{rema}{Remark}[section]
    \newtheorem{propo}[rema]{Proposition}
    \newtheorem{theo}[rema]{Theorem}
   \newtheorem{defi}[rema]{Definition}
    \newtheorem{lemma}[rema]{Lemma}
    \newtheorem{corol}[rema]{Corollary}
     \newtheorem{exam}[rema]{Example}
	\newcommand{\ba}{\begin{array}}
	\newcommand{\ea}{\end{array}}
        \newcommand{\be}{\begin{equation}}
        \newcommand{\ee}{\end{equation}}
	\newcommand{\bea}{\begin{eqnarray}}
	\newcommand{\eea}{\end{eqnarray}}
	\newcommand{\nno}{\nonumber}
	\newcommand{\lbar}{\bigg\vert}
	\newcommand{\p}{\partial}
	\newcommand{\dps}{\displaystyle}
	\newcommand{\bra}{\langle}
	\newcommand{\ket}{\rangle}
 \newcommand{\res}{\mbox{ \rm Res}}
\newcommand{\wt}{\mbox{ \rm wt}\ }
 \newcommand{\pf}{{\it Proof}\hspace{2ex}}
 \newcommand{\epf}{\hspace{2em}$\Box$}
 \newcommand{\epfv}{\hspace{1em}$\Box$\vspace{1em}}
	\renewcommand{\hom}{\mbox{\rm Hom}}

\hyphenation{Phil-a-del-phia}

 \makeatletter
\newlength{\@pxlwd} \newlength{\@rulewd} \newlength{\@pxlht}
\catcode`.=\active \catcode`B=\active \catcode`:=\active \catcode`|=\active
\def\sprite#1(#2,#3)[#4,#5]{
   \edef\@sprbox{\expandafter\@cdr\string#1\@nil @box}
   \expandafter\newsavebox\csname\@sprbox\endcsname
   \edef#1{\expandafter\usebox\csname\@sprbox\endcsname}
   \expandafter\setbox\csname\@sprbox\endcsname =\hbox\bgroup
   \vbox\bgroup
  \catcode`.=\active\catcode`B=\active\catcode`:=\active\catcode`|=\active
      \@pxlwd=#4 \divide\@pxlwd by #3 \@rulewd=\@pxlwd
      \@pxlht=#5 \divide\@pxlht by #2
      \def .{\hskip \@pxlwd \ignorespaces}
      \def B{\@ifnextchar B{\advance\@rulewd by \@pxlwd}{\vrule
         height \@pxlht width \@rulewd depth 0 pt \@rulewd=\@pxlwd}}
      \def :{\hbox\bgroup\vrule height \@pxlht width 0pt depth
0pt\ignorespaces}
      \def |{\vrule height \@pxlht width 0pt depth 0pt\egroup
         \prevdepth= -1000 pt}
   }
\def\endsprite{\egroup\egroup}
\catcode`.=12 \catcode`B=11 \catcode`:=12 \catcode`|=12\relax
\makeatother

\def\hboxtr{\FormOfHboxtr} 
\sprite{\FormOfHboxtr}(25,25)[0.5 em, 1.2 ex] 

:BBBBBBBBBBBBBBBBBBBBBBBBB |
:BB......................B |
:B.B.....................B |
:B..B....................B |
:B...B...................B |
:B....B..................B |
:B.....B.................B |
:B......B................B |
:B.......B...............B |
:B........B..............B |
:B.........B.............B |
:B..........B............B |
:B...........B...........B |
:B............B..........B |
:B.............B.........B |
:B..............B........B |
:B...............B.......B |
:B................B......B |
:B.................B.....B |
:B..................B....B |
:B...................B...B |
:B....................B..B |
:B.....................B.B |
:B......................BB |
:BBBBBBBBBBBBBBBBBBBBBBBBB |

\endsprite

\def\shboxtr{\FormOfShboxtr} 
\sprite{\FormOfShboxtr}(25,25)[0.3 em, 0.72 ex] 

:BBBBBBBBBBBBBBBBBBBBBBBBB |
:BB......................B |
:B.B.....................B |
:B..B....................B |
:B...B...................B |
:B....B..................B |
:B.....B.................B |
:B......B................B |
:B.......B...............B |
:B........B..............B |
:B.........B.............B |
:B..........B............B |
:B...........B...........B |
:B............B..........B |
:B.............B.........B |
:B..............B........B |
:B...............B.......B |
:B................B......B |
:B.................B.....B |
:B..................B....B |
:B...................B...B |
:B....................B..B |
:B.....................B.B |
:B......................BB |
:BBBBBBBBBBBBBBBBBBBBBBBBB |

\endsprite

\vspace{2em}

\begin{abstract}
This is the fourth part of a series of papers developing a tensor
product theory of modules for a vertex operator algebra. In this
paper, We establish the associativity of $P(z)$-tensor products for
nonzero complex numbers $z$ constructed in Part III of the present
series under suitable conditions. The associativity isomorphisms
constructed in this paper are analogous to associativity isomorphisms
for vector space tensor products in the sense that it relates the
tensor products of three elements in three modules taken in different
ways.  The main new feature is that they are controlled by the
decompositions of certain spheres with four punctures into spheres
with three punctures using a sewing operation.  We also show that
under certain conditions, the existence of the associativity
isomorphisms is equivalent to the associativity (or (nonmeromorphic)
operator product expansion in the language of physicists) for the
intertwining operators (or chiral vertex operators). Thus the
associativity of tensor products provides a means to establish the
(nonmeromorphic) operator product expansion.
\end{abstract}

\vspace{2em}

The present paper (Part IV) is the fourth in a series of papers developing
a theory of tensor products of modules for a vertex operator algebra. An
overview of the
theory being developed has been given in
\cite{HL5} and
the reader is referred to it for the motivation and the description of the main
results.

In  Part I (\cite{HL3}), the
notions of $P(z)$- and $Q(z)$-tensor product for any nonzero complex number $z$
 are introduced and two
constructions
of a $Q(z)$-tensor product are given based on certain results
proved in  Part II (\cite{HL4}).
In  Part III (\cite{HL6}), the notion of $P(z)$-tensor product
is discussed in the same way as in Section 4 of Part I for that
of $Q(z)$-tensor product, and two constructions
of a $P(z)$-tensor product are given
using the results for the $Q(z)$-tensor product.
Part III (\cite{HL6}) also contains a brief
description of the results in \cite{HL3} and \cite{HL4}.
In the present paper,
the associativity for  $P(\cdot)$-tensor products is formulated
and the associativity isomorphisms are constructed
under certain assumptions on the vertex operator algebra and
on the products or iterates of
two intertwining operators for the vertex
operator algebra. These assumptions are satisfied by familiar examples
and will be discussed in separate papers on applications of the theory
of tensor products developed in the present series of papers. We also show
that when the vertex operator algebra is rational and products of
two intertwining
operators are convergent in a suitable region, the existence of the
associativity isomorphisms
is equivalent to the associativity (or (nonmeromorphic) operator product
expansion in the language of physicists)  for the intertwining operators
(or chiral vertex operators).
See Theorems 14.11, 16.3 and 16.5 for the precise statements
of the main results
of the present paper.

Our conventions in this paper is the same as those in \cite{HL3},
\cite{HL4} and \cite{HL6}.  We also
add $y$ to our list of formal variables. The symbols $z_{1}, z_{2},
\dots$ will also denote nonzero complex numbers.
We fix a vertex operator algebra $V$.
The numberings of sections, formulas, etc.,
continues those of Part I, Part II and Part III.

Part IV is organized as follows:
 The associativity isomorphisms are constructed in
Section 14, based on certain assumptions and some lemmas. We also
prove in Section 14 that the existence of the associativity isomorphisms
is equivalent to the associativity of the intertwining operators.
The lemmas used in
Section 14 are proved in
Section 15. In Section 16,  we give some conditions and show that they imply
the assumptions used in the construction of the associativity isomorphisms.

\paragraph{Acknowledgments}
The present paper is one of the papers resulted {}from a long-term
project jointly
with James Lepowsky developing a tensor product theory of modules
for a vertex operator algebra.
I would like to express my gratitude to him for
collaborations and many discussions.
This work has been supported in part by NSF grants
DMS-9104519 and DMS-9301020 and by DIMACS, an
NSF Science and Technology Center funded under contract STC-88-09648.

\setcounter{section}{13}
\renewcommand{\theequation}{\thesection.\arabic{equation}}
\renewcommand{\therema}{\thesection.\arabic{rema}}
\setcounter{equation}{0}
\setcounter{rema}{0}

\section{Associativity isomorphisms for $P(\cdot )$-ten\-sor products}

In this section, we construct the associativity isomorphisms for
$P(\cdot)$-tensor products. To discuss the associativity we need to consider
the compositions of one $P(z_{1})$-intertwining map and one
$P(z_{2})$-intertwining map for suitable nonzero complex numbers $z_{1}$
and $z_{2}$. Since intertwining
maps are maps {}from $W_{1}\otimes W_{2}$ to $\overline{W}_{3}$, not to
$W_{3}$,
we have to give a precise meaning of the compositions of one
$P(z_{1})$-intertwining maps and one $P(z_{2})$-intertwining map.
Compositions of maps of this type have been defined precisely in \cite{HL1}
and \cite{HL2}. Here we repeat the definition for the concrete examples of
intertwining maps. Let
$F_{1}$ and $F_{2}$ be  $P(z_{1})$- and $P(z_{2})$-intertwining maps
of type ${W_{4}}\choose {W_{1}W_{5}}$ and ${W_{5}}\choose {W_{2}W_{3}}$,
respectively. If for any $w_{(1)}\in W_{1}$, $w_{(2)}\in W_{2}$,
$w_{(3)}\in W_{3}$ and $w'_{(4)}\in W'_{4}$, the series
\begin{equation}
\sum_{n\in {\Bbb C}}\langle w'_{(4)}, F_{1}(w_{(1)}\otimes
P_{n}(F_{2}(w_{(2)}\otimes w_{(3)})))\rangle_{W_{4}}
\end{equation}
is absolutely convergent, where for any $n\in {\Bbb C}$,
$P_{n}$ is the projection map {}from a
module to its homogeneous subspace of weight $n$, we say that {\it the product
of $F_{1}$ and $F_{2}$ exists} and we call the map {}from $W_{1}\otimes
W_{2}\otimes W_{3}$ to $\overline{W}_{4}$ defined by the  limits above the
{\it product} of $F_{1}$ and $F_{2}$. Similarly, let $F_{3}$ and
$F_{4}$ be $P(z_{3})$- and $P(z_{4})$-intertwining maps of type
${W_{5}}\choose {W_{1}W_{2}}$ and ${W_{4}}\choose {W_{5}W_{3}}$,
respectively. If for any $w_{(1)}\in W_{1}$, $w_{(2)}\in W_{2}$,
$w_{(3)}\in W_{3}$ and $w'_{(4)}\in W'_{4}$, the series
\begin{equation}
\sum_{n\in {\Bbb C}}\langle w'_{(4)}, F_{4}(P_{n}
(F_{3}(w_{(1)}\otimes w_{(2)}))\otimes w_{(3)})\rangle_{W_{4}}
\end{equation}
is absolutely convergent, we say that the {\it iterate
of $F_{3}$ and $F_{4}$ exists} and we call the map {}from $W_{1}\otimes
W_{2}\otimes W_{3}$ to $\overline{W}_{4}$ defined by the  limits above the
{\it iterate} of $F_{3}$ and $F_{4}$.

Recall that for any nonzero complex number $z$ and any integer $p$,
there is an isomorphism between
the space of $P(z)$-intertwining  maps and the space of intertwining operators
of the same type (see \cite{HL6}). Take $p=0$ and let
${\cal Y}_{1}$, ${\cal Y}_{2}$,
${\cal Y}_{3}$
and ${\cal Y}_{4}$ be the intertwining operators correspond to $F_{1}$,
$F_{2}$, $F_{3}$ and $F_{4}$, respectively.
(In this paper, we shall always use the isomorphism between the space of
$P(z)$-intertwining maps and the space of intertwining operators of
the same type associated to $p=0$.) Then for $i=1, 2, 3, 4$, we have
\begin{equation}
F_{i}(\cdot \otimes \cdot)={\cal Y}_{i}(\cdot, x)\cdot\lbar_{x^{n}=e^{n\log
z_{i}}, \;n\in {\Bbb C}}.
\end{equation}
The absolute convergence of the series (14.1) and (14.2) are equivalent to the
absolute convergence of the series
\begin{equation}
\langle w'_{(4)}, {\cal Y}_{1}(w_{(1)}, x_{1})
{\cal Y}_{2}(w_{(2)}, x_{2})w_{(3)}\rangle_{W_{4}}\lbar_{x^{n}_{1}=e^{n\log
z_{1}},\;
x^{n}_{2}=e^{n\log z_{2}}, \;n\in {\Bbb C}}
\end{equation}
and
\begin{equation}
\langle w'_{(4)}, {\cal Y}_{4}({\cal
Y}_{3}(w_{(1)},
x_{1})w_{(2)}, x_{2})w_{(3)}
\rangle_{W_{4}}\lbar_{x^{n}_{1}=e^{n\log z_{3}},\;
x^{n}_{2}=e^{n \log z_{4}}, \;n\in {\Bbb C}},
\end{equation}
 respectively. For convenience, we shall write the substitutions, say,
$$\lbar_{x^{n}_{1}=e^{n\log
z_{1}},\;
x^{n}_{2}=e^{n\log z_{2}},\; n\in {\Bbb C}}$$ simply as
$\lbar_{x_{1}=z_{1}, \;x_{2}=z_{2}}$.

Let $V$ be a vertex operator algebra such that for any $V$-modules $W_{1}$,
$W_{2}$, $W_{3}$, $W_{4}$ and $W_{5}$, any two nonzero complex
numbers $z_{1}$ and $z_{2}$ satisfying $|z_{1}|>|z_{2}|>0$ and any
$P(z_{1})$-intertwining map $F_{1}$ and $P(z_{2})$-intertwining map $F_{2}$
as above, (14.1) is absolutely convergent for all $w_{(1)}$,
$w_{(2)}$, $w_{(3)}$ and
$w'_{(4)}$.  Then (14.4) is also absolutely convergent.
Recall that we have an isomorphism
$\Omega_{-1}: {\cal V}_{W_{1} W_{5}}^{W_{4}}\to {\cal V}_{W_{5}W_{1}}^{W_{4}}$
defined by
$$\Omega_{-1}({\cal Y})(w_{(5)}, x)w_{(1)}=e^{xL(-1)}{\cal Y}(w_{(1)},
e^{-\pi i}x)w_{(5)}$$
for $w_{(1)}\in W_{1}$, $w_{(5)}\in W_{5}$ and ${\cal Y}\in
{\cal V}_{W_{1} W_{5}}^{W_{4}}$ (see \cite{FHL} and \cite{HL4}).
Its inverse is $\Omega_{0}: {\cal V}_{W_{5}W_{1}}^{W_{4}}
\to {\cal V}_{W_{1} W_{5}}^{W_{4}}$ defined by
$$\Omega_{0}({\cal Y})(w_{(1)}, x)w_{(5)}=e^{xL(-1)}{\cal Y}(w_{(5)},
e^{\pi i}x)w_{(1)}$$
for $w_{(1)}\in W_{1}$, $w_{(5)}\in W_{5}$ and ${\cal Y}\in
{\cal V}_{W_{5} W_{1}}^{W_{4}}$.
Let $F_{3}$ and
$F_{4}$ be $P(z_{3})$- and $P(z_{4})$-intertwining maps, respectively,
as above,
and ${\cal Y}_{3}$ and ${\cal Y}_{4}$
the corresponding intertwining operators. Thus
\begin{eqnarray}
\lefteqn{\langle w'_{(4)},
{\cal Y}_{4}({\cal Y}_{3}(w_{(1)},
x_{1})w_{(2)}, x_{2})w_{(3)}
\rangle_{W_{4}}\lbar_{x_{1}= z_{3},\;
x_{2}= z_{4}}=}\nno\\
&&=\langle w'_{(4)},
\Omega_{0}(\Omega_{-1}({\cal Y}_{4}))({\cal Y}_{3}(w_{(1)},
x_{1})w_{(2)}, x_{2})w_{(3)}
\rangle_{W_{4}}\lbar_{x_{1}= z_{3},\;
x_{2}= z_{4}}\nno\\
&&=\langle w'_{(4)},
e^{x_{2}L(-1)}\Omega_{-1}({\cal Y}_{4})
(w_{(3)}, e^{\pi i}x_{2})
{\cal Y}_{3}(w_{(1)},
x_{1})w_{(2)}\rangle_{W_{4}}\lbar_{x_{1}= z_{3},\;
x_{2}= z_{4}}\nno\\
&&=\langle e^{x_{2}L(1)}w'_{(4)},
\Omega_{-1}({\cal Y}_{4})
(w_{(3)}, e^{\pi i}x_{2})
{\cal Y}_{3}(w_{(1)},
x_{1})w_{(2)}\rangle_{W_{4}}\lbar_{x_{1}= z_{3},\;
x_{2}= z_{4}}\nno\\
&&=\langle e^{z_{4}L(1)}w'_{(4)},
\Omega_{-1}({\cal Y}_{4})
(w_{(3)}, x_{2})
{\cal Y}_{3}(w_{(1)},
x_{1})w_{(2)}\rangle_{W_{4}}\lbar_{x_{1}= z_{3},\;
x_{2}=e^{\pi i} z_{4}},\nno\\
&&
\end{eqnarray}
where for convenience, we have used $\lbar_{x_{1}= z_{3},\;
x_{2}=e^{\pi i} z_{4}}$ to denote $$\lbar_{x^{n}_{1}= e^{n\log z_{3}},\;
x^{n}_{2}=e^{n\pi i} e^{n\log z_{4}}, \; n\in {\Bbb C}}$$ (and we shall use
the similar notations below).
Using the isomorphism between the space of $P(-z_{4})$-intertwining
maps and the space of intertwining operators and the isomorphism
between the space of
$P(z_{3})$-intertwining maps and the space of intertwining operators,
we see that the right-hand side of (14.6) is equal to the product of a
$P(-z_{4})$-intertwining
map and a $P(z_{3})$-intertwining map evaluated at $w_{(3)}\otimes
w_{(1)}\otimes w_{(2)}\in W_{3}\otimes W_{1}\otimes W_{2}$
and paired with $e^{z_{4}L(1)}w'_{(4)}\in W'_{4}$. By assumption, the
right-hand side of (14.6) is convergent absolutely when
$|-z_{4}|>|z_{3}|>0$
or equivalently when $|z_{4}|>|z_{3}|>0$.
Thus we have proved half
of the following result:

\begin{propo}
For any vertex operator algebra $V$,  the following two properties
 are equivalent:
\begin{enumerate}
\item For any $V$-modules $W_{1}$,
$W_{2}$, $W_{3}$, $W_{4}$ and $W_{5}$, any nonzero complex
numbers $z_{1}$ and $z_{2}$ satisfying $|z_{1}|>|z_{2}|>0$ and any
$P(z_{1})$-intertwining map $F_{1}$ of type
${W_{4}}\choose {W_{1}W_{5}}$ and $P(z_{2})$-intertwining map $F_{2}$
of type ${W_{5}}\choose {W_{2}W_{3}}$, the product of $F_{1}$ and $F_{2}$
exists.

\item For any $V$-modules $W_{1}$,
$W_{2}$, $W_{3}$, $W_{4}$ and $W_{5}$, any nonzero complex
numbers $z_{3}$ and $z_{4}$ satisfying $|z_{4}|>|z_{3}|>0$ and any
$P(z_{3})$-intertwining map $F_{3}$ of type
${W_{4}}\choose {W_{5}W_{3}}$ and $P(z_{4})$-intertwining map $F_{4}$
of type ${W_{5}}\choose {W_{1}W_{2}}$,
the iterate of $F_{3}$ and $F_{4}$ exists. \epf
\end{enumerate}

\end{propo}

The  other half of this result can be proved similarly.

We need the following result on products of the $\delta$-function:

\begin{lemma}
For any two nonzero complex numbers $z_{1}$ and $z_{2}$ satisfying
\begin{equation}
|z_{1}|>|z_{2}|>|z_{1}-z_{2}|>0,
\end{equation}
 we have
\begin{eqnarray}
\lefteqn{x_{1}^{-1}\delta\left(\frac{
x^{-1}_{0}-z_{1}}
{x_{1}}\right)x_{2}^{-1}\delta\left(\frac{x^{-1}_{0}-z_{2}}{x_{2}}\right)=}
\nno\\
&&\hspace{2em}=x^{-1}_2\delta\left(\frac{x^{-1}_0-z_2}{x_2}\right)
x^{-1}_1\delta\left(\frac{x_2-(z_1-z_{2})}{x_1}\right),
\end{eqnarray}
\begin{eqnarray}
\lefteqn{z^{-1}_1\delta\left(\frac{x^{-1}_0-x_1}{z_1}\right)x^{-1}_2
\delta\left(\frac{x^{-1}_0-z_2}{x_2}\right)=}\nno\\
&&\hspace{2em}=z^{-1}_2\delta\left(\frac{x^{-1}_0-x_2}{z_2}\right)
 (z_1-z_2)^{-1}
\delta\left(\frac{x_2-x_1}{z_1-z_2}\right),
\end{eqnarray}
\begin{eqnarray}
\lefteqn{x^{-1}_1\delta\left(\frac{z_1-x^{-1}_0}{x_1}\right)
x^{-1}_2\delta\left(\frac{z_2-x^{-1}_0}{x_2}\right)=}\nno\\
&&\hspace{2em}=x^{-1}_2\delta\left(\frac{z_2-x^{-1}_0}{x_2}\right)
x^{-1}_1\delta\left(\frac{x_2-(z_1-z_2)}{x_1}\right),
\end{eqnarray}
\begin{eqnarray}
\lefteqn{x^{-1}_1\delta\left(\frac{z_1-x^{-1}_0}{x_1}\right)
z^{-1}_2\delta\left(\frac{x^{-1}_0-x_2}{z_2}\right)=}\nno\\
&&\hspace{2em}=z^{-1}_2\delta\left(\frac{x^{-1}_0-x_2}{z_2}\right)x^{-1}_1
\delta\left(\frac{z_1-z_2-x_2}{x_1}\right).
\end{eqnarray}
The left-hand sides and the right-hand sides of these formulas are
formal series in $x_{0}, x_{1}, x_{2}$ whose coefficients converge
absolutely when $|z_{1}|>|z_{2}|>0$ and when
$|z_{2}|>|z_{1}-z_{2}|>0$, respectively, to rational functions of
$z_{1}$ and $z_{2}$ with the only possible poles $z_{1}=\infty, 0$,
$z_{2}=\infty, 0$ and $z_{1}=z_{2}$.
\end{lemma}

This result will be proved in the next section.

Now we assume that $V$ satisfies either one of the properties in
Proposition 14.1.
By that proposition, $V$ satisfies both properties. For any
$P(z_{1})$-intertwining map $F_{1}$ and $P(z_{2})$-intertwining map
 $F_{2}$ as in Proposition 14.1,
using the notation in the theory of (partial) operads (see \cite{HL1}
and \cite{HL2}),
 we denote the product of $F_{1}$ and $F_{2}$
by $\gamma(F_{1}; I, F_{2})$. The reader should note the
well-definedness of the expressions
in the  calculations below.
Using the definitions of  $\gamma(F_{1}; I, F_{2})$ and
of $P(z_{1})$-intertwining map,  we have
\begin{eqnarray}
\lefteqn{\langle w'_{(4)}, x_{1}^{-1}\delta\left(\frac{ x_{0}-z_{1}}
{x_{1}}\right)x_{2}^{-1}\delta\left(\frac{ x_{0}-z_{2}}{x_{2}}\right)
Y_{4}(v, x_{0})\cdot}\nno\\
&&\hspace{4em}\cdot
\gamma(F_{1}; I, F_{2})(w_{(1)}\otimes w_{(2)}\otimes w_{(3)})
\rangle_{W_{4}}\nno\\
&&=\sum_{m\in {\Bbb C}}
\langle w'_{(4)}, x_{1}^{-1}\delta\left(\frac{ x_{0}-z_{1}}
{x_{1}}\right)x_{2}^{-1}\delta\left(\frac{
x_{0}-z_{2}}{x_{2}}\right)\cdot\nno\\
&&\hspace{4em}\cdot Y_{4}(v, x_{0})P_{m}(\gamma(F_{1}; I, F_{2})(w_{(1)}\otimes
w_{(2)}\otimes w_{(3)}))
\rangle_{W_{4}}\nno\\
&&=\sum_{m\in {\Bbb C}}\sum_{n\in {\Bbb C}}
\langle w'_{(4)}, x_{1}^{-1}\delta\left(\frac{ x_{0}-z_{1}}
{x_{1}}\right)x_{2}^{-1}\delta\left(\frac{
x_{0}-z_{2}}{x_{2}}\right)\cdot\nno\\
&&\hspace{4em}\cdot Y_{4}(v, x_{0})P_{m}(F_{1}(w_{(1)}\otimes P_{n}(F_{2}
(w_{(2)}\otimes w_{(3)}))))
\rangle_{W_{4}}\nno\\
&&=\sum_{m\in {\Bbb C}}\sum_{n\in {\Bbb C}}
\langle w'_{(4)},
x_{2}^{-1}\delta\left(\frac{x_{0}-z_{2}}{x_{2}}\right)
P_{m}\biggl(x_{1}^{-1}\delta\left(\frac{ x_{0}-z_{1}}{x_{1}}\right)\cdot\nno\\
&&\hspace{4em}\cdot
Y_{4}(v, x_{0})F_{1}(w_{(1)}\otimes P_{n}(F_{2}(w_{(2)}\otimes
w_{(3)})))\biggr)
\rangle_{W_{4}}\nno\\
&&=\sum_{m\in {\Bbb C}}\sum_{n\in {\Bbb C}}
\langle w'_{(4)},
x_{2}^{-1}\delta\left(\frac{x_{0}-z_{2}}{x_{2}}\right)
P_{m}\biggl(z_{1}^{-1}\delta\left(\frac{x_{0}-x_{1}}{z_{1}}\right)\cdot\nno\\
&&\hspace{4em}\cdot
F_{1}(Y_{1}(v, x_{1})w_{(1)}\otimes P_{n}(F_{2}(w_{(2)}\otimes
w_{(3)})))\biggr)
\rangle_{W_{4}}\nno\\
&&\quad +\sum_{m\in {\Bbb C}}\sum_{n\in {\Bbb C}}
\langle w'_{(4)},
x_{2}^{-1}\delta\left(\frac{x_{0}-z_{2}}{x_{2}}\right)
P_{m}\biggl(x_{1}^{-1}\delta\left(\frac{z_{1}-x_{0}}{-x_{1}}\right)\cdot\nno\\
&&\hspace{4em}\cdot
F_{1}(w_{(1)}\otimes Y_{5}(v, x_{0})P_{n}(F_{2}(w_{(2)}\otimes
w_{(3)})))\biggr)
\rangle_{W_{4}}\nno\\
&&=\langle w'_{(4)},
x_{2}^{-1}\delta\left(\frac{x_{0}-z_{2}}{x_{2}}\right)
z_{1}^{-1}\delta\left(\frac{x_{0}-x_{1}}{z_{1}}\right)\cdot\nno\\
&&\hspace{4em}\cdot
\gamma(F_{1}; I, F_{2})(Y_{1}(v, x_{1})w_{(1)}\otimes w_{(2)}\otimes w_{(3)})
\rangle_{W_{4}}\nno\\
&&\quad +\sum_{m\in {\Bbb C}}\sum_{n\in {\Bbb C}}
\langle w'_{(4)},
x_{2}^{-1}\delta\left(\frac{x_{0}-z_{2}}{x_{2}}\right)
P_{m}\biggl(x_{1}^{-1}\delta\left(\frac{z_{1}-x_{0}}{-x_{1}}\right)\cdot\nno\\
&&\hspace{4em}\cdot
F_{1}(w_{(1)}\otimes Y_{5}(v, x_{0})P_{n}(F_{2}(w_{(2)}\otimes
w_{(3)})))\biggr)
\rangle_{W_{4}}\hspace{3em}
\end{eqnarray}
for all $v\in V$, $w_{(1)}\in W_{1}$,
$w_{(2)}\in W_{2}$,
$w_{(3)}\in W_{3}$ and $w'_{(4)}\in W'_{4}$.
By the definition of $P(z_{2})$-intertwining map, the second term in the
right-hand side of  (14.12)
becomes
\begin{eqnarray}
\lefteqn{\sum_{m\in {\Bbb C}}\sum_{n\in {\Bbb C}}
\langle w'_{(4)},
P_{m}\biggl(x_{1}^{-1}\delta\left(\frac{z_{1}-x_{0}}{-x_{1}}\right)
F_{1}\biggl(w_{(1)}\otimes} \nno\\
&&\hspace{4em}\otimes
P_{n}\biggl(x_{2}^{-1}\delta\left(\frac{x_{0}-z_{2}}{x_{2}}\right)
Y_{5}(v, x_{0})F_{2}(w_{(2)}\otimes w_{(3)})\biggr)\biggr)\biggr)
\rangle_{W_{4}}\nno\\
&&=\sum_{m\in {\Bbb C}}\sum_{n\in {\Bbb C}}
\langle w'_{(4)},
P_{m}\biggl(x_{1}^{-1}\delta\left(\frac{z_{1}-x_{0}}{-x_{1}}\right)
F_{1}\biggl(w_{(1)}\otimes\nno\\
&&\hspace{4em}\otimes
P_{n}\biggl(z_{2}^{-1}\delta\left(\frac{x_{0}-x_{2}}{z_{2}}\right)
F_{2}(Y_{2}(v, x_{2})w_{(2)}\otimes w_{(3)})\biggr)\biggr)\biggr)
\rangle_{W_{4}}\nno\\
&& \quad +\sum_{m\in {\Bbb C}}\sum_{n\in {\Bbb C}}
\langle w'_{(4)},
P_{m}\biggl(x_{1}^{-1}\delta\left(\frac{z_{1}-x_{0}}{-x_{1}}\right)
F_{1}\biggl(w_{(1)}\otimes \nno\\
&&\hspace{4em}\otimes
P_{n}\biggl(x_{2}^{-1}\delta\left(\frac{z_{2}-x_{0}}{-x_{2}}\right)
F_{2}(w_{(2)}\otimes Y_{3}(v, x_{0})w_{(3)})\biggr)\biggr)\biggr)
\rangle_{W_{4}}\nno\\
&&=\langle w'_{(4)},
x_{1}^{-1}\delta\left(\frac{z_{1}-x_{0}}{-x_{1}}\right)
z_{2}^{-1}\delta\left(\frac{x_{0}-x_{2}}{z_{2}}\right)\cdot\nno\\
&&\hspace{4em}\cdot
\gamma(F_{1}; I, F_{2})(w_{(1)}\otimes
Y_{2}(v, x_{2})w_{(2)}\otimes w_{(3)})
\rangle_{W_{4}}\nno\\
&& \quad +\langle w'_{(4)},
x_{1}^{-1}\delta\left(\frac{z_{1}-x_{0}}{-x_{1}}\right)
x_{2}^{-1}\delta\left(\frac{z_{2}-x_{0}}{-x_{2}}\right)\cdot\nno\\
&&\hspace{4em}\cdot
\gamma(F_{1}; I, F_{2})(w_{(1)}\otimes
w_{(2)}\otimes Y_{3}(v, x_{0})w_{(3)})
\rangle_{W_{4}}.
\end{eqnarray}
Substituting (14.13) into (14.12), we obtain
\begin{eqnarray}
\lefteqn{\langle w'_{(4)}, x_{1}^{-1}\delta\left(\frac{ x_{0}-z_{1}}
{x_{1}}\right)x_{2}^{-1}\delta\left(\frac{ x_{0}-z_{2}}{x_{2}}\right)
Y_{4}(v, x_{0})\cdot}\nno\\
&&\hspace{4em}\cdot
\gamma(F_{1}; I, F_{2})(w_{(1)}\otimes w_{(2)}\otimes w_{(3)})
\rangle_{W_{4}}\nno\\
&&=\langle w'_{(4)},
x_{2}^{-1}\delta\left(\frac{x_{0}-z_{2}}{x_{2}}\right)
z_{1}^{-1}\delta\left(\frac{x_{0}-x_{1}}{z_{1}}\right)\cdot\nno\\
&&\hspace{4em}\cdot
\gamma(F_{1}; I, F_{2})(Y_{1}(v, x_{1})w_{(1)}\otimes w_{(2)}\otimes w_{(3)}))
\rangle_{W_{4}}\nno\\
&&\quad +\langle w'_{(4)},
x_{1}^{-1}\delta\left(\frac{z_{1}-x_{0}}{-x_{1}}\right)
z_{2}^{-1}\delta\left(\frac{x_{0}-x_{2}}{z_{2}}\right)\cdot\nno\\
&&\hspace{4em}\cdot
\gamma(F_{1}; I, F_{2})(w_{(1)}\otimes
Y_{2}(v, x_{2})w_{(2)}\otimes w_{(3)}))
\rangle_{W_{4}}\nno\\
&& \quad +\langle w'_{(4)},
x_{1}^{-1}\delta\left(\frac{z_{1}-x_{0}}{-x_{1}}\right)
x_{2}^{-1}\delta\left(\frac{z_{2}-x_{0}}{-x_{2}}\right)\cdot\nno\\
&&\hspace{4em}\cdot
\gamma(F_{1}; I, F_{2})(w_{(1)}\otimes
w_{(2)}\otimes Y_{3}(v, x_{0})w_{(3)}))
\rangle_{W_{4}}. \hspace{3em}
\end{eqnarray}
Since (14.14) holds for any $w'_{(4)}\in W'_{4}$, we obtain
\begin{eqnarray}
\lefteqn{x_{1}^{-1}\delta\left(\frac{ x_{0}-z_{1}}
{x_{1}}\right)x_{2}^{-1}\delta\left(\frac{ x_{0}-z_{2}}{x_{2}}\right)\cdot}
\nno\\
&&\cdot
(Y_{4}(v, x_{0})\gamma(F_{1}; I, F_{2}))(w_{(1)}\otimes w_{(2)}\otimes w_{(3)})
\nno\\
&&=
x_{2}^{-1}\delta\left(\frac{x_{0}-z_{2}}{x_{2}}\right)
z_{1}^{-1}\delta\left(\frac{x_{0}-x_{1}}{z_{1}}\right)\cdot\nno\\
&&\hspace{4em}\cdot
\gamma(F_{1}; I, F_{2})(Y_{1}(v, x_{1})w_{(1)}\otimes w_{(2)}\otimes
w_{(3)}))\nno\\
&&\quad +
x_{1}^{-1}\delta\left(\frac{z_{1}-x_{0}}{-x_{1}}\right)
z_{2}^{-1}\delta\left(\frac{x_{0}-x_{2}}{z_{2}}\right)\cdot\nno\\
&&\hspace{4em}\cdot \gamma(F_{1}; I, F_{2})(w_{(1)}\otimes
Y_{2}(v, x_{2})w_{(2)}\otimes w_{(3)}))\nno\\
&& \quad +
x_{1}^{-1}\delta\left(\frac{z_{1}-x_{0}}{-x_{1}}\right)
x_{2}^{-1}\delta\left(\frac{z_{2}-x_{0}}{-x_{2}}\right)\cdot\nno\\
&&\hspace{4em}\cdot
\gamma(F_{1}; I, F_{2})(w_{(1)}\otimes
w_{(2)}\otimes Y_{3}(v, x_{0})w_{(3)})).
\end{eqnarray}

We are mainly interested in (14.14). The left-hand side of (14.14) can be
written as
\begin{eqnarray}
\lefteqn{\langle x_{1}^{-1}\delta\left(\frac{ x_{0}-z_{1}}
{x_{1}}\right)x_{2}^{-1}\delta\left(\frac{ x_{0}-z_{2}}{x_{2}}\right)
\cdot}\nno\\
&&\cdot
Y'_{4}(e^{x_{0}L(1)}(-x_{0}^{2})^{-L(0)}v, x_{0}^{-1})w'_{(4)},
\gamma(F_{1}; I, F_{2})(w_{(1)}\otimes w_{(2)}\otimes w_{(3)})
\rangle_{W_{4}}.\nno\\
&&
\end{eqnarray}
First replacing $v$ by $(-x_{0}^{2})^{L(0)}e^{-x_{0}L(1)}v$ and then
replacing $x_{0}$ by $x_{0}^{-1}$ in both (14.16) and the right-hand side of
(14.14), we obtain
\begin{eqnarray}
\lefteqn{\langle x_{1}^{-1}\delta\left(\frac{ x^{-1}_{0}-z_{1}}
{x_{1}}\right)x_{2}^{-1}\delta\left(\frac{x^{-1}_{0}-z_{2}}{x_{2}}\right)
\cdot}\nno\\
&&\hspace{3em}\cdot Y'_{4}(v, x_{0})w'_{(4)},
\gamma(F_{1}; I, F_{2})(w_{(1)}\otimes w_{(2)}\otimes w_{(3)})
\rangle_{W_{4}}\nno\\
&&=\langle w'_{(4)},
x_{2}^{-1}\delta\left(\frac{x^{-1}_{0}-z_{2}}{x_{2}}\right)
z_{1}^{-1}\delta\left(\frac{x^{-1}_{0}-x_{1}}{z_{1}}\right)\cdot\nno\\
&&\hspace{3em}\cdot
\gamma(F_{1}; I, F_{2})(Y_{1}((-x_{0}^{-2})^{L(0)}e^{-x^{-1}_{0}L(1)}v,
x_{1})w_{(1)}\otimes w_{(2)}\otimes w_{(3)}))
\rangle_{W_{4}}\nno\\
&&\quad +\langle w'_{(4)},
x_{1}^{-1}\delta\left(\frac{z_{1}-x^{-1}_{0}}{-x_{1}}\right)
z_{2}^{-1}\delta\left(\frac{x^{-1}_{0}-x_{2}}{z_{2}}\right)\cdot\nno\\
&&\hspace{3em}\cdot
\gamma(F_{1}; I, F_{2})(w_{(1)}\otimes
Y_{2}((-x_{0}^{-2})^{L(0)}e^{-x^{-1}_{0}L(1)}v, x_{2})w_{(2)}\otimes w_{(3)}))
\rangle_{W_{4}}\nno\\
&& \quad +\langle w'_{(4)},
x_{1}^{-1}\delta\left(\frac{z_{1}-x^{-1}_{0}}{-x_{1}}\right)
x_{2}^{-1}\delta\left(\frac{z_{2}-x^{-1}_{0}}{-x_{2}}\right)\cdot\nno\\
&&\hspace{3em}\cdot
\gamma(F_{1}; I, F_{2})(w_{(1)}\otimes
w_{(2)}\otimes Y_{3}((-x_{0}^{-2})^{L(0)}e^{-x^{-1}_{0}L(1)}v, x^{-1}_{0})
w_{(3)}))
\rangle_{W_{4}}. \nno\\
&&
\end{eqnarray}
Note that the map $\gamma(F_{1}; I, F_{2}): W_{1}\otimes W_{2}\otimes W_{3}
\to \overline{W}_{4}$ amounts to a map {}from $W'_{4}$ to
$(W_{1}\otimes W_{2}\otimes
W_{3})^{*}$ and that the coefficients of
$$x_{1}^{-1}\delta\left(\frac{
x^{-1}_{0}-z_{1}}
{x_{1}}\right)x_{2}^{-1}\delta\left(\frac{x^{-1}_{0}-z_{2}}{x_{2}}\right)
Y_{t}(v, x_{0})$$
in powers of $x_{0}, x_{1}$ and $x_{2}$, for all $v\in V$, span
$$V \otimes \iota_{+}{\Bbb C}[t,t^{-
1},
(z_{1}^{-1}-t)^{-1}, (z_{2}^{-1}-t)^{-1}].$$
Thus the formula (14.17) motivates the definition of the following
action $\tau^{(1)}_{P(z_{1}, z_{2})}$ of $$V \otimes \iota_{+}{\Bbb C}[t,t^{-
1},
(z_{1}^{-1}-t)^{-1}, (z_{2}^{-1}-t)^{-1}]$$ on $(W_{1}\otimes W_{2}\otimes
W_{3})^{*}$:
\begin{eqnarray}
\lefteqn{(\tau^{(1)}_{P(z_{1}, z_{2})}(x_{1}^{-1}\delta\left(\frac{
x^{-1}_{0}-z_{1}}
{x_{1}}\right)x_{2}^{-1}\delta\left(\frac{x^{-1}_{0}-z_{2}}{x_{2}}\right)
Y_{t}(v, x_{0}))\lambda)(w_{(1)}\otimes w_{(2)}\otimes w_{(3)})}\nno\\
&&=
x_{2}^{-1}\delta\left(\frac{x^{-1}_{0}-z_{2}}{x_{2}}\right)
z_{1}^{-1}\delta\left(\frac{x^{-1}_{0}-x_{1}}{z_{1}}\right)\cdot\nno\\
&&\hspace{4em}\cdot
\lambda(Y_{1}((-x_{0}^{-2})^{L(0)}e^{-x^{-1}_{0}L(1)}v,
x_{1})w_{(1)}\otimes w_{(2)}\otimes w_{(3)}))\nno\\
&&\quad +
x_{1}^{-1}\delta\left(\frac{z_{1}-x^{-1}_{0}}{-x_{1}}\right)
z_{2}^{-1}\delta\left(\frac{x^{-1}_{0}-x_{2}}{z_{2}}\right)\cdot\nno\\
&&\hspace{4em}\cdot
\lambda(w_{(1)}\otimes
Y_{2}((-x_{0}^{-2})^{L(0)}e^{-x^{-1}_{0}L(1)}v,
x_{2})w_{(2)}\otimes w_{(3)}))\nno\\
&& \quad +
x_{1}^{-1}\delta\left(\frac{z_{1}-x^{-1}_{0}}{-x_{1}}\right)
x_{2}^{-1}\delta\left(\frac{z_{2}-x^{-1}_{0}}{-x_{2}}\right)\cdot\nno\\
&&\hspace{4em}\cdot
\lambda(w_{(1)}\otimes
w_{(2)}\otimes Y_{3}((-x_{0}^{-2})^{L(0)}e^{-x^{-1}_{0}L(1)}v,
x^{-1}_{0})w_{(3)}))
\hspace{4em}
\end{eqnarray}
for all $\lambda\in (W_{1}\otimes W_{2}\otimes
W_{3})^{*}$ and $w_{(1)}\in W_{1}$, $w_{(2)}\in W_{2}$ and $w_{(3)}\in W_{3}$.

Similarly, for any
$P(z_{3})$-intertwining map $F_{3}$ and $P(z_{4})$-intertwining map
 $F_{4}$ as in Proposition 14.1, we denote the iterate of $F_{3}$ and $F_{4}$
by $\gamma(F_{4};  F_{3}, I)$.  Similarly to the discussion for
$\gamma(F_{1}; I, F_{2})$, we can show that
\begin{eqnarray}
\lefteqn{\langle
x^{-1}_2\delta\left(\frac{x^{-1}_0-z_4}{x_2}\right)
x^{-1}_1\delta\left(\frac{x_2-z_3}{x_1}\right)\cdot}\nno\\
&&\hspace{2em}\cdot Y'_{4}(v, x_0)w'_{(4)},
\gamma(F_{4};  F_{3}, I)
(w_{(1)}\otimes w_{(2)}\otimes w_{(3)})\rangle_{W_{4}}\nno\\
&=&\langle w'_{(4)}, x^{-1}_2\delta\left(\frac{z_4-x^{-1}_0}{x_2}\right)
x^{-1}_1\delta\left(\frac{x_2-z_3}{x_1}\right)\cdot\nno\\
&&\hspace{2em}\cdot
\gamma(F_{4};  F_{3}, I)(w_{(1)}\otimes w_{(2)}\otimes
Y((-x_{0}^{-2})^{L(0)}e^{-x^{-1}_{0}L(1)}v, x^{-1}_0)
w_{(3)})\rangle_{W_{(4)}}\nno\\
&& + \langle w'_{(4)},
z^{-1}_4\delta\left(\frac{x^{-1}_0-x_2}{z_4}\right) z_3^{-1}
\delta\left(\frac{x_2-x_1}{z_3}\right)\cdot\nno\\
&&\hspace{2em}\cdot
\gamma(F_{4};  F_{3}, I)(Y((-x_{0}^{-2})^{L(0)}e^{-x^{-1}_{0}L(1)}
v, x_1) w_{(1)}\otimes w_{(2)}\otimes w_{(3)})\rangle_{W_{(4)}}\nno\\
&& + \langle w'_{(4)}, z^{-1}_4\delta\left(\frac{x^{-1}_0-x_2}{z_4}\right)
x^{-1}_1
\delta\left(\frac{z_3-x_2}{x_1}\right)\cdot\nno\\
&&\hspace{2em}\cdot
\gamma(F_{4};  F_{3}, I)(w_{(1)}\otimes Y((-x_{0}^{-2})^{L(0)}
e^{-x^{-1}_{0}L(1)}v, x_2) w_{(2)}\otimes w_{(3)})\rangle_{W_{(4)}}\nno\\
&&
\end{eqnarray}
for any $w_{(1)}\in W_{1}$, $w_{(2)}\in W_{2}$, $w_{(3)}\in W_{3}$ and
$w'_{(4)}\in W'_{4}$. Note that all terms in both sides of (14.19)
 are well defined.
Since the components of
$$x^{-1}_2\delta\left(\frac{x^{-1}_0-z_4}{x_2}\right)
 x^{-1}_1\delta\left(\frac{x_2-z_3}{x_1}\right)
Y_{t}(v, x_0)$$
in powers of $x_{0}, x_{1}$ and $x_{2}$ span
$$V \otimes \iota_{+}{\Bbb C}[t,t^{- 1},
((z_{3}+z_{4})^{-1}-t)^{-1}, (z_{4}^{-1}-t)^{-1}],$$
(14.19) motivates the following
definition of an
action $\tau^{(2)}_{P(z_{3}+z_{4}, z_{4})}$
of $$V \otimes \iota_{+}{\Bbb C}[t,t^{- 1},
((z_{3}+z_{4})^{-1}-t)^{-1}, (z_{4}^{-1}-t)^{-1}]$$ on
$\lambda\in (W_{1}\otimes W_{2}\otimes
W_{3})^{*}$:
\begin{eqnarray}
\lefteqn{\biggl(\tau^{(2)}_{P(z_{3}+z_{4}, z_{4})}
\biggl(x^{-1}_2\delta\left(\frac{x^{-1}_0-z_4}{x_2}\right)\cdot}\nno\\
&&\hspace{4em}\cdot x^{-1}_1\delta\left(\frac{x_2-z_3}{x_1}\right)
Y_{t}(v, x_0)\biggr)\lambda
\biggr)
(w_{(1)}\otimes w_{(2)}\otimes w_{(3)})\nno\\
&&=x^{-1}_2\delta\left(\frac{z_4-x^{-1}_0}{x_2}\right)
x^{-1}_1\delta\left(\frac{x_2-z_3}{x_1}\right)\cdot\nno\\
&&\hspace{4em}\cdot
\lambda(w_{(1)}\otimes w_{(2)}\otimes
Y((-x_{0}^{-2})^{L(0)}
e^{-x^{-1}_{0}L(1)}v, x^{-1}_0) w_{(3)})\nno\\
&&\quad +
z^{-1}_4\delta\left(\frac{x^{-1}_0-x_2}{z_4}\right) z_3^{-1}
\delta\left(\frac{x_2-x_1}{z_3}\right)\cdot\nno\\
&&\hspace{4em}\cdot
\lambda(Y((-x_{0}^{-2})^{L(0)}
e^{-x^{-1}_{0}L(1)}v, x_1) w_{(1)}\otimes w_{(2)}\otimes w_{(3)})\nno\\
&&\quad + z^{-1}_4\delta\left(\frac{x^{-1}_0-x_2}{z_4}\right)x^{-1}_1
\delta\left(\frac{z_3-x_2}{x_1}\right)\cdot\nno\\
&&\hspace{4em}\cdot
\lambda(w_{(1)}\otimes Y((-x_{0}^{-2})^{L(0)}
e^{-x^{-1}_{0}L(1)}v, x_2) w_{(2)}\otimes w_{(3)})\hspace{4em}
\end{eqnarray}
for all $\lambda\in (W_{1}\otimes W_{2}\otimes
W_{3})^{*}$ and $w_{(1)}\in W_{1}$, $w_{(2)}\in W_{2}$ and $w_{(3)}\in W_{3}$.

{}From (14.8), (14.9), (14.10) and (14.11), we see that
\begin{equation}
\tau^{(1)}_{P(z_1,z_2)}=\tau^{(2)}_{P(z_1, z_2)}
\end{equation}
when (14.7) holds.

When (14.7) holds, let
\begin{equation}
L'_{P(z_{1}, z_{2})}(0)=\tau^{(1)}_{P(z_{1}, z_{2})}(\omega\otimes t^{0})
=\tau^{(2)}_{P(z_{1}, z_{2})}(\omega\otimes t^{0})
\end{equation}
and
\begin{equation}
Y'_{P(z_{1}, z_{2})}(v, x)=\tau^{(1)}_{P(z_{1}, z_{2})}(Y_{t}(v, x))
=\tau^{(2)}_{P(z_{1}, z_{2})}(Y_{t}(v, x)).
\end{equation}
We call the eigenspaces of the operator $L'_{P(z_{1}, z_{2})}(0)$ the
{\it $P(z_{1}, z_{2})$-weight
subspaces} (or {\it
$P(z_{1}, z_{2})$-homogeneous subspaces}) of
$(W_{1}\otimes W_{2}\otimes W_{3})^{*}$, and we have the
corresponding notions of {\it $P(z_{1}, z_{2})$-weight vector} (or {\it
$P(z_{1}, z_{2})$-homogen\-eous
vector}) and  {\it $P(z_{1}, z_{2})$-weight}.

Let
$W_{1}$, $W_{2}$, $W_{3}$, $W_{4}$ and $W_{5}$ be  $V$-modules,
$F_{1}$, $F_{2}$, $F_{3}$
and $F_{4}$
$P(z_{1})$-, $P(z_{2})$-, $P(z_{1}-z_{2})$- and $P(z_{2})$-intertwining maps of
type ${W_{4}}\choose
{W_{1}W_{5}}$, ${W_{5}}\choose {W_{2}W_{3}}$, ${W_{5}}\choose
{W_{1}W_{2}}$ and ${W_{4}}\choose
{W_{5}W_{3}}$, respectively.
Then the restrictions to $W'_{4}$ of the adjoint maps of the product
$\gamma(F_{1}; I, F_{2})$ of
$F_{1}$ and $F_{2}$ when $|z_{1}|>|z_{2}|>0$ and
of the iterate $\gamma(F_{4};  F_{3}, I)$
of $F_{3}$ and $F_{4}$ when $|z_{2}|>|z_{1}-z_{2}|>0$
are maps $\gamma(F_{1}; I, F_{2})'$ and $\gamma(F_{4};  F_{3}, I)'$,
respectively,
{}from $W'_{4}$ to
$(W_{1}\otimes W_{2}\otimes W_{3})^{*}$ defined by
\begin{eqnarray}
\lefteqn{(\gamma(F_{1}; I, F_{2})'(w'_{(4)}))(w_{(1)}\otimes
w_{(2)}\otimes w_{(3)})=}\nno\\
&&=\langle w'_{(4)}, {\cal Y}_{1}(w_{(1)}, x_{1})
{\cal Y}_{2}(w_{(2)}, x_{2})w_{(3)}\rangle_{W_{4}}\lbar_{x_{1}
=z_{1}, \;x_{2}=z_{2}}
\end{eqnarray}
 and
\begin{eqnarray}
\lefteqn{(\gamma(F_{4};  F_{3}, I)'(w'_{(4)}))(w_{(1)}\otimes
w_{(2)}\otimes w_{(3)})=}\nno\\
&&=\langle w'_{(4)}, {\cal Y}_{4}(
{\cal Y}_{3}(w_{(1)}, x_{1})w_{(2)}, x_{2})w_{(3)})\rangle_{W_{4}}
\lbar_{x_{1}
=z_{3}, \;x_{2}=z_{4}},
\end{eqnarray}
 respectively. It is easy to verify that
for any $w'_{(4)}\in W'_{4}$, $\gamma(F_{1}; I, F_{2})'(w'_{(4)})$ and
$\gamma(F_{3};  F_{4}, I)'(w'_{(4)})$ satisfy the
following conditions for $\lambda \in (W_{1}\otimes W_{2}\otimes W_{3})^{*}$
when (14.7) holds:

\begin{description}
\item[The $P(z_{1}, z_{2})$-compatibility condition] \hfill

{\bf (a)} The  {\it $P(z_{1}, z_{2})$-lower
truncation condition}:
For all $v\in V$, the formal Laurent series $Y'_{P(z_{1}, z_{2})}(v, x)\lambda$
involves only finitely many negative
powers of $x$.

{\bf (b)} The following formula holds:
\begin{eqnarray}
\lefteqn{\tau^{(1)}_{P(z_1, z_2)}\left(x^{-1}_0\delta\left(\frac{x^{-1}_0-
z_1}{x_1}\right)x^{-1}_2\delta\left(\frac{x^{-1}_0-z_2}{x_2}\right)Y_{t}(v,
x_0)\right)\lambda=} \nno\\
&&=x^{-1}_1\delta\left(\frac{x^{-1}_0-
z_1}{x_1}\right)x^{-1}_2\delta\left(\frac{x^{-1}_0-
z_2}{x_2}\right) \tau^{(1)}_{P(z_{1}, z_{2})}(Y_{t}(v, x_{0}))\lambda\nno\\
&&\hspace{15em}\mbox{for all}\;\;\;v\in V,
\end{eqnarray}
or equivalently,
\begin{eqnarray}
\lefteqn{\tau^{(2)}_{P(z_1, z_2)}\left(x^{-1}_0\delta
\left(\frac{x^{-1}_0-z_1}{x_1}\right)
x^{-1}_2\delta\left(\frac{x^{-1}_0-z_2}{x_2}\right)Y_{t}(v,
x_0)\right)
\lambda=}
\nno\\
&&=x^{-1}_1\delta\left(\frac{x^{-1}_0-z_1}{x_1}\right)x^{-1}_2\delta
\left(\frac{x^{-1}_0-z_2}{x_2}\right) \tau^{(2)}_{P(z_1, z_2)}(Y_{t}(v, x_0))
\lambda\nno\\
&&\hspace{15em}\mbox{for all}\;\;\;v\in V.
\end{eqnarray}

\item[The $P(z_{1}, z_{2})$-local grading-restriction  condition]\hfill

{\bf (a)} The {\it $P(z_{1}, z_{2})$-grading condition}:
$\lambda$ is a (finite) sum of
eigenvectors of $L'_{P(z_{1}, z_{2})}(0)$.

{\bf (b)} Let $W_{\lambda}$ be the smallest subspace of $(W_{1}\otimes
W_{2}\otimes W_{3})^{*}$ containing $\lambda$ and stable under the component
operators $\tau^{(1)}_{P(z_{1}, z_{2})}(v\otimes t^{n})$ of
the operators $Y'_{P(z_{1}, z_{2})}(v,
x)$ for $v\in V$, $n\in {\Bbb Z}$. Then the weight spaces
$(W_{\lambda})_{(n)}$, $n\in {\Bbb C}$, of the (graded) space
$W_{\lambda}$ have the properties
\begin{eqnarray}
&\mbox{\rm dim}\ (W_{\lambda})_{(n)}<\infty \;\;\;\mbox{\rm for}\
n\in {\Bbb C},&\\
&(W_{\lambda})_{(n)}=0 \;\;\;\mbox{\rm for $n$ whose real part is
sufficiently small.}&
\end{eqnarray}
\end{description}

These two conditions are the analogues of the corresponding conditions for
$\tau_{P(z)}$ and $\tau_{Q(z)}$. But in this case,
we also have to consider
 two additional conditions. To state the conditions,
it is convenient to introduce the following concepts:
Let $W$ be a vector space and
$W^{*}$ its dual. A formal sum $\sum_{n\in {\Bbb C}}w^{*}_{n}$
where $w^{*}_{n}\in W^{*}$, $n\in {\Bbb C}$, is called a {\it series}
in $W^{*}$. A series $\sum_{n\in {\Bbb C}}w^{*}_{n}$ in $W^{*}$ is said to
be {\it absolutely convergent} if for any $w\in W$, the sum
$\sum_{n\in {\Bbb C}}|w_{n}^{*}(w)|$ is convergent. If $\sum_{n\in
{\Bbb C}}w^{*}_{n}$ is absolutely convergent, then for any $w\in W$,
the sum $\sum_{n\in {\Bbb C}}w^{*}_{n}(w)$ is also convergent and
the limits for all $w\in W$ define an element of $W^{*}$. This element
is called the {\it limit} of the series $\sum_{n\in {\Bbb C}}w^{*}_{n}$.
If a series $\sum_{n\in {\Bbb C}}w_{n}^{*}$
of homogeneous vectors in a ${\Bbb C}$-graded subspace of $W^{*}$
have the property that for any $w\in W$,
the weights of all those vectors $w_{n}^{*}$ in the series such that
$w_{n}^{*}(w)\ne 0$ can be arranged to give
a sequence of strictly increasing real numbers, we say that this series
is a
{\it series indexed by sequences of strictly increasing real numbers}.
We consider the following two conditions for
$\lambda\in (W_{1}\otimes W_{2}\otimes W_{3})^{*}$:

\begin{description}
\item[The $P(z_{2})$-local grading restriction condition]\hfill

{\bf  (a)} The  {$P(z_{2})$-grading condition}: For any $w_{(1)}\in W_{1}$, the
element
$\mu^{(1)}_{\lambda, w_{(1)}}\in (W_{2}\otimes W_{3})^{*}$ defined by
\begin{equation}
\mu^{(1)}_{\lambda, w_{(1)}}(w_{(2)}\otimes w_{(3)})
=\lambda(w_{(1)}\otimes w_{(2)}\otimes w_{(3)})
\end{equation}
for all $w_{(2)}\in W_{2}$ and $w_{(3)}\in W_{3}$, is the  limit of  an
absolutely convergent series  of $P(z_{2})$-weight vectors
in $(W_{2}\otimes
W_{3})^{*}$ indexed by
sequences of strictly increasing real numbers.

{\bf (b)} For any $w_{(1)}\in W_{1}$, let $W^{(1)}_{\lambda, w_{(1)}}$ be
the smallest subspace
of  $(W_{2}\otimes W_{3})^{*}$ containing the $P(z_{2})$-weight  vectors in the
series
absolutely convergent to $\mu^{(1)}_{\lambda, w_{(1)}}$ and stable under the
component operators
$\tau_{P(z_{2})}(v\otimes t^{n})$ of
the operators $Y'_{P(z_{2})}(v,
x)$ for $v\in V$, $n\in {\Bbb Z}$. Then the weight spaces
$(W^{(1)}_{\lambda, w_{(1)}})_{(n)}$, $n\in {\Bbb C}$, of the (graded) space
$W^{(1)}_{\lambda, w_{(1)}}$ have the properties
\begin{eqnarray}
&\mbox{\rm dim}\ (W^{(1)}_{\lambda, w_{(1)}})_{(n)}<\infty \;\;\;\mbox{\rm
for}\
n\in {\Bbb C},&\\
&(W^{(1)}_{\lambda, w_{(1)}})_{(n)}=0 \;\;\;\mbox{\rm for $n$ whose real part
is
sufficiently small.}\;\;\;\;\;\;&
\end{eqnarray}

\item[The $P(z_{1}-z_{2})$-local grading restriction condition]\hfill

{\bf  (a)} The  {$P(z_{1}-z_{2})$-grading condition}: For any $w_{(3)}\in
W_{3}$, the element
$\mu^{(2)}_{\lambda,  w_{(3)}}\in (W_{1}\otimes W_{2})^{*}$ defined by
\begin{equation}
\mu^{(2)}_{\lambda, w_{(3)}}(w_{(1)}\otimes w_{(2)})
=\lambda(w_{(1)}\otimes w_{(2)}\otimes w_{(3)})
\end{equation}
for all $w_{(1)}\in W_{1}$ and $w_{(2)}\in W_{2}$, is the  limit of  an
absolutely convergent  series
of  $P(z_{1}-z_{2})$-weight vectors in
$(W_{1}\otimes
W_{2})^{*}$ indexed by
sequences of strictly increasing real numbers.

{\bf (b)} For any $w_{(3)}\in W_{3}$, let $W^{(2)}_{\lambda, w_{(3)}}$ be
the smallest subspace
of  $(W_{1}\otimes W_{2})^{*}$ containing the $P(z_{1}-z_{2})$-weight  vectors
in the series
absolutely convergent to $\mu^{(2)}_{\lambda, w_{(3)}}$ and stable under the
component operators
$$\tau_{P(z_{1}-z_{2})}(v\otimes t^{n})$$ of
the operators $Y'_{P(z_{1}-z_{2})}(v,
x)$ for $v\in V$, $n\in {\Bbb Z}$. Then the weight spaces
$(W^{(2)}_{\lambda, w_{(3)}})_{(n)}$, $n\in {\Bbb C}$, of the (graded) space
$W^{(2)}_{\lambda, w_{(3)}}$ have the properties
\begin{eqnarray}
&\mbox{\rm dim}\ (W^{(2)}_{\lambda, w_{(3)}})_{(n)}<\infty \;\;\;\mbox{\rm
for}\
n\in {\Bbb C},&\\
&(W^{(2)}_{\lambda, w_{(3)}})_{(n)}=0 \;\;\;\mbox{\rm for $n$ whose real part
is
sufficiently small.}\;\;\;\;\;\;&
\end{eqnarray}
\end{description}

We have:

\begin{lemma}
Let $z_{1}, z_{2}$ be complex numbers satisfying (14.7).
Assume that $\lambda\in (W_{1}\otimes W_{2}\otimes W_{3})^{*}$ satisfies the
$P(z_{1}, z_{2})$-compatibil\-ity condition. Then for any $v\in V$
and $w_{(3)}\in W_{3}$,
$$x^{-1}_1 \delta\left(\frac{x^{-1}-(z_1-z_{2})}{x_1}\right)
Y'_{P(z_{1}-z_{2})}(v, x)
\mu^{(2)}_{\lambda, w_{(3)}}$$
is well defined and we have
\begin{eqnarray}
\lefteqn{\tau_{P(z_{1}-z_{2})}\biggl(
x^{-1}_1 \delta\left(\frac{x^{-1}-(z_1-z_{2})}{x_1}\right)
Y_{t}(v, x)\biggr)
\mu^{(2)}_{\lambda, w_{(3)}}=}\nno\\
&&=x^{-1}_1 \delta\left(\frac{x^{-1}-(z_1-z_{2})}{x_1}\right)
Y'_{P(z_{1}-z_{2})}(v, x)
\mu^{(2)}_{\lambda, w_{(3)}}.
\end{eqnarray}
Similarly, for any $v\in V$ and any $w_{(1)}\in W_{1}$,
$$x^{-1}_1 \delta\left(\frac{x^{-1}-z_{2}}{x_1}\right)
Y'_{P(z_{2})}(v, x)
\mu^{(1)}_{\lambda, w_{(1)}}$$
is well defined and we have
\begin{eqnarray}
\lefteqn{\tau_{P(z_{2})}\biggl(
x^{-1}_1 \delta\left(\frac{x^{-1}-z_{2}}{x_1}\right)
Y_{t}(v, x)\biggr)
\mu^{(1)}_{\lambda, w_{(1)}}=}\nno\\
&&=x^{-1}_1 \delta\left(\frac{x^{-1}-z_{2}}{x_1}\right)
Y'_{P(z_{1}-z_{2})}(v, x)
\mu^{(1)}_{\lambda, w_{(1)}}.
\end{eqnarray}
\end{lemma}

The proof of this result will be given in the next section.

We now assume that $V$ is rational and all irreducible $V$-modules are
${\Bbb R}$-graded. Note that in this case all $V$-modules are
${\Bbb R}$-graded. We
 show that $\gamma(F_{1}; I, F_{2})'(w'_{(4)})$ satisfies
the $P(z_{2})$-local grading-restriction
condition. We assume that $w'_{(4)}$ is homogeneous.
For any homogeneous element $w_{(1)}\in W_{1}$,
let $\alpha_{n}(w'_{(4)}, w_{(1)})$, $n\in {\Bbb C}$, be the elements
of $(W_{2}\otimes W_{3})^{*}$ defined by
\begin{equation}
(\alpha_{n}(w'_{(4)}, w_{(1)}))(w_{(2)}\otimes w_{(3)})=
\langle w'_{(4)}, F_{1}(w_{(1)}
\otimes P_{n}(F_{2}(w_{(2)}\otimes w_{(3)})))\rangle_{W_{4}}.
\end{equation}
By the commutator formula for $L(0)$ and $F_{2}$, the definition
$L'_{P(z_{2})}(0)$ and the definition of intertwining map,
it is easy to verify that $\alpha_{n}(w'_{(4)}, w_{(1)})$,
$n\in {\Bbb C}$, are  $P(z_{2})$-weight vectors and are in $W_{2}
\hboxtr_{P(z_{2})}W_{3}$. {}From our assumption,
the series $\sum_{n\in {\Bbb C}}\alpha_{n}(w'_{(4)}, w_{(1)})$ is absolutely
convergent to $\mu^{(1)}_{\gamma(F_{1}; I, F_{2})'(w'_{(4)}), w_{(1)}}$.
To show that this series is indexed by sequences of strictly
increasing real numbers,
we need the following lemma, whose proof will be given in the next
section:

\begin{lemma}
Let $W_{1}$, $W_{2}$ and $W_{3}$ be three modules for a rational vertex
operator algebra $V$ and ${\cal Y}: W_{1}\otimes W_{2}\to W_{3}\{ x\}$
an intertwining operator of type  ${W_{3}}\choose {W_{1}W_{2}}$. Then
there exist finitely many complex numbers $h_{1}, \dots, h_{p}$ such that
${\cal Y}$ is in fact a map {}from $W_{1}\otimes W_{2}$ to
$x^{h_{1}}W_{3}[[x, x^{-1}]]+ \cdots +x^{h_{p}}W_{3}[[x, x^{-1}]]$.
If all irreducible $V$-modules are ${\Bbb R}$-graded, $h_{1}, \dots, h_{p}$
are real numbers.
\end{lemma}

Combining this lemma with the lower truncation condition for the intertwining
operator ${\cal Y}_{2}$ corresponding to $F_{2}$, we see that there exists
 a sequence $\{n_{i}\}_{i\in {\Bbb Z}_{+}}$  of strictly increasing
real numbers
such that
\begin{equation}
{\cal Y}_{2}(w_{(2)}, x_{2})w_{(3)}
=\sum_{i\in {\Bbb Z}_{+}}(w_{(2)})_{-n_{i}-1}w_{(3)}
x_{2}^{n_{i}}.
\end{equation}
Thus for any $w_{(2)}\in W_{2}$ and $w_{(3)}\in W_{3}$,
 if $n\ne \wt w_{(2)}+n_{i}+1+\wt w_{(3)}$, $i\in {\Bbb Z}_{+}$,
\begin{eqnarray}
\lefteqn{(\alpha_{n}(w'_{(4)}, w_{(1)}))(w_{(2)}\otimes w_{(3)})=}\nno\\
&&=\langle w'_{(4)}, F_{1}(w_{(1)}
\otimes P_{n}(F_{2}(w_{(2)}\otimes w_{(3)})))\rangle_{W_{4}}\nno\\
&&=\langle w'_{(4)}, F_{1}(w_{(1)}
\otimes P_{n}({\cal Y}_{2}(w_{(2)}, x_{2})w_{(3)})))\rangle_{W_{4}}
\lbar_{x_{2}=z_{2}}\nno\\
&&=\langle w'_{(4)}, F_{1}(w_{(1)}
\otimes P_{n}(\sum_{i\in {\Bbb Z}_{+}}(w_{(2)})_{-n_{i}-1}w_{(3)}
x_{2}^{n_{i}})))\rangle_{W_{4}}
\lbar_{x_{2}=z_{2}}\nno\\
&&=0.
\end{eqnarray}
So $\sum_{n\in {\Bbb C}}\alpha_{n}(w'_{(4)}, w_{(1)})$ is indexed
by sequences of strictly increasing real numbers.
This proves that $\gamma(F_{1}; I, F_{2})'(w'_{(4)})$ satisfies
the $P(z_{2})$-local grading-restriction condition. Similarly we can show that
$\gamma(F_{4}; F_{3}, I)'(w'_{(4)})$
satisfies the $P(z_{1}-z_{2})$-local grading restriction condition.

Assume in addition that
the vertex operator algebra $V$ has the property that for any modules
$W_{1}$, $W_{2}$, $W_{3}$, $W_{4}$ and $W_{5}$, any $z_{1}, z_{2}\in {\Bbb C}$
satisfying (14.7),
any
$P(z_{1})$- and $P(z_{2})$-intertwining maps $F_{1}$ and $F_{2}$ of the types
as above and any $w'_{(4)}\in W'_{4}$,
the element $\gamma(F_{1}; I, F_{2})'(w'_{(4)})$  in
$(W_{1}\otimes W_{2}\otimes W_{3})^{*}$  satisfies the
$P(z_{1}-z_{2})$-local grading-restriction condition.
Let ${\cal Y}_{1}$ and ${\cal Y}_{2}$ be the intertwining operators
corresponding to $F_{1}$ and $F_{2}$, respectively.
For $z\in {\Bbb C}^{\times}$,
we define the $P(z_{1}-z_{2}+zz_{2})$- and
$P(zz_{2})$-intertwining maps $F^{z}_{1}$ and $F^{z}_{2}$
of the same types as those of
$F_{1}$ and $F_{2}$, respectively,  by
\begin{eqnarray}
F_{1}^{z}&=&{\cal Y}(\cdot, x)\cdot \lbar_{x=z_{1}-z_{2}+zz_{2}},\nno\\
F_{2}^{z}&=&{\cal Y}(\cdot, x)\cdot \lbar_{x=zz_{2}}.
\end{eqnarray}
If $0<|z|< \frac{|z_{1}-z_{2}|}{2|z_{2}|}$, we have
$|z_{1}-z_{2}+zz_{2}|>|zz_{2}|>0$. Thus the product
$\gamma(F^{z}_{1}; I, F^{z}_{2})$ exists.

By assumption, there is a series
$\sum_{n\in {\Bbb C}}\alpha_{n}(w'_{(4)}, w_{(3)})$ indexed by
sequences of strictly increasing real numbers
of $P(z_{1}-z_{2})$-weight vectors  in $(W_{1}\otimes W_{2})^{*}$
converging  absolutely to the element
\begin{equation}
\mu^{(2)}_{\gamma(F_{1}; I, F_{2})'
(w'_{(4)}), w_{(3)}}
\end{equation}
of $(W_{1}\otimes W_{2})^{*}$ when (14.7) holds. By the definition of
(14.42), $L'_{P(z_{1}-z_{2})}(0)$ and $\gamma(F_{1}; I, F_{1})'$, we see that
for any $w_{(1)}\in W_{1}$,
$w_{(2)}\in W_{2}$, $w_{(3)}\in W_{3}$, $w'_{(4)}\in W'_{4}$ and formal
variable $x$,
\begin{eqnarray}
\lefteqn{((1-x)^{-L'_{P(z_{1}-z_{2})}(0)}
\mu^{(2)}_{\gamma(F_{1}; I, F_{2})'
(w'_{(4)}), w_{(3)}})(w_{(1)}\otimes w_{(2)})=}\nno\\
&&=\mu^{(2)}_{\gamma(F_{1}; I, F_{2})'
(w'_{(4)}), w_{(3)}}((1-x)^{-(z_{1}-z_{2})L(-1)-L(0)}w_{(1)}\otimes
(1-x)^{-L(0)}w_{(2)})\nno\\
&&=(\gamma(F_{1}; I, F_{1})'(w'_{(4)}))((1-x)^{-(z_{1}-z_{2})L(-1)-L(0)}
w_{(1)}\otimes\nno\\
&&\hspace{5em}\otimes
(1-x)^{-L(0)}w_{(2)}\otimes w_{(3)})\nno\\
&&=\langle w'_{(4)}, {\cal Y}_{1}((1-x)^{-(z_{1}-z_{2})L(-1)-L(0)}
w_{(1)}, x_{1})\cdot\nno\\
&&\hspace{5em}\cdot {\cal Y}_{2}(
(1-x)^{-L(0)}w_{(2)}, x_{2})
w_{(3)}\rangle_{W_{4}}\lbar_{x_{1}=z_{1}, x_{2}=z_{2}}.
\end{eqnarray}
Using the commutator formulas for $L(0)$, $L(-1)$ and intertwining operators,
it is easy to show by direct calculations that the right-hand side of (14.43)
is equal to
\begin{eqnarray}
\lefteqn{\langle w'_{(4)}, (1-x)^{-(L(0)-z_{2}L(-1))}{\cal Y}_{1}(
w_{(1)}, x_{1})\cdot}\nno\\
&&\hspace{3em}\cdot{\cal Y}_{2}(
w_{(2)}, x_{2})(1-x)^{L(0)-z_{2}L(-1)}
w_{(3)}\rangle_{W_{4}}\lbar_{x_{1}=z_{1}, x_{2}=z_{2}}.
\end{eqnarray}
By Lemma 9.3 in \cite{HL4} which gives the formula
$$(1-x)^{L(0)-z_{2}L(-1)}=e^{z_{2}xL(-1)}(1-x)^{L(0)},$$
(14.44) is equal to
\begin{eqnarray}
\lefteqn{\langle w'_{(4)}, (1-x)^{-L(0)}e^{-z_{2}xL(-1)}{\cal Y}_{1}(
w_{(1)}, x_{1})\cdot}\nno\\
&&\hspace{3em}\cdot {\cal Y}_{2}(
w_{(2)}, x_{2})e^{z_{2}xL(-1)}(1-x)^{L(0)}
w_{(3)}\rangle_{W_{4}}\lbar_{x_{1}=z_{1}, x_{2}=z_{2}}\nno\\
&&=\langle (1-x)^{-L(0)}w'_{(4)}, {\cal Y}_{1}(
w_{(1)}, x_{1}-z_{2}x)\cdot\nno\\
&&\hspace{3em}\cdot{\cal Y}_{2}(
w_{(2)}, x_{2}-z_{2}x)(1-x)^{L(0)}
w_{(3)}\rangle_{W_{4}}\lbar_{x_{1}=z_{1}, x_{2}=z_{2}}\nno\\
&&=\langle (1-x)^{-L(0)}w'_{(4)}, {\cal Y}_{1}(
w_{(1)}, x_{1})\cdot\nno\\
&&\hspace{3em}\cdot{\cal Y}_{2}(
w_{(2)}, x_{2})(1-x)^{L(0)}
w_{(3)}\rangle_{W_{4}}\lbar_{x_{1}=z_{1}-z_{2}x, x_{2}=z_{2}-z_{2}x}.
\end{eqnarray}

The calculations above shows that the left-hand side of (14.43) is equal
to the right-hand side of (14.45).
Since in the right-hand side of (14.45), we can substitute $1-e^{\log z}$ for
$x$ when $0<|z|<\frac{|z_{1}-z_{2}|}{2|z_{2}|}$,
we can also substitute $1-e^{\log z}$ for
$x$ when $0<|z|<\frac{|z_{1}-z_{2}|}{2|z_{2}|}$
in the left-hand side of (14.43), and after the substitutions, the left-hand
side of (14.43) and the right-hand side of (14.45) are still equal. So we have
\begin{eqnarray}
\lefteqn{(e^{-(\log z)L'_{P(z_{1}-z_{2})}(0)}
\mu^{(2)}_{\gamma(F_{1}; I, F_{2})'
(w'_{(4)}), w_{(3)}})(w_{(1)}\otimes w_{(2)})=}\nno\\
&&=\langle e^{-(\log z)L'(0)}w'_{(4)}, {\cal Y}_{1}(
w_{(1)}, x_{1})\cdot\nno\\
&&\hspace{3em}\cdot{\cal Y}_{2}(
w_{(2)}, x_{2})e^{(\log z)L(0)}
w_{(3)}\rangle_{W_{4}}\lbar_{x_{1}=z_{1}-z_{2}+zz_{2}, x_{2}=zz_{2}}\nno\\
&&=(\mu^{(2)}_{\gamma(F^{z}_{1}; I, F^{z}_{2})'
(e^{-(\log z)L'(0)}w'_{(4)}), e^{(\log z)L(0)}w_{(3)}})(w_{(1)}\otimes
w_{(2)}).
\end{eqnarray}
Since $\alpha_{n}(w'_{(4)}, w_{(3)})$, $n\in {\Bbb C}$, are weight vectors,
we have
\begin{eqnarray}
\lefteqn{e^{-(\log z)L'_{P(z_{1}-z_{2})}(0)}
\mu^{(2)}_{\gamma(F_{1}; I, F_{2})'
(w'_{(4)}), w_{(3)}}=}\nno\\
&&=\sum_{n\in {\Bbb C}}\alpha_{n}(w'_{(4)}, w_{(3)})e^{-\mbox{\scriptsize wt}\;
 \alpha_{n}(w'_{(4)},
w_{(3)})\log z}
\end{eqnarray}
when (14.7) holds and $0<|z|< \frac{|z_{1}-z_{2}|}{2|z_{2}|}$.
By (14.46) and (14.47), we obtain
\begin{eqnarray}
\lefteqn{\mu^{(2)}_{\gamma(F^{z}_{1}; I, F^{z}_{2})'
(w'_{(4)}), w_{(3)}}=}\nno\\
&&=\sum_{n\in {\Bbb C}}\alpha_{n}(w'_{(4)}, w_{(3)})
e^{-(\mbox{\scriptsize wt}\; \alpha_{n}(w'_{(4)},
w_{(3)})+\mbox{\scriptsize wt}\; w'_{(4)}-\mbox{\scriptsize wt}\;
 w_{(3)})\log z}
\end{eqnarray}
when (14.7) holds and $0
<|z|< \frac{|z_{1}-z_{2}|}{2|z_{2}|}$. When
$$|z_{1}-z_{2}+zz_{2}|>|zz_{2}|>|z_{1}-z_{2}|>0,$$ by Lemma 14.3,
$\mu^{(2)}_{\gamma(F^{z}_{1}; I, F^{z}_{2})'
(w'_{(4)}), w_{(3)}}$
satisfies (14.36). Since the open set on the $z$-plane given by
$|z_{1}-z_{2}+zz_{2}|>|zz_{2}|>0$ is connected and
$\mu^{(2)}_{\gamma(F^{z}_{1}; I, F^{z}_{2})'
(w'_{(4)}), w_{(3)}}$ is defined and analytic in $z$ when
$|z_{1}-z_{2}+zz_{2}|>|zz_{2}|>0$, the coefficients in powers of
$x$ and $x_{1}$
of both sides of
(14.36) with $\lambda=\gamma(F^{z}_{1}; I, F^{z}_{2})'
(w'_{(4)})$, as (multi-valued) functions of $z$ can be
analytically extended to $|z_{1}-z_{2}+zz_{2}|>|zz_{2}|>0$
and are still equal.  When (14.7) holds and  $0
<|z|< \frac{|z_{1}-z_{2}|}{2|z_{2}|}$, these coefficients can be expanded as
series in powers of $e^{\log z}$. By part (b) of the $P(z_{1}-z_{2})$-local
grading-restriction condition which we assumed to be satisfied by
$\gamma(F_{1}; I, F_{2})$,
$\alpha_{n}(w'_{(4)},
w_{(3)})$, $n\in {\Bbb C}$, satisfies the $P(z_{1}-z_{2})$-lower truncation
condition. Thus for any $n\in {\Bbb C}$,
$$x^{-1}_1 \delta\left(\frac{x^{-1}-(z_1-z_{2})}{x_1}\right)
Y'_{P(z_{1}-z_{2})}(v, x)\alpha_{n}(w'_{(4)}, w_{(3)})$$
exists. By (14.48), the coefficients in powers of $x$ and $x_{1}$ of
\begin{eqnarray*}
\lefteqn{x^{-1}_1 \delta\left(\frac{x^{-1}-(z_1-z_{2})}{x_1}\right)
Y'_{P(z_{1}-z_{2})}(v, x)\cdot}\\
&&\quad \cdot
\left(\sum_{n\in {\Bbb C}}\alpha_{n}(w'_{(4)}, w_{(3)})
e^{-(\mbox{\scriptsize $\wt$} \alpha_{n}(w'_{(4)},
w_{(3)})+\mbox{\scriptsize $\wt$} w'_{(4)}-
\mbox{\scriptsize $\wt$} w_{(3)})\log z}\right)
\end{eqnarray*}
exist since they are in fact the expansions  as series in powers
of $e^{\log z}$
of the analytic extensions of
the coefficients in powers of
$x$ and $x_{1}$ of
the right-hand side of (14.36) with $\lambda=\gamma(F^{z}_{1}; I, F^{z}_{2})'
(w'_{(4)})$.
We obatin
\begin{eqnarray}
\lefteqn{\tau_{P(z_{1}-z_{2})}\biggl(
x^{-1}_1 \delta\left(\frac{x^{-1}-(z_1-z_{2})}{x_1}\right)
Y_{t}(v, x)\biggr)\cdot}\nno\\
&&\quad \cdot \left(\sum_{n\in {\Bbb C}}\alpha_{n}(w'_{(4)}, w_{(3)})
e^{-(\mbox{\scriptsize $\wt$} \alpha_{n}(w'_{(4)},
w_{(3)})+\mbox{\scriptsize $\wt$} w'_{(4)}-
\mbox{\scriptsize $\wt$} w_{(3)})\log z}\right)\nno\\
&&=x^{-1}_1 \delta\left(\frac{x^{-1}-(z_1-z_{2})}{x_1}\right)
Y'_{P(z_{1}-z_{2})}(v, x)\cdot\nno\\
&&\quad \cdot
\left(\sum_{n\in {\Bbb C}}\alpha_{n}(w'_{(4)}, w_{(3)})
e^{-(\mbox{\scriptsize $\wt$} \alpha_{n}(w'_{(4)},
w_{(3)})+\mbox{\scriptsize $\wt$} w'_{(4)}-
\mbox{\scriptsize $\wt$} w_{(3)})\log z}\right)\nno\\
&&
\end{eqnarray}
 when
(14.7) holds and $0<|z|< \frac{|z_{1}-z_{2}|}{2|z_{2}|}$.
 We need the following
lemma:

\begin{lemma}
Let $f(z)$ be a multi-valued complex analytic function of $z$ in a neighborhood
of $0$ with $0$ deleted. If for $z$ in this neighborhood such that
$\arg z$ satisfies
$0\le \arg z<2\pi$, $f(z)$ can be expanded as
\begin{equation}
\sum_{i\in {\Bbb Z}_{+}}a_{i}e^{m_{i}\log z}
\end{equation}
where $\{m_{i}\}_{i\in {\Bbb Z}_{+}}$ is a
sequence of strictly increasing real numbers,
 then the sequence $\{m_{i}\}_{i\in {\Bbb Z}_{+}}$ and the coefficients
$a_{i}$, $i\in {\Bbb Z}_{+}$,  are  uniquely determined {}by $f(z)$.
\end{lemma}

The proof of this result will be given in the next section.

By assumption, we see that for $z$ in this neighborhood such that
$\arg z$ satisfies
$0\le \arg z<2\pi$,
 $$\biggl(\sum_{n\in {\Bbb C}}\alpha_{n}(w'_{(4)}, w_{(3)})
e^{-(\mbox{\scriptsize $\wt$} \alpha_{n}(w'_{(4)},
w_{(3)})+\mbox{\scriptsize $\wt$}
w'_{(4)}-\mbox{\scriptsize $\wt$}
w_{(3)})\log z}\biggr)(w_{(1)}\otimes w_{(2)})$$
for any
$w_{(1)}\in W_{1}$ and $w_{(2)}\in W_{2}$, as a multi-valued function of
$z$, satisfies the condition in Lemma 14.5.
Thus the coefficients of both sides of (14.49) in powers
of $x$ and $x_{1}$ applied to $w_{(1)}\otimes w_{(2)}$ also satisfies
the conditions of Lemma 14.5. By Lemma 14.5, the coefficients of
the left-hand side and of the right-hand  side of (14.49) in powers of
$e^{\log z}$ are equal. Since we have already shown that
for any $n\in {\Bbb C}$,
$$x^{-1}_1 \delta\left(\frac{x^{-1}-(z_1-z_{2})}{x_1}\right)
Y'_{P(z_{1}-z_{2})}(v, x)\alpha_{n}(w'_{(4)}, w_{(3)})$$
exists, we obtain
\begin{eqnarray}
\lefteqn{\tau_{P(z_{1}-z_{2})}\biggl(
x^{-1}_1 \delta\left(\frac{x^{-1}-(z_1-z_{2})}{x_1}\right)
Y_{t}(v, x)\biggr)
\alpha_{n}(w'_{(4)}, w_{(3)})
=}\nno\\
&&=x^{-1}_1 \delta\left(\frac{x^{-1}-(z_1-z_{2})}{x_1}\right)
Y'_{P(z_{1}-z_{2})}(v, x)
\alpha_{n}(w'_{(4)}, w_{(3)}).
\end{eqnarray}
proving that $\alpha_{n}(w'_{(4)},
w_{(3)})$, $n\in {\Bbb C}$, satisfy the $P(z_{1}-z_{2})$-compatibility
condition.

So $\alpha_{n}(w'_{(4)}, w_{(3)})$,
$n\in {\Bbb C}$,
are in fact elements of $W_{1}\hboxtr_{P(z_{1}-z_{2})}W_{2}$ and (14.42)
can be regarded as an element of
$\overline{W_{1}\hboxtr_{P(z_{1}-z_{2})}W_{2}}$.
Let ${\cal Y}_{3}$
be  the intertwining operator correponding to the
$P(z_{1}-z_{2})$-intertwining map $\boxtimes_{P(z_{1}-z_{2})}$.
We have:

\begin{lemma}
Let  $G \in \hom(W_{1}\boxtimes_{P(z_{1}-z_{2})} W_{2},
(W'_{(4)}\otimes W_{3})^{*})$ be defined by
\begin{equation}
(G(w))(w'_{(4)}\otimes w_{(3)})=\langle \mu^{(2)}_{\gamma(F_{1}; I, F_{2})'
(w'_{(4)}), w_{(3)}}, \nu\rangle_{W_{1}\shboxtr_{P(z_{1}-z_{2})} W_{2}}
\end{equation}
for $w\in W_{1}\boxtimes_{P(z_{1}-z_{2})} W_{2}$, $w'_{(4)}\in W'_{4}$ and
$w_{(3)}\in W_{3}$. Then $G$ intertwinines the two actions
$\tau_{W_{1}\shboxtr_{P(z_{1}-z_{2})} W_{2}}$ and $\tau_{Q(z_{2})}$ of
$V\otimes \iota_{+} {\Bbb C}[t, t^{-1}, (z+t)^{-1}]$ on
$W_{1}\hboxtr_{P(z_{1}-z_{2})} W_{2}$ and on $(W'_{(4)}\otimes W_{3})^{*}$.
\end{lemma}

This result will be proved in the next section.

By this lemma and Remark 5.4 in \cite{HL3} (Remark I.10 in \cite{HL6}),
there exists an intertwining operator ${\cal Y}_{4}$ of type
${W_{4}}\choose {W_{1}\boxtimes_{P(z_{1}-z_{2})} W_{2}\;W_{3}}$ such that
\begin{eqnarray}
(G(w))(w'_{(4)}\otimes w_{(3)})=
\langle w'_{(4)}, {\cal Y}_{4}(w, x_{2})w_{(3)}\rangle_{W_{4}}
\lbar_{x_{2}=z_{2}}
\end{eqnarray}
for $w\in W_{1}\hboxtr_{P(z_{1}-z_{2})} W_{2}$ and $w'_{(4)}\in W'_{4}$ and
$w_{(3)}\in W_{3}$.
By Remark 13.2 in \cite{HL6}, (14.52), (14.53) and the absolute convergence
of iterates of intertwining operators,
\begin{eqnarray}
\lefteqn{\mu^{(2)}_{\gamma(F_{1}; I, F_{2})'
(w'_{(4)}), w_{(3)}}(w_{(1)}\otimes w_{(2)})}\nno\\
&&=\langle \mu^{(2)}_{\gamma(F_{1}; I, F_{2})'
(w'_{(4)}), w_{(3)}}, {\cal Y}_{3}(w_{(1)}, x_{0})w_{(2)}
\rangle_{W_{1}\hboxtr_{P(z_{1}-z_{2})} W_{2}}\lbar_{x_{0}=z_{1}-z_{2}}\nno\\
&&=(G({\cal Y}_{3}(w_{(1)}, x_{0})w_{(2)}))(w'_{(4)}, w_{(3)})
\lbar_{x_{0}=z_{1}-z_{2}}\nno\\
&&=\langle w'_{(4)}, {\cal Y}_{4}({\cal Y}_{3}(w_{(1)}
x_{0})w_{(2)}, x_{2})w_{(3)}\rangle_{W_{4}}\lbar_{x_{0}=z_{1}-z_{2},
x_{2}=z_{2}}
\end{eqnarray}
for any $w'_{(4)}\in W_{4}$,
$w_{(1)}\in W_{1}$, $w_{(2)}\in W_{2}$ and $w_{(3)}\in W_{3}$ and any
$z_{1}, z_{2}\in {\Bbb C}$ satisfying (14.7).
Let ${\cal Y}_{1}$ and ${\cal Y}_{2}$ be the intertwining operators
corresponding to $F_{1}$ and $F_{2}$, respectively. Then
for any $w'_{(4)}\in W_{4}$,
$w_{(1)}\in W_{1}$, $w_{(2)}\in W_{2}$ and $w_{(3)}\in W_{3}$ and any
$z_{1}, z_{2}\in {\Bbb C}$ satisfying (14.7), we obtain
\begin{eqnarray}
\lefteqn{\langle w'_{(4)}, {\cal Y}_{1}(w_{(1)},
x_{1}){\cal Y}_{2}(w_{(2)}, x_{2})w_{(3)}\rangle_{W_{4}}\lbar_{x_{1}=z_{1},
x_{2}=z_{2}}}
\nno\\
&&=\langle w'_{(4)}, {\cal Y}_{4}({\cal Y}_{3}(w_{(1)},
x_{0})w_{(2)}, x_{2})w_{(3)}\rangle_{W_{4}}\lbar_{x_{0}=z_{1}-z_{2},
x_{2}=z_{2}}.
\end{eqnarray}

Similarly, if in addtion to the assumptions that $V$ has one of the
properties in Proposition 14.1, that $V$ is rational and that
all irreducible $V$-modules are ${\Bbb R}$-graded, we assume that
the vertex operator algebra $V$ has the property that for any modules
$W_{1}$, $W_{2}$, $W_{3}$, $W_{4}$ and $W_{5}$, any $z_{1}, z_{2}\in {\Bbb C}$
satisfying (14.7),
any
$P(z_{1}-z_{2})$- and $P(z_{2})$-intertwining maps $F_{3}$ and $F_{4}$ of
the types
as above and any $w'_{(4)}\in W'_{4}$,
the element $\gamma(F_{4}; F_{3}, I)'(w'_{(4)})$  of
$(W_{1}\otimes W_{2}\otimes W_{3})^{*}$ satisfies the
$P(z_{2})$-local grading-restriction condition, we can show that
for any intertwining operators ${\cal Y}_{3}$ and ${\cal Y}_{4}$ of the same
types as those for $F_{3}$ and $F_{4}$, respectively,
we can show that there exist intertwining operators
${\cal Y}_{1}$ and ${\cal Y}_{2}$ such that (14.55) holds
for any $w'_{(4)}\in W_{4}$,
$w_{(1)}\in W_{1}$, $w_{(2)}\in W_{2}$ and $w_{(3)}\in W_{3}$ and any
$z_{1}, z_{2}\in {\Bbb C}$ satisfying (14.7).
{}.

We still assume, in addtion to the assumptions that $V$ has one of the
properties in Proposition 14.1, that $V$ is rational and
that all irreducible $V$-modules are ${\Bbb R}$-graded, that
the vertex operator algebra $V$ has the property that for any modules
$W_{1}$, $W_{2}$, $W_{3}$, $W_{4}$ and $W_{5}$, any $z_{1}, z_{2}\in {\Bbb C}$
satisfying (14.7),
any
$P(z_{1})$- and $P(z_{2})$-intertwining maps $F_{1}$ and $F_{2}$ of the types
as above and any $w'_{(4)}\in W'_{4}$,
the element $\gamma(F_{1}; I, F_{2})'(w'_{(4)})$  in
$(W_{1}\otimes W_{2}\otimes W_{3})^{*}$  satisfies the
$P(z_{1}-z_{2})$-local grading restriction condition.
For  any $P(z_{1}-z_{2})$-intertwining map $F_{3}$ and
$P(z_{2})$-intertwining map
$F_{4}$ of the types as above, when (14.7) holds,
using (14.6) and $\Omega_{0}$ and $\Omega_{-1}$, we have
\begin{eqnarray}
\lefteqn{(\gamma(F_{4}; F_{3}, I)'(w'_{(4)}))(w_{(1)}\otimes w_{(2)}\otimes
w_{(3)})=}\nno\\
&&=\langle e^{z_{2}L'(1)}w'_{(4)},
\Omega_{-1}({\cal Y}_{4})
(w_{(3)}, x_{2})\cdot\nno\\
&&\hspace{4em}\cdot
\Omega_{0}(\Omega_{-1}({\cal Y}_{3}))(w_{(1)},
x_{0})w_{(2)}\rangle_{W_{4}}\lbar_{x_{0}=z_{1}-z_{2},\;
x_{2}=e^{\pi i}z_{2}}\nno\\
&&=\langle e^{z_{2}L'(1)}w'_{(4)},
\Omega_{-1}({\cal Y}_{4})
(w_{(3)}, x_{2})e^{x_{0}L(-1)}\cdot \nno\\
&&\hspace{4em}\cdot \Omega_{-1}({\cal Y}_{3})(w_{(2)}, e^{\pi i}x_{0})
w_{(1)}\rangle_{W_{4}}\lbar_{x_{0}=z_{1}-z_{2},\;
x_{2}=e^{\pi i}z_{2}}\nno\\
&&=\langle e^{z_{2}L'(1)}w'_{(4)},
e^{x_{0}L(-1)}\Omega_{-1}({\cal Y}_{4})
(w_{(3)}, x_{2}-x_{0})\cdot \nno\\
&&\hspace{4em}\cdot\Omega_{-1}({\cal Y}_{3})(w_{(2)}, e^{\pi i}x_{0})
w_{(1)}\rangle_{W_{4}}\lbar_{x_{0}=z_{1}-z_{2},\;
x_{2}=e^{\pi i}z_{2}}\nno\\
&&=\langle e^{z_{1}L'(1)}w'_{(4)},
\Omega_{-1}({\cal Y}_{4})
(w_{(3)}, x_{1})\cdot \nno\\
&&\hspace{4em}\cdot\Omega_{-1}({\cal Y}_{3})(w_{(2)}, x_{0})
w_{(1)}\rangle_{W_{4}}\lbar_{x_{1}=e^{\pi i}z_{1},\;
x_{0}=e^{\pi i}(z_{1}-z_{2})}.
\end{eqnarray}
By the proof of (14.55) above, there exist a module $W_{6}$ and
intertwining operators
${\cal Y}_{5}$ and ${\cal Y}_{6}$ of type ${W_{6}}\choose {W_{3}W_{2}}$
and ${W_{4}}\choose {W_{6}W_{1}}$, respectively, such that when
$|z_{1}|>|z_{1}-z_{2}|>|z_{2}|>0$,
\begin{eqnarray}
\lefteqn{\langle e^{z_{1}L'(1)}w'_{(4)},
{\cal Y}_{6}({\cal Y}_{5}(w_{(3)}, x_{2})w_{(2)}, x_{0})
w_{(1)}\rangle_{W_{4}}\lbar_{x_{2}=e^{\pi i}z_{2},\;
x_{0}=e^{\pi i}(z_{1}-z_{2})}}\nno\\
&&=\langle e^{z_{1}L'(1)}w'_{(4)},
\Omega_{-1}({\cal Y}_{4})
(w_{(3)}, x_{1})\cdot \nno\\
&&\hspace{4em}\cdot\Omega_{-1}({\cal Y}_{3})(w_{(2)}, x_{0})
w_{(1)}\rangle_{W_{4}}\lbar_{x_{1}=e^{\pi i}z_{1},\;
x_{0}=e^{\pi i}(z_{1}-z_{2})}.
\end{eqnarray}
On the other hand, when $|z_{1}|>|z_{1}-z_{2}|>|z_{2}|>0$,
\begin{eqnarray}
\lefteqn{\langle e^{z_{1}L'(1)}w'_{(4)},
{\cal Y}_{6}({\cal Y}_{5}(w_{(3)}, x_{2})w_{(2)}, x_{0})
w_{(1)}\rangle_{W_{4}}\lbar_{x_{2}=e^{\pi i}z_{2},\;
x_{0}=e^{\pi i}(z_{1}-z_{2})}}\nno\\
&&=\langle e^{z_{1}L'(1)}w'_{(4)},
{\cal Y}_{6}(\Omega_{0}(\Omega_{-1}({\cal Y}_{5}))(w_{(3)},
x_{2})\cdot\nno\\
&&\hspace{4em}\cdot w_{(2)}, x_{0})
w_{(1)}\rangle_{W_{4}}\lbar_{x_{2}=e^{\pi i}z_{2},\;
x_{0}=e^{\pi i}(z_{1}-z_{2})}\nno\\
&&=\langle e^{z_{1}L'(1)}w'_{(4)},
{\cal Y}_{6}(e^{x_{2}L(-1)}\Omega_{-1}({\cal Y}_{5})(w_{(2)},
e^{\pi i}x_{2})\cdot \nno\\
&&\hspace{4em}\cdot w_{(3)}, x_{0})
w_{(1)}\rangle_{W_{4}}\lbar_{x_{2}=e^{\pi i}z_{2},\;
x_{0}=e^{\pi i}(z_{1}-z_{2})}\nno\\
&&=\langle e^{z_{1}L'(1)}w'_{(4)},
{\cal Y}_{6}(\Omega_{-1}({\cal Y}_{5})(w_{(2)},
e^{\pi i}x_{2})\cdot \nno\\
&&\hspace{4em}\cdot w_{(3)}, x_{0}+x_{2})
w_{(1)}\rangle_{W_{4}}\lbar_{x_{2}=e^{\pi i}z_{2},\;
x_{0}=e^{\pi i}(z_{1}-z_{2})}\nno\\
&&=\langle e^{z_{1}L'(1)}w'_{(4)},
{\cal Y}_{6}(\Omega_{-1}({\cal Y}_{5})(w_{(2)},
x_{2})w_{(3)}, x_{1})
w_{(1)}\rangle_{W_{4}}\lbar_{x_{1}=e^{\pi i}z_{1}, \; x_{2}=e^{2\pi i}z_{2}}.
\nno\\
&&
\end{eqnarray}
By the $L(-1)$-derivative property for intertwining
operators, both sides
of (14.56), (14.57) and (14.58) are analytic (multi-valued)
functions of $z_{1}$
and $z_{2}$. By these three formulas, we see that the left-hand side of
(14.56) is equal to the value of an branch of the analytic extension of
 the right-hand side of (14.58) at $(z_{1}, z_{2})$ satisfying (14.7).
Thus we can find $p, q\in {\Bbb Z}$ such that if we let
\begin{eqnarray}
{\cal Y}_{7}(\cdot, x)&=&\Omega_{-1}({\cal Y}_{5})(\cdot, e^{2\pi p i}x),\\
{\cal Y}_{8}(\cdot, x)&=&{\cal Y}_{6}(\cdot, e^{2\pi q i}x)
\end{eqnarray}
which are obviously also intertwining operators and
let $F_{7}$ and $F_{8}$ the $P(-z_{1})$- and $P(z_{2})$-intertwining maps
corresponding to ${\cal Y}_{7}$ and ${\cal Y}_{8}$, respectively,
then when (14.7) holds,
\begin{eqnarray}
\lefteqn{(\gamma(F_{4}; F_{3}, I)'(w'_{(4)}))(w_{(1)}\otimes w_{(2)}\otimes
w_{(3)})=}\nno\\
&&=\langle e^{z_{1}L'(1)}w'_{(4)},
{\cal Y}_{8}({\cal Y}_{7}(w_{(2)},
x_{2})w_{(3)}, x_{1})
w_{(1)}\rangle_{W_{4}}\lbar_{x_{1}=e^{\pi i}z_{1}, \; x_{2}=z_{2}}.
\nno\\
&&=(\gamma(F_{8}; F_{7}, I)'(e^{z_{1}L'(1)}w'_{(4)}))(w_{(2)}\otimes w_{(3)}
\otimes
w_{(1)})
\end{eqnarray}
Since  $(\gamma(F_{8}; F_{7}, I)'(e^{z_{1}L'(1)}w'_{(4)}))$ satisfies the
$P(z_{2})$-local grading-restriction condition, by
this formula, we see that $\gamma(F_{3}; F_{4}, I)'(w'_{(4)})$ satisfies
the $P(z_{2})$-local grading-restriction condition.
We have shown  half
of the following result:

\begin{propo}
Assume that $V$ is a rational vertex operator algebra, that all irreducible
$V$-modules are ${\Bbb R}$-greded and that $V$  has one
of the properties in Proposition 14.1.
Then the following two properties
 are equivalent:
\begin{enumerate}
\item For any $V$-modules $W_{1}$,
$W_{2}$, $W_{3}$, $W_{4}$ and $W_{5}$, any  nonzero complex
numbers $z_{1}$ and $z_{2}$ satisfying (14.7), any
$P(z_{1})$-intertwining map $F_{1}$ of type
${W_{4}}\choose {W_{1}W_{5}}$ and $P(z_{2})$-intertwining map $F_{2}$
of type ${W_{5}}\choose {W_{2}W_{3}}$ and any $w'_{(4)}\in W'_{4}$,
$\gamma(F_{1}; I, F_{2})'(w'_{(4)})\in (W_{1}\otimes W_{2}\otimes W_{3})^{*}$
 satisfies
the $P(z_{1}-z_{2})$-local grading-restriction condition.

\item For any $V$-modules $W_{1}$,
$W_{2}$, $W_{3}$, $W_{4}$ and $W_{5}$, any  nonzero complex
numbers $z_{1}$ and $z_{2}$ satisfying (14.7) and any
$P(z_{1}-z_{2})$-intertwining map $F_{3}$ of type
${W_{4}}\choose {W_{5}W_{3}}$ and $P(z_{2})$-intertwining map $F_{4}$
of type ${W_{5}}\choose {W_{1}W_{2}}$ and any $w'_{(4)}\in W'_{4}$,
$\gamma(F_{4}; F_{3}, I)
(w'_{(4)})\in (W_{1}\otimes W_{2}\otimes W_{3})^{*}$  satisfies
the $P(z_{2})$-local grading-restriction condition. \epf
\end{enumerate}

\end{propo}

The  other half of this result can be proved similarly.

In fact  much more has been proved. We have proved half of the following
result on the associativity of intertwining operators:

\begin{theo}
Assume that $V$ is a rational vertex operator algebra, that all irreducible
$V$-modules are ${\Bbb R}$-greded and that $V$  has one
of the properties in Proposition 14.1. If $V$ also has one of
the properties in Proposition 14.7, then the intertwining operators for
$V$ have the following two
associativity properties:
\begin{enumerate}

\item For any modules
$W_{1}$, $W_{2}$, $W_{3}$, $W_{4}$ and $W_{5}$ and
any
intertwining operators ${\cal Y}_{1}$ and ${\cal Y}_{2}$ of
 type ${W_{4}}\choose {W_{1}W_{5}}$ and ${W_{5}}\choose {W_{2}W_{3}}$,
respectively,
there exist a module $W_{6}$ and intertwining operators ${\cal Y}_{3}$
and ${\cal Y}_{4}$ of type ${W_{6}}\choose {W_{1}W_{2}}$ and
${W_{4}}\choose {W_{6}W_{3}}$, respectively, such that
for any $z_{1}, z_{2}\in {\Bbb C}$
satisfying (14.7), (14.55) holds for any
$w'_{(4)}\in W'_{4}$, $w_{(1)}\in W_{1}$,
$w_{(2)}\in W_{2}$ and $w_{(3)}\in W_{3}$.

\item For any modules
$W_{1}$, $W_{2}$, $W_{3}$, $W_{4}$ and $W_{6}$ and
any intertwining operators ${\cal Y}_{3}$
and ${\cal Y}_{4}$ of type ${W_{6}}\choose {W_{1}W_{2}}$ and
${W_{4}}\choose {W_{6}W_{3}}$, respectively, there exist a module $W_{5}$
and intertwining operators ${\cal Y}_{1}$ and ${\cal Y}_{2}$ of
 type ${W_{4}}\choose {W_{1}W_{5}}$ and ${W_{5}}\choose {W_{2}W_{3}}$,
respectively, such that
for any $z_{1}, z_{2}\in {\Bbb C}$
satisfying (14.7), (14.55) holds for any
$w'_{(4)}\in W'_{4}$, $w_{(1)}\in W_{1}$,
$w_{(2)}\in W_{2}$ and $w_{(3)}\in W_{3}$.

\end{enumerate}
Conversely,
If the intertwining operators for $V$ have
one of the associativity properties above, then $V$ has both of
the properties of Proposition 14.7, and in particular,
the intertwining operators for $V$ have both
of the associativity properties above.
\end{theo}
\pf
We need only to prove the second half of the result.
Assume that the intertwining operators for $V$ have the first associativity
property in the theorem. For any $V$-modules $W_{1}$,
$W_{2}$, $W_{3}$, $W_{4}$ and $W_{5}$, any  nonzero complex
numbers $z_{1}$ and $z_{2}$ satisfying (14.7), any
$P(z_{1})$-intertwining map $F_{1}$ of type
${W_{4}}\choose {W_{1}W_{5}}$ and $P(z_{2})$-intertwining map $F_{2}$
of type ${W_{5}}\choose {W_{2}W_{3}}$ and any $w'_{(4)}\in W'_{4}$,
by the first associativity property and the isomorphism between the
space of intertwining operators and the space of intertwining maps,
there exist a $V$-module $W_{6}$ and
$P(z_{1}-z_{2})$- and $P(z_{2})$-intertwining maps $F_{3}$ and $F_{4}$ of
type ${W_{4}}\choose {W_{6}W_{3}}$ and
${W_{6}}\choose {W_{1}W_{2}}$ such that
\begin{equation}
\gamma(F_{1}; I, F_{2})'=\gamma(F_{4}; F_{3}, I)'.
\end{equation}
We already know that for any $w'_{(4)}$, $\gamma(F_{4}; F_{3}, I)'(w'_{(4)})$
satisfies the $P(z_{1}-z_{2})$-local grading-restriction condition. So
$\gamma(F_{1}; I, F_{2})'(w'_{(4)})$ also satisfies the $P(z_{1}-z_{2})$-local
grading-restriction condition. Similarly we can prove the second half
when the intertwining operators for $V$ have
 the other associativity property.
\epfv

Now we assume that (14.7) holds and that $V$ is rational, that all irreducible
$V$-modules are ${\Bbb R}$-graded and that $V$
has either one of the properties in
Proposition 14.1 and also either one of the properties in Proposition 14.7.
Thus $V$ satisfies all the four properties by these propositions.

We denote the graded space of all elements of
$(W_{1}\otimes W_{2}\otimes W_{3})^{*}$
satisfying the $P(z_{1}, z_{2})$-compatibility condition,
the $P(z_{1}, z_{2})$-,
the $P(z_{2})$- and the
$P(z_{1}-z_{2})$-local grading-restriction conditions by $W_{P(z_{1}, z_{2})}$.

The product of the $P(z_{1})$-intertwining map $\boxtimes_{P(z_{1})}$
of type
$${W_{1}\boxtimes_{P(z_{1})}(W_{2}\boxtimes_{P(z_{2})}W_{3})}\choose
{W_{1}\quad\quad\quad W_{2}\boxtimes_{P(z_{2})}W_{3}}$$
and the $P(z_{2})$-intertwining map $\boxtimes_{P(z_{2})}$
of the type ${W_{2}\boxtimes_{P(z_{2})}W_{3}}\choose {W_{2}
\quad\quad\quad W_{3}}$
gives a linear map
$\Psi^{(1)}_{P(z_{1}, z_{2})}$ {}from
$$W_{1}\hboxtr_{P(z_{1})}(W_{2}\boxtimes_{P(z_{2})}W_{3})
=(W_{1}\boxtimes_{P(z_{1})}(W_{2}\boxtimes_{P(z_{2})}W_{3}))'$$
to $(W_{1}\otimes W_{2}\otimes W_{3})^{*}$ such that
\begin{eqnarray}
\lefteqn{\Psi^{(1)}_{P(z_{1}, z_{2})}(\nu)(w_{(1)}\otimes w_{(2)}\otimes
w_{(3)})}\nno\\
&&=\langle \nu, w_{(1)}\boxtimes_{P(z_{1})}(w_{(2)}
\boxtimes_{P(z_{2})}w_{(3)})\rangle_{W_{1}\boxtimes_{P(z_{1})}(W_{2}
\boxtimes_{P(z_{2})}W_{3})}
\end{eqnarray}
for all $\nu \in W_{1}\hboxtr_{P(z_{1})}(W_{2}\boxtimes_{P(z_{2})}W_{3})$,
$w_{(1)}\in W_{1}$, $w_{(2)}\in W_{2}$ and $w_{(3)}\in W_{3}$,
where $w_{(1)}\boxtimes_{P(z_{1})}(w_{(2)}
\boxtimes_{P(z_{2})}w_{(3)})$  is the image of
$w_{(1)}\otimes w_{(2)}\otimes w_{(3)}$ under the product of
$\boxtimes_{P(z_{1})}$ and $\boxtimes_{P(z_{2})}$. Since the images of
elements of $W_{1}\hboxtr_{P(z_{1})}(W_{2}\boxtimes_{P(z_{2})}W_{3})$
under $\Psi^{(1)}_{P(z_{1}, z_{2})}$ satisfy the the $P(z_{1},
z_{2})$-compatibility condition and the $P(z_{1}, z_{2})$-, the
$P(z_{2})$- and the $P(z_{1}-z_{2})$-local grading-restriction
conditions, $\Psi^{(1)}_{P(z_{1}, z_{2})}$ is in fact a map to
$W_{P(z_{1}, z_{2})}$. {}From the definitions of $L'_{P(z_{1},
z_{2})}(0)$ and of $Y'_{P(z_{1}, z_{2})}$, we see that
$\Psi^{(1)}_{P(z_{1}, z_{2})}$ is also a graded linear map and we have
\begin{equation}
\Psi^{(1)}_{P(z_{1}, z_{2})}Y'_{W_{1}\hboxtr_{P(z_{1})}
(W_{2}\boxtimes_{P(z_{2})}W_{3})}(v, x)
(\Psi^{(1)}_{P(z_{1}, z_{2})})^{-1}
=Y'_{P(z_{1}, z_{2})}(v, x)
\end{equation}
for all $v\in V$. Thus the image of
$\Psi^{(1)}_{P(z_{1}, z_{2})}$ equipped with the vertex operator map
$Y'_{P(z_{1}, z_{2})}$ is a $V$-module and
$\Psi^{(1)}_{P(z_{1}, z_{2})}$ is a
homomorphism of $V$-modules
{}from
$W_{1}\hboxtr_{P(z_{1})}(W_{2}\boxtimes_{P(z_{2})}W_{3})$ to this $V$-module.
We now want to show that the image of $\Psi^{(1)}_{P(z_{1}, z_{2})}$
is equal to
$W_{P(z_{1}, z_{2})}$ and
$\Psi^{(1)}_{P(z_{1}, z_{2})}$ is an isomorphism of $V$-modules {}from
$W_{1}\hboxtr_{P(z_{1})}(W_{2}\boxtimes_{P(z_{2})}W_{3})$ to
$W_{P(z_{1}, z_{2})}$.

We shall show this by constructing a linear map
$(\Psi^{(1)}_{P(z_{1}, z_{2})})^{-1}$ {}from
$W_{P(z_{1}, z_{2})}$ to
$W_{1}\hboxtr_{P(z_{1})}(W_{2}\boxtimes_{P(z_{2})}W_{3})$
and showing that it is the inverse of $\Psi^{(1)}_{P(z_{1}, z_{2})}$.
For any $w_{(1)}\in W_{1}$ and any $\lambda\in W_{P(z_{1}, z_{2})}$,
$\mu^{(1)}_{\lambda, w_{(1)}}$ is in $\overline{W_{2}\hboxtr_{P(z_{2})}W_{3}}$
by the
$P(z_{1}, z_{2})$-compatibility condition, the $P(z_{1}, z_{2})$-local
grading-restriction condition and the $P(z_{2})$-local
grading-restriction condition. We define an element
$(\Psi^{(1)}_{P(z_{1}, z_{2})})^{-1}(\lambda)\in
(W_1\otimes(W_2\boxtimes_{P(z_2)}W_3))^*$ by
\begin{equation}
(\Psi^{(1)}_{P(z_{1}, z_{2})})^{-1}(\lambda)(w_{(1)}\otimes w)=\bra w,
\mu^{(1)}_{\lambda,
w_{(1)}}\ket_{W_{2}\shboxtr_{P(z_{2})}W_{3}}
\end{equation}
for all $w_{(1)}\in W_{1}$ and $w\in W_{2}\boxtimes_{P(z_{2})}W_{3}$. By
the $P(z_{1}, z_{2})$-compatibility condition and the $P(z_{1}, z_{2})$-local
grading-restriction condition for $\lambda$,
$(\Psi^{(1)}_{P(z_{1}, z_{2})})^{-1}(\lambda)$ is in
fact in $W_1\hboxtr_{P(z_{1})} (W_2\boxtimes_{P(z_2)}W_3)$.
Thus $\lambda\mapsto(\Psi^{(1)}_{P(z_{1}, z_{2})})^{-1}(\lambda)$
defines a map $(\Psi^{(1)}_{P(z_{1}, z_{2})})^{-1}$ {}from
$W_{P(z_{1}, z_{2})}$ to $W_1\hboxtr_{P(z_{1})}
(W_2\boxtimes_{P(z_2)}W_3)$.

For any $\lambda\in W_{P(z_{1}, z_{2})}$, using the definitions of
$\Psi^{(1)}_{P(z_{1}, z_{2})}$,
$(\Psi^{(1)}_{P(z_{1}, z_{2})})^{-1}$, $\boxtimes_{P(z_{1})}$ and
$\hboxtr_{P(z_{1})}$, we have
\begin{eqnarray}
\lefteqn{\Psi^{(1)}_{P(z_{1}, z_{2})}
((\Psi^{(1)}_{P(z_{1}, z_{2})})^{-1}(\lambda))(w_{(1)}\otimes
w_{(2)}\otimes w_{(3)})}\nno\\
&&=\sum _{n\in \Bbb{C}}\langle (\Psi^{(1)}_{P(z_{1}, z_{2})})^{-1}(\lambda),
w_{(1)}
\boxtimes_{P(z_{1})}
P_{n}(w_{(2)}\boxtimes_{P(z_{2})} w_{(3)})\rangle_{W_{1}\boxtimes_{P(z_{1})}
(W_{2}\boxtimes_{P(z_{2})}W_{3})} \nno\\
&&=\sum _{n\in \Bbb{C}}\hboxtr_{P(z_{1})}
((\Psi^{(1)}_{P(z_{1}, z_{2})})^{-1}(\lambda))(w_{(1)}
\otimes P_{n}(w_{(2)}\boxtimes_{P(z_{2})} w_{(3)}))  \nno\\
&&=\sum _{n\in \Bbb{C}}((\Psi^{(1)}_{P(z_{1}, z_{2})})^{-1}(\lambda))(w_{(1)}
\otimes P_{n}(w_{(2)}\boxtimes_{P(z_{2})} w_{(3)}))  \nno\\
&&=\sum _{n\in \Bbb{C}}\langle
P_{n}(w_{(2)}\boxtimes_{P(z_{2})} w_{(3)}), \mu^{(1)}_{\lambda, w_{(1)}}
 \rangle_{W_{2}\shboxtr_{P(z_{2})}W_{3}}\nno\\
&&=\sum _{n\in \Bbb{C}}\langle P_{n}(\mu^{(1)}_{\lambda, w_{(1)}}),
w_{(2)}\boxtimes_{P(z_{2})}
w_{(3)}\rangle_{W_{2}\boxtimes_{P(z_{2})}W_{3}}\nno\\
&&=\sum _{n\in \Bbb{C}}(\hboxtr_{P(z_{2})}(P_{n}(\mu^{(1)}_{\lambda,
w_{(1)}})))
(w_{(2)}\otimes w_{(3)})\nno\\
&&=\sum _{n\in \Bbb{C}}(P_{n}(\mu^{(1)}_{\lambda, w_{(1)}}))
(w_{(2)}\otimes w_{(3)})\nno\\
&&=\mu^{(1)}_{\lambda, w_{(1)}}
(w_{(2)}\otimes w_{(3)})\nno\\
&&=\lambda(w_{(1)}\otimes w_{(2)}\otimes w_{(3)}),
\end{eqnarray}
proving
\begin{equation}
\Psi^{(1)}_{P(z_{1}, z_{2})}(\Psi^{(1)}_{P(z_{1}, z_{2})})^{-1}=1.
\end{equation}

Next we want to show that $\Psi^{(1)}_{P(z_{1}, z_{2})}$ is injective.
Let $\nu_{1},
\nu_{2}$ be two elements of $W_1\hboxtr_{P(z_{1})} (W_2\boxtimes_{P(z_2)}W_3)$
such that
\begin{equation}
\Psi^{(1)}_{P(z_{1}, z_{2})}(\nu_{1})=\Psi^{(1)}_{P(z_{1}, z_{2})}(\nu_{2}),
\end{equation}
that is,
\begin{eqnarray}
\lefteqn{\langle \nu_{1}, w_{(1)}\boxtimes_{P(z_{1})}(w_{(2)}
\boxtimes_{P(z_{2})}w_{(3)})\rangle_{W_1\boxtimes_{P(z_{1})}
(W_2\boxtimes_{P(z_2)}W_3)}}\nno\\
&&=\langle \nu_{2}, w_{(1)}\boxtimes_{P(z_{1})}(w_{(2)}
\boxtimes_{P(z_{2})}w_{(3)})\rangle_{W_1\boxtimes_{P(z_{1})}
(W_2\boxtimes_{P(z_2)}W_3)}
\end{eqnarray}
for all $w_{(1)}\in W_{1}$, $w_{(2)}\in W_{2}$ and $w_{(3)}\in W_{3}$.
Let ${\cal Y}_{1}$ and ${\cal Y}_{2}$ be the intertwining operators of type
${W_{1}\boxtimes_{P(z_{1})} (W_{2}\boxtimes_{P(z_{2})} W_{3})}
\choose {W_{1}\; W_{2}\boxtimes_{P(z_{2})} W_{3}}$
and ${W_{2}\boxtimes_{P(z_{2})} W_{3}}\choose {W_{2}\;W_{3}}$ corresponding to
the intertwining maps $\boxtimes_{P(z_{1})}$ and $\boxtimes_{P(z_{2})}$,
respectively. Then (14.69) becomes
\begin{eqnarray}
\lefteqn{\langle \nu_{1}, {\cal Y}_{1}(w_{(1)}, x_{1}){\cal Y}_{2}(w_{(2)},
x_{2})
w_{(3)}\rangle_{W_1\boxtimes_{P(z_{1})}
(W_2\boxtimes_{P(z_2)}W_3)}\lbar_{x_{1}=z_{1}, x_{2}
= z_{2}}}\nno\\
&&=\langle \nu_{2}, {\cal Y}_{1}(w_{(1)}, x_{1})\cdot\nno\\
&&\hspace{4em}\cdot{\cal Y}_{2}(w_{(2)}, x_{2})
w_{(3)}\rangle_{W_1\boxtimes_{P(z_{1})}
(W_2\boxtimes_{P(z_2)}W_3)}\lbar_{x_{1}=z_{1}, x_{2}
= z_{2}}.
\end{eqnarray}
By (14.39), both sides of (14.70), or equivalently
of (14.69), can be expanded as series in $e^{\log z_{2}}$ of
the form
\begin{equation}
\sum_{i\in {\Bbb Z}_{+}}a_{i}(z_{1})e^{m_{i}\log z_{2}}
\end{equation}
where $\{m_{i}\}_{i\in {\Bbb Z}_{+}}$ is a
sequence with of strictly increasing real numbers.
For a multi-valued function $f(z)$ which can be expanded as
$\sum_{i\in {\Bbb Z}_{+}}a_{i}e^{m_{i}\log z}$ in a
neighborhood of $0$ with $0$ deleted, we define
$\res_{z}f(z)$ to be $a_{i_{0}}$ if there is a positive integer $i_{0}$ such
that $m_{i_{0}}=-1$ and to be $0$ if there is no such $i_{0}$.
Lemma 14.5 implies that $\res_{z}f(z)$ is well
defined when $\{m_{l}\}_{l\in {\Bbb Z}_{+}}$
is a sequence of strictly increasing real numbers.

Let $f_{l}(z)$ and $f_{r}(z)$ be the multi-valued functions of $z$ obtained
{}from the left-hand side and the right-hand side, repectively, of (14.70)
by replacing
$z_{2}$ be $z$.
{}From the assumption that the product of a
$P(z_{1})$-intertwining map and a $P(z)$-intertwining map is convergent
when $0<|z|<|z_{1}|$, we see that  $f_{l}(z)$ and $f_{r}(z)$ are
well-deined multi-valued functions
of $z$ when $0<|z|<|z_{1}|$.
By the $L(-1)$-derivative property for intertwining
operators, we see that $f_{l}(z)$ and $f_{r}(z)$ are
 analytic when $0<|z|<|z_{1}|$.

Now we can apply Lemma 14.5  to $f_{l}(z)$ and $f_{r}(z)$ since both can be
expanded in the form of (14.71) with $z_{2}$ repalced by $z$.
By this lemma, the expansion coefficients
$\res_{z}e^{m\log z}f_{l}(z)$ and $\res_{z}e^{m\log z}f_{l}(z)$,
$m\in {\Bbb C}$, of
$f_{l}(z)$ and $f_{r}(z)$
 are determined uniquely by $f_{l}(z)$ and $f_{r}(z)$,
respectively. But $f_{l}(z)$ and $f_{r}(z)$ are equal. So their expansion
coefficients are also equal. Writing down these expansion coefficients
explicitly, we obtain
\begin{eqnarray}
\lefteqn{\langle \nu_{1}, {\cal Y}_{1}(w_{(1)},
x_{1})(w_{(2)})_{m}
w_{(3)}\rangle_{W_1\boxtimes_{P(z_{1})}
(W_2\boxtimes_{P(z_2)}W_3)}\lbar_{x^{n}_{1}=e^{n\log z_{1}}, n\in {\Bbb
C}}}\nno\\
&&=\langle \nu_{2}, {\cal Y}_{1}(w_{(1)},
x_{1})(w_{(2)})_{m}
w_{(3)}\rangle_{W_1\boxtimes_{P(z_{1})}
(W_2\boxtimes_{P(z_2)}W_3)}\lbar_{x^{n}_{1}=e^{n\log z_{1}}, n\in {\Bbb C}}
\;\;\;\;\;\;\;
\end{eqnarray}
or equivalently
\begin{eqnarray}
\lefteqn{\langle \nu_{1}, w_{(1)}\boxtimes_{P(z_{1})}(w_{(2)})_{m}
w_{(3)})\rangle_{W_1\boxtimes_{P(z_{1})}
(W_2\boxtimes_{P(z_2)}W_3)}}\nno\\
&&=\langle \nu_{2}, w_{(1)}\boxtimes_{P(z_{1})}(w_{(2)})_{m}
w_{(3)})\rangle_{W_1\boxtimes_{P(z_{1})}
(W_2\boxtimes_{P(z_2)}W_3)}
\end{eqnarray}
for all $m\in {\Bbb C}$ and all $w_{(1)}\in W_{1}$, $w_{(2)}\in W_{2}$ and
$w_{(3)}\in W_{3}$. The conclusion $\nu_{1}=\nu_{2}$ which we need
follows {}from the result below:

\begin{lemma}
The module $W_{2}\boxtimes_{P(z_{2})}W_{3}$ is spanned by the
homogeneous components
 of the elements of $\overline{W_{2}\boxtimes_{P(z_{2})}W_{3}}$ of
the form $w_{(2)}\boxtimes_{P(z_{2})} w_{(3)}$, for all $w_{(2)}\in W_{2}$
and $w_{(3)}\in W_{3}$.
\end{lemma}

We shall prove this in the next section.

We have shown that $\Psi^{(1)}_{P(z_{1}, z_{2})}$ is injective.
Combining this result with
(14.53), we see that $\Psi^{(1)}_{P(z_{1}, z_{2})}$ is an isomorphism of
$V$-modules {}from
$$W_{1}\hboxtr_{P(z_{1})}(W_{2}\boxtimes_{P(z_{2})}W_{3})$$ to
$W_{P(z_{1}, z_{2})}$.

Similarly,
the iterate of the $P(z_{1}-z_{2})$-intertwining map
$\boxtimes_{P(z_{1}-z_{2})}$ of type
$${(W_{1}\boxtimes_{P(z_{1}-z_{2})}W_{2})
\boxtimes_{P(z_{2})}W_{3}}\choose {W_{1}\boxtimes_{P(z_{1}-z_{2})}W_{2}
\quad\quad\quad
W_{3}}$$
and the $P(z_{4})$-intertwining map $\boxtimes_{P(z_{4})}$
of type ${W_{1}\boxtimes_{P(z_{1}-z_{2})}W_{2}}\choose {W_{1}
\quad\quad\quad W_{2}}$
gives a linear map
$\Psi^{(2)}_{P(z_{1}, z_{2})}$ {}from
$$(W_{1}\boxtimes_{P(z_{1}-z_{2})}W_{2})\hboxtr_{P(z_{2})}W_{3}
=((W_{1}\boxtimes_{P(z_{1}-z_{2})}W_{2})\boxtimes_{P(z_{2})}W_{3})'$$
to $(W_{1}\otimes W_{2}\otimes W_{3})^{*}$ such that
\begin{eqnarray}
\lefteqn{\Psi^{(2)}_{P(z_{1}, z_{2})}(\nu)(w_{(1)}\otimes w_{(2)}
\otimes w_{(3)})=}\nno\\
&&=\langle \nu, (w_{(1)}\boxtimes_{P(z_{1}-z_{2})}w_{(2)})
\boxtimes_{P(z_{2})}w_{(3)}\rangle_{(W_{1}\boxtimes_{P(z_{1}-z_{2})}W_{2})
\boxtimes_{P(z_{2})}W_{3}}
\end{eqnarray}
for all
$$\nu \in (W_{1}\boxtimes_{P(z_{1}z_{2})}W_{2})\hboxtr_{P(z_{2})}W_{3},$$
$w_{(1)}\in W_{1}$, $w_{(2)}\in W_{2}$ and $w_{(3)}\in W_{3}$, where
$(w_{(1)}\boxtimes_{P(z_{1}-z_{2})}w_{(2)})
\boxtimes_{P(z_{2})}w_{(3)}$ is the image of $w_{(1)}\otimes
w_{(2)}\otimes w_{(3)}$ under the iterate of $\boxtimes_{P(z_{1}-z_{2})}$
and $\boxtimes_{P(z_{2})}$.
We can show similarly that this map
is in fact an isomorphism of $V$-modules {}from
$$(W_{1}\boxtimes_{P(z_{1}-z_{2})}W_{2})\hboxtr_{P(z_{2})}W_{3}$$ to
$W_{P(z_{1}, z_{2})}$.

Let
$${\cal A}_{P(z_{1}), P(z_{2})}
^{P(z_{1}-z_{2}), P(z_{2})}: W_{1}\boxtimes_{P(z_{1})}(W_{2}
\boxtimes_{P(z_{2})}W_{3})\to (W_{1}\boxtimes_{P(z_{1}-z_{2})}W_{2})
\boxtimes_{P(z_{2})}W_{3}$$ be the graded adjoint map of
$(\Psi^{(1)}_{P(z_{1}, z_{2})})^{-1}\circ
\Psi^{(2)}_{P(z_{1}, z_{2})}$.
Then ${\cal A}_{P(z_{1}), P(z_{2})}
^{P(z_{1}-z_{2}), P(z_{2})}$ is an isomorphism and is called the
{\it associativity isomorphism for  $W_{1}$, $W_{2}$ and $W_{3}$
associated with $(P(z_{1}), P(z_{2}); P(z_{1}-z_{2}), P(z_{2}))$}. {}From
(14.63) and (14.74), we obtain
\begin{eqnarray}
\lefteqn{\langle \nu, \overline{{\cal A}}_{P(z_{1}),
P(z_{2})}^{P(z_{1}-z_{2}),
P(z_{2})}(w_{(1)}\boxtimes_{P(z_{1})}(w_{(2)}
\boxtimes_{P(z_{2})}w_{(3)}))\rangle_{W_{1}\boxtimes_{P(z_{1})}(W_{2}
\boxtimes_{P(z_{2})}W_{3})}}\nno\\
&&=\langle \nu, (w_{(1)}\boxtimes_{P(z_{1}-z_{2})}w_{(2)})
\boxtimes_{P(z_{2})}w_{(3)}\rangle_{(W_1 \boxtimes_{P(z_1-z_2)}W_2)
\boxtimes_{P(z_2)}W_3}\;\;\;\;\;\;\;\;\;
\end{eqnarray}
for all $\nu\in (W_1 \boxtimes_{P(z_1-z_2)}W_2)\hboxtr_{P(z_2)}W_3$,
$w_{(1)}\in W_{1}$, $w_{(2)}\in W_{2}$ and $w_{(3)}\in W_{3}$,
where
$$\overline{{\cal A}}_{P(z_{1}), P(z_{2})}^{P(z_{1}-z_{2}),
P(z_{2})}: \overline{W_1\boxtimes_{P(z_1)}(W_2 \boxtimes_{P(z_2)}W_3)}
\to \overline{(W_1 \boxtimes_{P(z_1-z_2)}W_2)\boxtimes_{P(z_2)}W_3}$$
is the unique extension of ${\cal A}_{P(z_{1}),
P(z_{2})}^{P(z_{1}-z_{2}), P(z_{2})}$. Equivalently, we obtain
\begin{eqnarray}
\lefteqn{\overline{{\cal A}}_{P(z_{1}), P(z_{2})}^{P(z_{1}-z_{2}),
P(z_{2})}(w_{(1)}\boxtimes_{P(z_{1})}(w_{(2)}
\boxtimes_{P(z_{2})}w_{(3)}))}\nno\\
&&=(w_{(1)}\boxtimes_{P(z_{1}-z_{2})}w_{(2)})
\boxtimes_{P(z_{2})}w_{(3)}
\end{eqnarray}
for all $w_{(1)}\in W_{1}$, $w_{(2)}\in W_{2}$ and $w_{(3)}\in W_{3}$.

We summarize what we obtained above in half of the following theorem:

\begin{theo}
Assume that $V$ is a rational vertex operator algebra, that all irreducible
$V$-modules are ${\Bbb R}$-greded and that $V$  has one
of the properties in Proposition 14.1. If $V$
 also has either one of the properties
in Propositon 14.7, then for any $V$-module $W_{1}$, $W_{2}$ and $W_{3}$
and any complex numbers $z_{1}$ and $z_{2}$ satisfying (14.7),
there exists a unique isomorphism
${\cal A}_{P(z_{1}), P(z_{2})}^{P(z_{1}-z_{2}), P(z_{2})}$ {}from
$W_1\boxtimes_{P(z_1)}(W_2 \boxtimes_{P(z_2)}W_3)$
to $(W_1 \boxtimes_{P(z_1-z_2)}W_2)\boxtimes_{P(z_2)}W_3$ such that
(14.76) holds for all $w_{(1)}\in W_{1}$, $w_{(2)}\in W_{2}$
and $w_{(3)}\in W_{3}$. Conversely, if there exists such an isomorphism
for any $V$-module $W_{1}$, $W_{2}$ and $W_{3}$
and any complex numbers $z_{1}$ and $z_{2}$ satisfying (14.7), then $V$
has both poperties in Proposition 14.7.
\end{theo}
\pf
We need only to show that if there exists such an isomorphism for
any $V$-module $W_{1}$, $W_{2}$ and $W_{3}$
and any complex numbers $z_{1}$ and $z_{2}$ satisfying (14.7), $V$ has the
first property in Proposition 14.7.

 Let $W_{P(z_{1}), P(z_{2})}$ be the
subspace of $(W_{1}\otimes W_{2}\otimes W_{3})^{*}$ consisting of elements
satisfying the $P(z_{1}, z_{2})$-compatibility condition, the
$P(z_{1}, z_{2})$-local grading-restriction condition and the
$P(z_{2})$-local grading-restriction condition and
$W^{P(z_{1}-z_{2}), P(z_{2})}$ the
subspace of $(W_{1}\otimes W_{2}\otimes W_{3})^{*}$ consisting of elements
satisfying the $P(z_{1}, z_{2})$-compatibility condition, the
$P(z_{1}, z_{2})$-local grading-restriction condition and the
$P(z_{1}-z_{2})$-local grading-restriction condition. Since
products and iterates of intertwining maps give elements
in $W_{P(z_{1}), P(z_{2})}$ and $W^{P(z_{1}-z_{2}), P(z_{2})}$,
respectively, we need  only to prove
\begin{equation}
W_{P(z_{1}), P(z_{2})}=W^{P(z_{1}-z_{2}), P(z_{2})}(=W_{P(z_{1}, z_{2})}).
\end{equation}
The maps $\Psi^{(1)}_{P(z_{1}, z_{2})}$ and
$(\Psi^{(1)}_{P(z_{1}, z_{2})})^{-1}$ in this case are still well defined and
are  maps {}from $$W_{1}\hboxtr_{P(z_{1})}
(W_{2}\boxtimes_{P(z_{2})}W_{3})$$
 to $W_{P(z_{1}), P(z_{2})}$ and {}from
$W_{P(z_{1}), P(z_{2})}$ to $$W_{1}\hboxtr_{P(z_{1})}
(W_{2}\boxtimes_{P(z_{2})}W_{3}),$$ respectively, and they are
still the inverses of each other. Similarly we have the maps
$\Psi^{(2)}_{P(z_{1}, z_{2})}$ and
$(\Psi^{(2)}_{P(z_{1}, z_{2})})^{-1}$ {}from $$(W_{1}\boxtimes_{P(z_{1}-z_{2})}
W_{2})\hboxtr_{P(z_{2})}W_{3})$$
 to $W^{P(z_{1}-z_{2}), P(z_{2})}$ and {}from
$W^{P(z_{1}-z_{2}), P(z_{2})}$ to $$(W_{1}\boxtimes_{P(z_{1}-z_{2})}
W_{2})\hboxtr_{P(z_{2})}W_{3}),$$ respectively, and they are the inverses
of each other. Let
$$({\cal A}_{P(z_{1}), P(z_{2})}^{P(z_{1}-z_{2}), P(z_{2})})':
(W_{1}\boxtimes_{P(z_{1}-z_{2})}
W_{2})\hboxtr_{P(z_{2})}W_{3}) \to W_{1}\hboxtr_{P(z_{1})}
(W_{2}\boxtimes_{P(z_{2})}W_{3})$$ be the graded
adjoint of ${\cal A}_{P(z_{1}), P(z_{2})}^{P(z_{1}-z_{2}), P(z_{2})}$.
Then
$$\Psi^{(1)}_{P(z_{1}, z_{2})}\circ
({\cal A}_{P(z_{1}), P(z_{2})}^{P(z_{1}-z_{2}), P(z_{2})})'
\circ(\Psi^{(2)}_{P(z_{1}, z_{2})})^{-1}$$
is an isomorphism {}from $W^{P(z_{1}-z_{2}), P(z_{2})}$ to
$W_{P(z_{1}), P(z_{2})}$. But on the other hand, by the definitions of
$(\Psi^{(2)}_{P(z_{1}, z_{2})})^{-1}$ and $\Psi^{(1)}_{P(z_{1}, z_{2})}$
and (14.75), for any $\lambda\in W^{P(z_{1}-z_{2}), P(z_{2})}$,
$$\Psi^{(1)}_{P(z_{1}, z_{2})}\circ
({\cal A}_{P(z_{1}), P(z_{2})}^{P(z_{1}-z_{2}), P(z_{2})})'
\circ(\Psi^{(2)}_{P(z_{1}, z_{2})})^{-1}(\lambda)$$
as an element of $(W_{1}\otimes W_{2}\otimes W_{3})^{*}$ is equal to
$\lambda$ itself, Thus we have (14.77).
\epfv

By Proposition 14.8 and this theorem, we obtain:

\begin{theo}
Assume that $V$ is a rational vertex operator algebra, that all irreducible
$V$-modules are ${\Bbb R}$-greded and that $V$  has one
of the properties in Proposition 14.1. Then for any $V$-module $W_{1}$,
$W_{2}$ and $W_{3}$
and any complex numbers $z_{1}$ and $z_{2}$ satisfying (14.7),
there exists a unique isomorphism
${\cal A}_{P(z_{1}), P(z_{2})}^{P(z_{1}-z_{2}), P(z_{2})}$ {}from
$W_1\boxtimes_{P(z_1)}(W_2 \boxtimes_{P(z_2)}W_3)$
to $(W_1 \boxtimes_{P(z_1-z_2)}W_2)\boxtimes_{P(z_2)}W_3$ such that
(14.76) holds for all $w_{(1)}\in W_{1}$, $w_{(2)}\in W_{2}$
and $w_{(3)}\in W_{3}$ if and only if
the intertwining operators for $V$ have the associativity properties
in Proposition 14.8.
\end{theo}

\renewcommand{\theequation}{\thesection.\arabic{equation}}
\renewcommand{\therema}{\thesection.\arabic{rema}}
\setcounter{equation}{0}
\setcounter{rema}{0}

\section{Proofs of the lemmas in Section 14}

\subsection{Proof of Lemma 14.2}

We prove only (14.9). The Proofs of (14.8), (14.10) and (14.11) are similar.
In the
process of the proof,
we shall see that both sides of (14.9) exist and their coefficients
are in fact expansions of rational
functions as indicated.

By the definition of the $\delta$-function, the right-hand side of  (14.9) is
equal to
\begin{eqnarray}
\lefteqn{\sum_{n, m\in {\Bbb Z}}z_{2}^{-n-1}(x_{0}^{-1}-x_{2})^{n}
(z_{1}-z_{2})^{-m-1}(x_{2}-x_{1})^{m}}\nno\\
&&= \sum_{n, m\in {\Bbb Z}}\sum_{k, l\in {\Bbb N}}
{{n} \choose {k}} {{m} \choose {l}}
z_{2}^{-n-1}(z_{1}-z_{2})^{-m-1}x_{0}^{-n+k}(-1)^{k}
x_{2}^{k}x_{2}^{m-l}(-1)^{l}x_{1}^{l}\nno\\
&&=\sum_{n, m\in {\Bbb Z}}\sum_{k, l\in {\Bbb N}}
{{n} \choose {k}} {{m} \choose {l}}(-1)^{k+l}
z_{2}^{-n-1}(z_{1}-z_{2})^{-m-1}x_{0}^{-n+k}
x_{1}^{l}x_{2}^{k+m-l}.\nno\\
&&
\end{eqnarray}
The coefficient of the right-hand side of (15.1) in
$x_{0}^{r}x_{1}^{s}x_{2}^{t}$ is
\begin{eqnarray}
\lefteqn{\sum_{k\in {\Bbb N}}{{k-r} \choose {k}} {{t-k+s} \choose
{s}}(-1)^{k+s} z_{2}^{-k+r-1}(z_{1}-z_{2})^{-t+k-s-1}}\nno\\
&&=\sum_{k\in {\Bbb N}}{{r-1} \choose {k}} {{k-t-1} \choose {s}}
z_{2}^{-k+r-1}(z_{1}-z_{2})^{-t+k-s-1}\nno\\
&&=\sum_{k\in {\Bbb N}}{{r-1} \choose {k}}
z_{2}^{-k+r-1}\frac{1}{s!}\frac{d^{s}}{dz_{1}^{s}}(z_{1}-z_{2})^{-t+k-1}\nno\\
&&=\frac{1}{s!}\frac{d^{s}}{dz_{1}^{s}}\biggl((z_{1}-z_{2})^{-t-1}\sum_{k\in
{\Bbb N}}{{r-1} \choose {k}}  z_{2}^{-k+r-1}(z_{1}-z_{2})^{k}\biggr).
\end{eqnarray}
Since we
have the inequality $|z_{2}|>|z_{1}-z_{2}|$, the right-hand side of (15.2)
is absolutely convergent to
\begin{equation}
\frac{1}{s!}\frac{d^{s}}{dz_{1}^{s}}\left((z_{1}-z_{2})^{-t-1}
z_{1}^{r-1}\right).
\end{equation}

On the other hand, the left-hand side of (14.9) is equal to
\begin{eqnarray}
\lefteqn{\sum_{n, m\in {\Bbb Z}}z_{1}^{-n-1}(x_{0}^{-1}-x_{1})^{n}x_{2}^{-m-1}
(x_{0}^{-1}-z_{2})^{m}}\nno\\
&&=\sum_{n, m\in {\Bbb Z}}\sum_{k, l\in {\Bbb N}}{{n}\choose {k}}{{m}\choose
{l}} z_{1}^{-n-1}x_{0}^{-n+k}(-1)^{k}x_{1}^{k}x_{2}^{-m-1}
x_{0}^{-m+l}(-1)^{l}z_{2}^{l}\nno\\
&&=\sum_{n, m\in {\Bbb Z}}\sum_{k, l\in {\Bbb N}}{{n}\choose {k}}{{m}\choose
{l}}(-1)^{k+l} z_{1}^{-n-1}z_{2}^{l}x_{0}^{-n+k-m+l}x_{1}^{k}x_{2}^{-m-1}.
\end{eqnarray}
The coefficient of the right-hand side of (15.4) in
$x_{0}^{r}x_{1}^{s}x_{2}^{t}$ is
\begin{eqnarray}
\lefteqn{\sum_{l\in {\Bbb N}}{{-r+s+t+1+l}\choose {s}}{{-t-1}\choose
{l}}(-1)^{s+l} z_{1}^{r-s-t-l-2}z_{2}^{l}}\nno\\
&&=\sum_{l\in {\Bbb N}}\frac{1}{s!}\frac{d^{s}}{dz_{1}^{s}}{{-t-1}\choose
{l}}(-1)^{l} z_{1}^{r-t-l-2}z_{2}^{l}\nno\\
&&=\frac{1}{s!}\frac{d^{s}}{dz_{1}^{s}}\biggl(z_{1}^{r-1}\sum_{l\in {\Bbb N}}
{{-t-1}\choose {l}}(-1)^{l} z_{1}^{-t-1-l}z_{2}^{l}\biggr).
\end{eqnarray}
Since we also have the inequality $|z_{1}|>|z_{2}|$, the right-hand side of
(15.5) is absolutely convergent to (15.3), proving (14.9).

\subsection{Proof of Lemma 14.3}

We prove only (14.36). The definition (14.18) can be written as
\begin{eqnarray}
\lefteqn{z^{-1}_2\delta\left(\frac{x^{-1}_0-x_2}{z_2}\right) (z_1-z_{2})^{-1}
\delta\left(\frac{x_2-x_1}{z_1-z_{2}}\right)\cdot}\nno\\
&&\hspace{4em}\cdot
\lambda(Y_{1}(e^{x_0L(1)}(-x^2_0)^{-
L(0)}v, x_1) w_{(1)}\otimes w_{(2)}\otimes w_{(3)})\nno\\
&&\quad + z^{-1}_2\delta\left(\frac{x^{-1}_0-x_2}{z_2}\right)x^{-1}_1
\delta\left(\frac{z_1-z_{2}-x_2}{x_1}\right)\cdot\nno\\
&&\hspace{4em}\cdot
\lambda(w_{(1)}\otimes Y_{2}(e^{x_0L(1)}(-
x^2_0)^{-L(0)}v, x_2) w_{(2)}\otimes w_{(3)})\nno\\
&&=\biggl(\tau^{(2)}_{P(z_{1}, z_{2})}
\biggl(x^{-1}_2\delta\left(\frac{x^{-1}_0-z_2}{x_2}\right)\cdot\nno\\
&&\hspace{4em}\cdot x^{-1}_1\delta\left(\frac{x_2-(z_1-z_{2})}{x_1}\right)
Y_{t}(v, x_0)\biggr)\lambda
\biggr)
(w_{(1)}\otimes w_{(2)}\otimes w_{(3)})\nno\\
&&\quad +x^{-1}_2\delta\left(\frac{z_2-x^{-1}_0}{x_2}\right)
x^{-1}_1\delta\left(\frac{x_2-(z_1-z_{2})}{x_1}\right)\cdot\nno\\
&&\hspace{4em}\cdot
\lambda(w_{(1)}\otimes w_{(2)}\otimes
Y_{3}(e^{x_0L(1)}(-x^2_0)^{-L(0)}v, x^{-1}_0) w_{(3)})
\end{eqnarray}
On the both sides of (15.6),
first replacing $v$  by
$$(-x^2_0)^{L(0)}e^{-x_0L(1)}e^{x^{-1}_2L(1)}(-x^{2}_2)^{L(0)}v$$
 and then taking
$\res_{x_{0}^{-1}}$, we obtain
\begin{eqnarray}
\lefteqn{(z_1-z_{2})^{-1}
\delta\left(\frac{x_2-x_1}{z_1-z_{2}}\right)\cdot}\nno\\
&&\hspace{2em}\cdot
\lambda(Y_{1}(e^{x^{-1}_2L(1)}(-x^{2}_2)^{L(0)}v, x_1) w_{(1)}
\otimes w_{(2)}\otimes w_{(3)})\nno\\
&&\quad + x^{-1}_1
\delta\left(\frac{z_1-z_{2}-x_2}{x_1}\right)\cdot\nno\\
&&\hspace{2em}\cdot
\lambda(w_{(1)}\otimes Y_{2}(e^{x^{-1}_2L(1)}(-x^{2}_2)^{L(0)}v,
x_2) w_{(2)}\otimes w_{(3)})\nno\\
&&=\res_{x_{0}^{-1}}\biggl(\tau^{(2)}_{P(z_{1}, z_{2})}
\biggl(x^{-1}_2\delta\left(\frac{x^{-1}_0-z_2}{x_2}\right)
 x^{-1}_1\delta\left(\frac{x_2-(z_1-z_{2})}{x_1}\right)\cdot\nno\\
&&\hspace{2em}\cdot
Y_{t}((-x^2_0)^{L(0)}e^{-x_0L(1)}e^{x^{-1}_2L(1)}(-x^{2}_2)^{L(0)}v,
x_0)\biggr)\lambda\biggr)
(w_{(1)}\otimes w_{(2)}\otimes w_{(3)})\nno\\
&&\quad +\res_{x_{0}^{-1}}x^{-1}_2\delta\left(\frac{z_2-x^{-1}_0}{x_2}\right)
x^{-1}_1\delta\left(\frac{x_2-(z_1-z_{2})}{x_1}\right)\cdot\nno\\
&&\hspace{2em}\cdot
\lambda(w_{(1)}\otimes w_{(2)}\otimes
Y_{3}(e^{x^{-1}_2L(1)}(-x^{2}_2)^{L(0)}v, x^{-1}_0) w_{(3)}).
\end{eqnarray}
By the definition of $\tau_{P(z_{1}-z_{2})}$ and (15.7),
\begin{eqnarray}
\lefteqn{\biggl(\tau_{P(z_{1}-z_{2})}\biggl(x^{-1}_1
\delta\left(\frac{x_2-(z_1-z_{2})}{x_1}\right)Y_{t}(v, x^{-1}_{2})\biggr)
\mu^{(2)}_{\lambda, w_{(3)}}\biggr)(w_{(1)}\otimes w_{(2)})}\nno\\
&&=\res_{x_{0}^{-1}}\biggl(\tau^{(2)}_{P(z_{1}, z_{2})}
\biggl(x^{-1}_2\delta\left(\frac{x^{-1}_0-z_2}{x_2}\right)
 x^{-1}_1\delta\left(\frac{x_2-(z_1-z_{2})}{x_1}\right)\cdot\nno\\
&&\hspace{2em}\cdot
Y_{t}((-x^2_0)^{L(0)}e^{-x_0L(1)}e^{x^{-1}_2L(1)}(-x^{2}_2)^{L(0)}v,
x_0)\biggr)\lambda\biggr)
(w_{(1)}\otimes w_{(2)}\otimes w_{(3)})\nno\\
&&\quad +\res_{x_{0}^{-1}}x^{-1}_2\delta\left(\frac{z_2-x^{-1}_0}{x_2}\right)
x^{-1}_1\delta\left(\frac{x_2-(z_1-z_{2})}{x_1}\right)\cdot\nno\\
&&\hspace{2em}\cdot
\lambda(w_{(1)}\otimes w_{(2)}\otimes
Y_{3}(e^{x^{-1}_2L(1)}(-x^{2}_2)^{L(0)}v, x^{-1}_0) w_{(3)}).
\end{eqnarray}
Note that the right-hand side of (15.8) exists.
Since $\lambda$ satisfies the $P(z_{1}, z_{2})$-compatibility condition, the
right-hand side of (15.8) is equal to
\begin{eqnarray}
\lefteqn{x^{-1}_1\delta\left(\frac{x_2-(z_1-z_{2})}{x_1}\right)
\res_{x_{0}^{-1}} \biggl(\tau^{(2)}_{P(z_{1}, z_{2})}
\biggl(x^{-1}_2\delta\left(\frac{x^{-1}_0-z_2}{x_2}\right)\cdot}\nno\\
&&\hspace{2em}\cdot
Y_{t}((-x^2_0)^{L(0)}e^{-x_0L(1)}e^{x^{-1}_2L(1)}(-x^{2}_2)^{L(0)}v,
x_0)\biggr)\lambda\biggr)
(w_{(1)}\otimes w_{(2)}\otimes w_{(3)})\nno\\
&& +x^{-1}_1\delta\left(\frac{x_2-(z_1-z_{2})}{x_1}\right)
\res_{x_{0}^{-1}}x^{-1}_2\delta\left(\frac{z_2-x^{-1}_0}{x_2}\right)
\cdot\nno\\
&&\hspace{2em}\cdot
\lambda(w_{(1)}\otimes w_{(2)}\otimes
Y_{3}(e^{x^{-1}_2L(1)}(-x^{2}_2)^{L(0)}v, x^{-1}_0) w_{(3)})
\end{eqnarray}
Taking $\res_{x_{1}}$ in (15.9) and then multiplying it by
$x^{-1}_1\delta\left(\frac{x_2-(z_1-z_{2})}{x_1}\right)$, we obtain
(15.9) itself. But by the calculations (15.8)--(15.9),
(15.9) is equal to the left-hand side of
(15.8). So we see that the left-hand side of (15.8) is
equal to
\begin{eqnarray}
\lefteqn{x^{-1}_1\delta\left(\frac{x_2-(z_1-z_{2})}{x_1}\right)
(\tau_{P(z_{1}-z_{2})}(Y_{t}(v, x^{-1}_{2})
\mu^{(2)}_{\lambda, w_{(3)}})(w_{(1)}\otimes w_{(2)})}\nno\\
&&=x^{-1}_1\delta\left(\frac{x_2-(z_1-z_{2})}{x_1}\right)
(Y_{P(z_{1}-z_{2})}(v, x^{-1}_{2})
\mu^{(2)}_{\lambda, w_{(3)}})(w_{(1)}\otimes w_{(2)}),\;\;\;\;\;\;\;\;
\end{eqnarray}
proving (14.36).

\subsection{Proof of  Lemma 14.4}

 Since the vertex operator algebra
is rational, the modules $W_{1}$, $W_{2}$ and $W_{3}$
can be decomposed as
direct sums of finitely many irreducible modules.
Let
\begin{eqnarray}
W_{1}&=&\coprod_{i=1}^{r}M^{(1)}_{i},\\
W_{2}&=&\coprod_{i=1}^{s}M^{(2)}_{i},\\
W_{3}&=&\coprod_{i=1}^{t}M^{(3)}_{i}
\end{eqnarray}
where  $M^{(1)}_{1}, \dots, M^{(1)}_{r}$,
$M^{(2)}_{1}, \dots, M^{(2)}_{s}$ and
$M^{(3)}_{1}, \dots, M^{(3)}_{t}$ are irreducible
$V$-modules.
Let $\eta_{i}$,
$\eta_{j}$ and $\eta_{k}$ be the projections
{}from $W_{2}$ to $M_{i}$, {}from $W_{3}$ to $M_{j}$ and {}from
$W_{2}\boxtimes W_{3}$ to $M_{k}$, respectively,
for any positive integers  $i$, $j$ and $k$
less than or equal to $r$, $s$
and $t$, respectively.
Then for such $i$, $j$ and $k$,
it is clear that $\eta_{k}\circ {\cal Y}\circ (\eta_{i}\otimes \eta_{j})$
is an intertwining operator of type ${M_{k}}\choose
{M_{i}M_{j}}$. Moreover, we have
\begin{equation}
{\cal Y}=\sum_{i=1}^{r}\sum_{j=1}^{s}\sum_{k=1}^{t}
\eta_{k}\circ {\cal Y}\circ (\eta_{i}\otimes \eta_{j}).
\end{equation}
Since $\eta_{k}\circ {\cal Y}\circ (\eta_{i}\otimes \eta_{j})$, $1\le i\le r$,
$1\le j\le s$, $1\le k\le t$, are intertwining operators among irreducible
modules, there exist complex numbers $h_{ijk}$, $1\le i\le r$,
$1\le j\le s$, $1\le k\le t$, such that
$\eta_{k}\circ {\cal Y}\circ (\eta_{i}\otimes \eta_{j})$
are in fact maps {}from
$M_{i}\otimes M_{j}$ to $x^{h{ijk}}M_{k}[[x, x^{-1}]]$, $1\le i\le r$,
$1\le j\le s$, $1\le k\le t$ (see \cite{FHL}). Also
if all irreducible $V$-modules are ${\Bbb R}$-graded,
$h_{ijk}$, $1\le i\le r$,
$1\le j\le s$, $1\le k\le t$, are in fact real numbers.
These facts together with
(15.14) proves the first part of Lemma 14.4.

\subsection{Proof of Lemma 14.5}

We have
\begin{eqnarray}\label{14.5-0}
e^{-m_{1}(\log |z|+i\arg z)}f(z)&=&\sum_{l\in {\Bbb Z}_{+}}a_{l}
e^{(m_{l}-m_{1})
(\log |z|+i\arg z)}\nno\\
&=&a_{1}+e^{(m_{2}-m_{1})(\log |z|+i\arg z)}g(z)
\end{eqnarray}
where
\begin{equation}\label{14.5-1}
g(z)=\sum_{l\in {\Bbb Z}_{+}\setminus \{1\}}a_{l}e^{(m_{l}-m_{2})
(\log |z|+i\arg z)}.
\end{equation}
Since $\sum_{l\in {\Bbb Z}_{+}}a_{l}e^{m_{l}
(\log |z|+i\arg z)}$ converges absolutely to $f(z)$, the right-hand side of
(\ref{14.5-1}) coverges absolutely to a multi-valued function
$g(z)$. Since $m_{l}$, $l\in
{\Bbb Z}_{+}$, are real, the absolute convergence of
$\sum_{l\in {\Bbb Z}_{+}}a_{l}e^{m_{l}
(\log |z|+i\arg z)}$ when $|z|$ is small
 means that $\sum_{l\in {\Bbb Z}_{+}}|a_{l}||e^{m_{l}
(\log |z|}|$ is convergent when $|z|$ is small. Thus
\begin{equation}\label{14.5-2}
\sum_{l\in {\Bbb Z}_{+}}|a_{l}|
|e^{(m_{l}-m_{2})
(\log |z|}|
\end{equation}
 is also convergent when $|z|$ is small. Since $m_{l}\ge m_{2}$, $l\in
{\Bbb Z}_{+}\setminus \{1\}$, we see that the limit of (\ref{14.5-2})
 when $|z|\to 0$ exists. In particular, (\ref{14.5-2}) is bounded.
Since
\begin{eqnarray}
|g(z)|&=&\left|\sum_{l\in {\Bbb Z}_{+}\setminus \{1\}}a_{l}e^{(m_{l}-m_{2})
(\log |z|+i\arg z)}\right|\nno\\
&\le&\sum_{l\in {\Bbb Z}_{+}\setminus \{1\}}|a_{l}||e^{(m_{l}-m_{2})
(\log |z|}|,
\end{eqnarray}
we see that $g(z)$ is also bounded. Taking  $\lim_{|z|\to 0}$ on both sides
of (\ref{14.5-0}), we obtain
$a_{1}=\lim_{|z|\to 0}e^{-m_{1}(\log |z|+i\arg z)}f(z).$
So $a_{1}$ is uniquely determined by $f(z)$ and $m_{1}$. Using
induction, we can
show that for any $l\in {\Bbb Z}_{+}$, $a_{l}$ is uniquely determined by
$f(z)$ and $m_{1}, \dots, m_{l}$.

We now show that $\{m_{l}\}_{l\in {\Bbb Z}_{+}}$ is uniquely determined by
$f(z)$. Assume that there is another sequence  $\{n_{l}\}_{l\in {\Bbb Z}_{+}}$
of strictly increasing real numbers such that $f(z)$ can be expanded as
$\sum_{l\in {\Bbb Z}_{+}}b_{l}e^{n_{l}
(\log |z|+i\arg z)}$. If $\{m_{l}\}_{l\in {\Bbb Z}_{+}}$ and
$\{n_{l}\}_{l\in {\Bbb Z}_{+}}$ are different, we can find $l_{0}$ such
that $m_{l_{0}}\ne n_{l_{0}}$ but $m_{l}=n_{l}$, $l< l_{0}$. Let
\begin{equation}
\tilde{f}(z)=\sum_{l\in {\Bbb Z}_{+}\setminus \{1, \dots, l_{0}-1\}}
b_{l}e^{n_{l}(\log |z|+i\arg z)}.
\end{equation}
Since $f(z)$ is also equal to $\sum_{l\in {\Bbb Z}_{+}}a_{l}e^{m_{l}
(\log |z|+i\arg z)}$ and $m_{l}=n_{l}$, $l< l_{0}$, we have
\begin{equation}\label{14.5-6}
\tilde{f}(z)=\sum_{l\in {\Bbb Z}_{+}\setminus \{1, \dots, l_{0}-1\}}
a_{l}e^{m_{l}(\log |z|+i\arg z)}.
\end{equation}
Assume that $m_{l_{0}}<n_{l_{0}}$. By the discussion above we see that
\begin{equation}\label{14.5-7}
b_{l_{0}}=\lim_{|z|\to 0}e^{-n_{l_{0}}(\log |z|+i\arg z)}\tilde{f}(z).
\end{equation}
But by (\ref{14.5-6}) and the assumption $m_{l_{0}}<n_{l_{0}}$, we see that
the right-hand side of (\ref{14.5-7}) does not exists. Contradiction.
So $\{m_{l}\}_{l\in {\Bbb Z}_{+}}$ and
$\{n_{l}\}_{l\in {\Bbb Z}_{+}}$ are equal.

\subsection{Proof of Lemma 14.6}

For any $w\in W_{1}\boxtimes_{P(z_{1}-z_{2})}W_{2}$,
$$\langle w, \mu^{(2)}_{\gamma(F_{1}; I, F_{2})(\cdot), \cdot}
\rangle_{W_{1}\shboxtr_{P(z_{1}-z_{2})}W_{2}}$$
is an element of $(W'_{(4)}\otimes W_{3})^{*}$. We want to show
\begin{eqnarray}
\lefteqn{\langle z^{-1}_2\delta\left(\frac{x^{-1}_0-x_2}{z_2}\right)
Y_{P(z_{1}-z_{2})}(v, x)w, \mu^{(2)}_{\gamma(F_{1}; I, F_{2})(\cdot), \cdot}
\rangle_{W_{1}\shboxtr_{P(z_{1}-z_{2})}W_{2}}}\nno\\
&&=\tau_{Q(z_{2})}\biggl(
z^{-1}_2\delta\left(\frac{x^{-1}_0-x_2}{z_2}\right)
Y_{t}(v, x)\biggr)(\langle w, \mu^{(2)}_{\gamma(F_{1}; I, F_{2})(\cdot), \cdot}
\rangle_{W_{1}\shboxtr_{P(z_{1}-z_{2})}W_{2}}).\nno\\
&&
\end{eqnarray}
By (15.6) and (15.7) with  $\lambda=\gamma(F_{1}; I, F_{2})(w'_{(4)})$,
we have
\begin{eqnarray}
\lefteqn{\biggl(\tau^{(2)}_{P(z_{1}, z_{2})}
\biggl(x^{-1}_2\delta\left(\frac{x^{-1}_0-z_2}{x_2}\right)
\cdot}\nno\\
&&\hspace{2em}\cdot
Y_{t}((-x^2_0)^{L(0)}e^{-x_0L(1)}v,
x_0)\biggr)(\gamma(F_{1}; I, F_{2})'(w'_{(4)}))\biggr)
(w_{(1)}\otimes w_{(2)}\otimes w_{(3)})\nno\\
&&\quad +x^{-1}_2\delta\left(\frac{z_2-x^{-1}_0}{x_2}\right)
\cdot\nno\\
&&\hspace{2em}\cdot
(\gamma(F_{1}; I, F_{2})'(w'_{(4)}))(w_{(1)}\otimes w_{(2)}\otimes
Y(v, x^{-1}_0) w_{(3)})\nno\\
&&=z^{-1}_2\delta\left(\frac{x^{-1}_0-x_2}{z_2}\right)
\res_{x_{0}^{-1}}(\tau^{(2)}_{P(z_{1}, z_{2})}
(x^{-1}_2\delta\left(\frac{x^{-1}_0-z_{2}}{x_2}\right) \cdot\nno\\
&&\hspace{2em}\cdot Y'_{t}((-x^2_0)^{L(0)}e^{-x_0L(1)}v,
x_0))
(\gamma(F_{1}; I, F_{2})'(w'_{(4)})))
(w_{(1)}\otimes w_{(2)}\otimes w_{(3)})\nno\\
&&\quad +z^{-1}_2\delta\left(\frac{x^{-1}_0-x_2}{z_2}\right)
\res_{x_{0}^{-1}}x^{-1}_2\delta\left(z_{2}-\frac{x^{-1}_0}{x_2}
\right)\cdot\nno\\
&&\hspace{2em}\cdot
(\gamma(F_{1}; I, F_{2})'(w'_{(4)}))(w_{(1)}\otimes w_{(2)}\otimes
Y(v, x^{-1}_0) w_{(3)}).
\end{eqnarray}
By (14.14), the definition of $\gamma(F_{1}; I, F_{2})'$ and the definitions
of $\tau_{P(z_{1}, z_{2})}$ and $\mu_{\cdot, \cdot}^{(2)}$,
the left-hand side of (15.23) is equal to
\begin{eqnarray}
\lefteqn{x^{-1}_2\delta\left(\frac{x^{-1}_0-z_2}{x_2}\right)
\cdot}\nno\\
&&\hspace{2em}\cdot
(\gamma(F_{1}; I, F_{2})'(Y'_{4}((-x^2_0)^{L(0)}e^{-x_0L(1)}v,
x_0)w'_{(4)})))
(w_{(1)}\otimes w_{(2)}\otimes w_{(3)})\nno\\
&&\quad +x^{-1}_2\delta\left(\frac{z_2-x^{-1}_0}{x_2}\right)
\cdot\nno\\
&&\hspace{2em}\cdot
(\gamma(F_{1}; I, F_{2})'(w'_{(4)}))(w_{(1)}\otimes w_{(2)}\otimes
Y(v, x^{-1}_0) w_{(3)})\nno\\
&&=x^{-1}_2\delta\left(\frac{x^{-1}_0-z_2}{x_2}\right)
\mu^{(2)}_{\gamma(F_{1}; I, F_{2})'(Y^{\prime *}_{4}(v,
x^{-1}_0)w'_{(4)}), w_{(3)}}
(w_{(1)}\otimes w_{(2)})\nno\\
&&\quad +x^{-1}_2\delta\left(\frac{z_2-x^{-1}_0}{x_2}\right)
\mu^{(2)}_{\gamma(F_{1}; I, F_{2})'(w'_{(4)}), Y(v, x^{-1}_0) w_{(3)}}
(w_{(1)}\otimes w_{(2)}).
\end{eqnarray}
On the other hand, taking $\res_{x_{1}}$ in (15.8) with
$\lambda=\gamma(F_{1}; I, F_{2})(w'_{(4)})$,
 replacing $v$ in the result  by
$(-x^{-2}_2)^{L(0)}e^{-x^{-1}_2L(1)}v$ and then multiplying the result
by $z^{-1}_2\delta\left(\frac{x^{-1}_0-x_2}{z_2}\right)$, we obatin
\begin{eqnarray}
\lefteqn{\biggl(z^{-1}_2\delta\left(\frac{x^{-1}_0-x_2}{z_2}\right)
Y'_{P(z_{1}-z_{2})}((-x^{-2}_2)^{L(0)}e^{-x^{-1}_2L(1)}v, x_{2}^{-1})
\cdot}\nno\\
&&\hspace{2em}\cdot
\mu^{(2)}_{\gamma(F_{1}; I, F_{2})'(w'_{(4)}),  w_{(3)}}\biggr)
(w_{(1)}\otimes w_{(2)})\nno\\
&&=z^{-1}_2\delta\left(\frac{x^{-1}_0-x_2}{z_2}\right)
\res_{x_{0}^{-1}}(\tau^{(2)}_{P(z_{1}, z_{2})}
(x^{-1}_2\delta\left(\frac{x^{-1}_0-z_{2}}{x_2}\right) \cdot\nno\\
&&\hspace{2em}\cdot Y'_{t}((-x^2_0)^{L(0)}e^{-x_0L(1)}v,
x_0))
(\gamma(F_{1}; I, F_{2})'(w'_{(4)})))
(w_{(1)}\otimes w_{(2)}\otimes w_{(3)})\nno\\
&&\quad +z^{-1}_2\delta\left(\frac{x^{-1}_0-x_2}{z_2}\right)
\res_{x_{0}^{-1}}x^{-1}_2\delta\left(z_{2}-\frac{x^{-1}_0}{x_2}
\right)\cdot\nno\\
&&\hspace{2em}\cdot
(\gamma(F_{1}; I, F_{2})'(w'_{(4)}))(w_{(1)}\otimes w_{(2)}\otimes
Y(v, x^{-1}_0) w_{(3)}).
\end{eqnarray}
{}From the calculations (15.23)--(15.25),
we obtain
\begin{eqnarray}
\lefteqn{\biggl(z^{-1}_2\delta\left(\frac{x^{-1}_0-x_2}{z_2}\right)
Y'_{P(z_{1}-z_{2})}((-x^{-2}_2)^{L(0)}e^{-x^{-1}_2L(1)}v, x_{2}^{-1})
\mu^{(2)}_{\gamma(F_{1}; I, F_{2})'(w'_{(4)}),  w_{(3)}}\biggr)}\nno\\
&&=x^{-1}_2\delta\left(\frac{x^{-1}_0-z_2}{x_2}\right)
\mu^{(2)}_{\gamma(F_{1}; I, F_{2})'(Y^{*}_{4}(v,
x^{-1}_0)w'_{(4)}), w_{(3)}}\nno\\
&&\quad +x^{-1}_2\delta\left(\frac{z_2-x^{-1}_0}{x_2}\right)
\mu^{(2)}_{\gamma(F_{1}; I, F_{2})'(w'_{(4)}), Y(v, x^{-1}_0) w_{(3)}}.
\hspace{10em}
\end{eqnarray}
This formula is equivalent to (15.22).

\subsection{Proof of Lemma 14.9}

Assume that the module spanned by the homogeneous components of the elements of
$W_{2}\boxtimes_{P(z_{2})} W_{3}$ of the form
$w_{(2)}\boxtimes_{P(z_{2})}w_{(3)}$ for all $w_{(2)}\in W_{2}$ and $w_{(3)}\in
W_{3}$ is $W_{0}$. We want to show that $W_{0}=W_{2}\boxtimes_{P(z_{2})}
W_{3}$.

 If not, the quotient module
\begin{equation}
W=(W_{2}\boxtimes_{P(z_{2})}
W_{3})/W_{0}
\end{equation}
 is nontrivial.
Let $P_{W}$ be the projection {}from $W_{2}\boxtimes_{P(z_{2})} W_{3}$ to $W$.
Since $W$ is nontrivial, $P_{W}$ is a nontrivial module map. By Proposition
12.3 in \cite{HL6},
$\eta\to \overline{\eta}\circ \boxtimes_{P(z_{2})}$ is a linear isomorphism
{}from the space of module maps {}from $W_{2}\boxtimes_{P(z_{2})} W_{3}$ to $W$
to
the space of $P(z_{2})$-intertwining maps of type ${W}\choose {W_{2}W_{3}}$.
Thus $\overline{P}_{W}\circ \boxtimes_{P(z_{2})}$ is a nontrivial
$P(z_{2})$-intertwining maps of type ${W}\choose {W_{2}W_{3}}$. But since the
image of $\boxtimes_{P(z_{2})}$ is in $\overline{W}_{0}$,
$\overline{P}_{W}\circ \boxtimes_{P(z_{2})}$ is the trivial map. We have a
contradiction.

\renewcommand{\theequation}{\thesection.\arabic{equation}}
\renewcommand{\therema}{\thesection.\arabic{rema}}
\setcounter{equation}{0}
\setcounter{rema}{0}

\section{Sufficient conditions for the existence of
associativity isomorphisms}

In Section 14, we have constructed the associativity isomorphisms for
$P(\cdot)$-tensor products under the assumption that the vertex
operator algebra has the properties in Propositions 14.1 and 14.7. We
have shown that these properties are in fact also necessary conditions
for the existence of the associativity isomorphisms. In this section,
we give some sufficient conditions for a vertex operator algebra to
have the properties in Propositions 14.1 and 14.7. For concrete vertex
operator algebras, for example, the vertex operator algebras
associated to Wess-Zumino-Novikov-Witten models and minimal models,
instead of proving the properties in Propositions 14.1 and 14.7, we
can prove the sufficient conditions given in this section.

Let $V$ be a rational vertex operator algebra having the property that
all irreducible
$V$-modules are ${\Bbb R}$-graded. We see that any $V$-module must also
be ${\Bbb R}$-graded, that is, the weight of an element in a $V$-module
is always a real number.

Given any $V$-modules $W_{1}$, $W_{2}$, $W_{3}$, $W_{4}$ and $W_{5}$,
let ${\cal Y}_{1}$, ${\cal Y}_{2}$, ${\cal Y}_{3}$ and ${\cal Y}_{4}$
be intertwining operators of type ${W_{4}}\choose {W_{1}W_{5}}$,
${W_{5}}\choose {W_{2}W_{3}}$, ${W_{5}}\choose {W_{1}W_{2}}$ and
${W_{4}}\choose {W_{5}W_{3}}$, respectively. Consider the following
conditions for the product of ${\cal Y}_{1}$ and ${\cal Y}_{2}$ and
for the iterate of ${\cal Y}_{3}$ and ${\cal Y}_{4}$, respectively:

\begin{description}

\item[Convergence and extension property for products]
 There exists an integer $N$
depending only on ${\cal Y}_{1}$ and ${\cal Y}_{2}$, and
for any $w_{(1)}\in W_{1}$,
$w_{(2)}\in W_{2}$, $w_{(3)}\in W_{3}$, $w'_{(4)}\in W'_{4}$, there exist
$j\in {\Bbb N}$, $r_{i}, s_{i}\in {\Bbb R}$, $i=1, \dots, j$, and analytic
functions $f_{i}(z)$ on $|z|<1$, $i=1, \dots, j$,
satisfying
\begin{equation}
\wt w_{(1)}+\wt w_{(2)}+s_{i}>N,\;\;\;i=1, \dots, j,
\end{equation}
such that
\begin{equation}
\langle w'_{(4)}, {\cal Y}_{1}(w_{(1)}, x_{2})
{\cal Y}_{2}(w_{(2)}, x_{2})w_{(3)}\rangle_{W_{4}}
\lbar_{x_{1}= z_{1}, \;x_{2}=z_{2}}
\end{equation}
is convergent when $|z_{1}|>|z_{2}|>0$ and can be analytically extended to
the multi-valued analytic function
\begin{equation}
\sum_{i=1}^{j}z_{2}^{r_{i}}(z_{1}-z_{2})^{s_{i}}
f_{i}\left(\frac{z_{1}-z_{2}}{z_{2}}\right)
\end{equation}
when $|z_{2}|>|z_{1}-z_{2}|>0$.

\item[Convergence and extension property for iterates] There exists an integer
$\tilde{N}$
depending only on ${\cal Y}_{3}$ and ${\cal Y}_{4}$, and
for any $w_{(1)}\in W_{1}$,
$w_{(2)}\in W_{2}$, $w_{(3)}\in W_{3}$, $w'_{(4)}\in W'_{4}$, there exist
$k\in {\Bbb N}$, $\tilde{r}_{i}, \tilde{s}_{i}\in {\Bbb R}$, $i=1, \dots, k$,
and analytic
functions $\tilde{f}_{i}(z)$ on $|z|<1$, $i=1, \dots, k$,
satisfying
\begin{equation}
\wt w_{(2)}+\wt w_{(3)}+\tilde{s}_{i}>\tilde{N},\;\;\;i=1, \dots, k,
\end{equation}
 such that
\begin{equation}
\langle w'_{(4)},
{\cal Y}_{4}({\cal Y}_{3}(w_{(1)}, x_{0})w_{(2)}, x_{2})w_{(3)}\rangle_{W_{4}}
\lbar_{x_{0}=z_{1}-z_{2},\;x_{2}=z_{2}}
\end{equation}
is convergent when $|z_{2}|>|z_{1}-z_{2}|>0$ and can be analytically extended
to the multi-valued analytic function
\begin{equation}
\sum_{i=1}^{k}z_{1}^{\tilde{r}_{i}}z_{2}^{\tilde{s}_{i}}
\tilde{f}_{i}\left(\frac{z_{2}}{z_{1}}\right)
\end{equation}
when $|z_{1}|>|z_{2}|>0$.

\end{description}

\begin{rema}
{\rm
If $V$ is rational, we can always choose $j$,
$r_{i}, s_{i}$,
$i=1, \dots, j$, and $f_{i}(z)$, $i=1, \dots, j$, such that
\begin{equation}
r_{i}+s_{i}=\Delta,\;\;\;i=1, \dots, j
\end{equation}
where
\begin{equation}
\Delta=\wt w_{(1)}+\wt w_{(2)} +\wt w_{(3)}-\wt w'_{(4)}.
\end{equation}
Similarly, for the convergence and extension property for iterates,
we can always choose
$k$, $\tilde{r}_{i}, \tilde{s}_{i}$, $i=1, \dots, k$, and $\tilde{f}_{i}$,
$i=1, \dots, k$,
such that
\begin{equation}
\tilde{r}_{i}+\tilde{s}_{i}=\Delta,\;\;\;i=1, \dots, k.
\end{equation}}
\end{rema}

If for any $V$-modules $W_{1}$, $W_{2}$, $W_{3}$, $W_{4}$ and $W_{5}$ and
 any intertwining operators ${\cal Y}_{1}$ and ${\cal Y}_{2}$
of the types as above, the convergence and extension property for
products holds,
we say that
{\it  the products of the
 intertwining operators for $V$ have the convergence and extension property}.
Similarly we can define what {\it the iterates of the intertwining
operators for $V$ have the convergence and extension property}
means.

We also need the following concept: If a generalized
$V$-module $W=\coprod_{n\in {\Bbb C}}W_{(n)}$ satisfies the condition that
$W_{(n)}=0$ for $n$ whose real real part is sufficiently small, we
say that $W$ is a {\it lower-truncated generalized $V$-module}.

The following theorem is the main result of this section:

\begin{theo}
Assume that $V$ is a rational vertex operator algebra and
that all irreducible $V$-modules are ${\Bbb R}$-graded.  Also assume that $V$
satisfies the following conditions:
\begin{enumerate}

\item Every finitely-generated lower-truncated generalized $V$-module
is a $V$-module.

\item The products or the iterates of
the intertwining operators for $V$ have the convergence and extension property.
\end{enumerate}
Then $V$ has  all the properties in
Propositions 14.1 and 14.7.
\end{theo}
\pf
By the convergence and extension property, $V$ has the properties in
Proposition 14.1. By
Proposition 14.7, we need only to prove that $V$ has the first property in
Proposition 14.7, that is, for any $V$-modules $W_{1}$, $W_{2}$, $W_{3}$,
$W_{4}$ and $W_{5}$, any $z_{1}$ and $z_{2}$ satisfying   (14.7)
and any $P(z_{1})$- and $P(z_{2})$-intertwining maps $F_{1}$ and $F_{2}$
of the types as above, the product
$\gamma(F_{1}; I, F_{2})$ of $F_{1}$ and $F_{2}$ satisfies the
$P(z_{1}-z_{2})$-local grading-restriction condition.

By assumption, there exist $j\in {\Bbb N}$,
$r_{i}, s_{i}\in {\Bbb C}$,
$i=1, \dots, j$, and analytic functions $f_{i}(z)$ on $|z|<1$,
$i=1, \dots, j$, satisfying (16.1) and
 (16.7) such that when $|z_{1}|>|z_{2}|>0$,
(14.1), or equivalently (14.4), is absolutely convergent
and it has an analytic extension in the region $|z_{2}|>|z_{1}-z_{2}|>0$ of
the form
\begin{equation}
\sum_{i=1}^{j}e^{r_{i}\log z_{2}}e^{s_{i}\log (z_{1}-z_{2})}
f_{i}\left(\frac{z_{1}-z_{2}}{z_{2}}\right)
\end{equation}
for any homogeneous $w'_{(4)}\in W'_{4}$,
$w_{(1)}\in W_{1}$, $w_{(2)}\in W_{2}$ and $w_{(3)}\in W_{3}$.
Expanding $f_{i}$, $i=1, \dots, j$, we can write (16.10) as
\begin{equation}
\sum_{i=1}^{j}\sum_{m\in {\Bbb N}}C_{im}(w'_{(4)}, w_{(1)}, w_{(2)}, w_{(3)})
e^{(r_{i}-m)\log z_{2}}e^{(s_{i}+m)
\log (z_{1}-z_{2})}.
\end{equation}
For any $n\in {\Bbb C}$, let
\begin{equation}
a_{n}(w'_{(4)}, w_{(1)}, w_{(2)}, w_{(3)})=
\sum_{-r_{i}+m-1=n}
C_{im}(w'_{(4)}, w_{(1)}, w_{(2)}, w_{(3)}).
\end{equation}
Then (16.11) can be written as
\begin{equation}
\sum_{n\in {\Bbb C}}a_{n}(w'_{(4)},
w_{(1)},
w_{(2)}, w_{(3)})e^{(\Delta+n+1)\log (z_{1}-z_{2})}e^{(-n-1)\log z_{2}}.
\end{equation}
If  $n\in {\Bbb C}$ satisfies
\begin{equation}
n\ne -r_{i}+m-1=-\Delta+s_{i}+m-1
\end{equation}
for any
$i$, $1\le i\le j$, and any $m\in {\Bbb N}$,
then by definition
\begin{equation}
a_{n}(w'_{(4)}, w_{(1)},
w_{(2)}, w_{(3)})=0.
\end{equation}
Since  $m\ge 0$, by (16.1), (16.7) and (16.14)
 we see that for $n\in {\Bbb C}$, if
\begin{equation}
n+1+ \wt w'_{(4)}-
 \wt w_{(3)}\le N,
\end{equation}
 (16.15) holds.

When (14.7) holds,
\begin{equation}
\sum_{n\in {\Bbb C}}a_{n}(w'_{(4)},
w_{(1)},
w_{(2)}, w_{(3)})e^{(-\Delta+n+1)\log (z_{1}-z_{2})}e^{(-n-1)\log z_{2}}
\end{equation}
converges absolutely to (14.1) or (14.4). {}From the definition, we see that
the series (16.17) is of the form (14.50) as a function of $z_{2}$.
For $w'_{(4)}\in W'_{4}$ and $w_{(3)}\in W_{3}$,
let $\beta_{n}(w'_{(4)}, w_{(3)})\in (W_{1}\otimes W_{2})^{*}$ be defined by
\begin{equation}
\beta_{n}(w'_{(4)}, w_{(3)})(w_{(1)}\otimes w_{(2)})
=a_{n}(w'_{(4)}, w_{(1)},
w_{(2)}, w_{(3)})
\end{equation}
for all $w_{(1)}\in W_{1}$ and $w_{(2)}\in W_{2}$.
By definition the series
$$\sum_{\Re(n\in {\Bbb C}}\beta_{n}
(w'_{(4)}, w_{(3)})e^{(\Delta+n+1)\log (z_{1}-z_{2})}e^{(-n-1)\log z_{2}}$$
is absolutely
convergent to $\mu^{(2)}_{\gamma(F_{1}; I, F_{2})'(w'_{(4)}), w_{(3)}}$
and is indexed by sequences of strictly increasing real numbers.
To show that $\gamma(F_{1}; I, F_{2})'(w'_{(4)})$ satisfies the
$P(z_{1}-z_{2})$-local grading-restriction
condition,
we need to calculate
$$(v_{1})_{m_{1}}\cdots (v_{r})_{m_{r}}\beta_{n}(w'_{(4)},
w_{(3)})$$ and its weight for any $r\in {\Bbb N}$, $v_{1}, \dots, v_{r}\in V,$
 $m_{1}, \dots, m_{r}\in {\Bbb Z}$,
$n\in {\Bbb C}$, where $(v_{1})_{m_{1}}, \dots,  (v_{r})_{m_{r}}$, $m_{1},
\cdots, m_{r}\in {\Bbb Z}$, are the components of
$Y'_{P(z_{1}-z_{2})}(v_{1}, x)$, $\dots\ $, $Y'_{P(z_{1}-z_{2})}(v_{r}, x)$,
respectively,  on $(W_{1}\otimes W_{2})^{*}$. For convenience, we
instead calculate
$$(v^{*}_{1})_{m_{1}}\cdots (v^{*}_{r})_{m_{r}}
\beta_{n}(w'_{(4)}, w_{(3)}),$$ where
$(v^{*}_{1})_{m_{1}}, \dots,  (v^{*}_{r})_{m_{r}}$,
$m_{1}, \cdots, m_{r}\in {\Bbb Z}$,  are the components of the opposite
 vertex operators
$Y^{\prime *}_{P(z_{1}-z_{2})}(v_{1}, x)$, $\dots\ $,
$Y^{\prime *}_{P(z_{1}-z_{2})}(v_{r}, x)$, respectively.
By the definition of
$Y'_{P(z_{1}-z_{2})}(v, x)$ and
$$Y^{\prime *}_{P(z_{1}-z_{2})}(v, x)=Y'_{P(z_{1}-z_{2})}
(e^{xL(1)}(-x^{-2})^{L(0)}v, x^{-1}),$$ we have
\begin{eqnarray}
\lefteqn{(Y^{\prime *}_{P(z_{1}-z_{2})}(v, x)
\beta_{n}(w'_{(4)}, w_{(3)}))(w_{(1)}\otimes w_{(2)})=}\nno\\
&&=\res_{y}(z_{1}-z_{2})^{-1}\delta\left(\frac{x-y}{z_{1}-z_{2}}
\right)
(\beta_{n}(w'_{(4)}, w_{(3)}))(Y_{1}(v, y)
w_{(1)}\otimes w_{(2)})
\nno\\
&&\quad +
(\beta_{n}(w'_{(4)}, w_{(3)}))(w_{(1)}\otimes Y_{2}(v, x)
w_{(2)})\nno\\
&&=\res_{y}(z_{1}-z_{2})^{-1}\delta\left(\frac{x-y}{z_{1}-z_{2}}
\right)
a_{n}(w'_{(4)}, Y_{1}(v, y)
w_{(1)}, w_{(2)}, w_{(3)})
\nno\\
&&\quad +
a_{n}(w'_{(4)}, w_{(1)}, Y_{2}(v, x)
w_{(2)}, w_{(3)}).
\end{eqnarray}
On the other hand, since (16.17) is equal to
(14.4)
when (14.7) holds, we have
\begin{eqnarray}
\lefteqn{\res_{y}(z_{1}-z_{2})^{-1}\delta\left(\frac{x-y}{z_{1}-z_{2}}
\right)\sum_{r\in {\Bbb C}}a_{r}(w'_{(4)}, Y_{1}(v, y)
w_{(1)}, w_{(2)}, w_{(3)})\cdot}\nno\\
&&\quad \quad\cdot
 e^{(\Delta+r+1)\log (z_{1}-z_{2})}e^{(-r-1)\log z_{2}}
\nno\\
&&+\sum_{r\in {\Bbb C}}a_{r}(w'_{(4)},
w_{(1)}, Y_{1}(v, y)w_{(2)}, w_{(3)}) e^{(\Delta+r+1)\log (z_{1}-z_{2})}
e^{(-r-1)\log z_{2}}\nno\\
&&=\res_{y}(z_{1}-z_{2})^{-1}\delta\left(\frac{x-y}{z_{1}-z_{2}}
\right)
\langle w'_{(4)}, {\cal Y}_{1}(Y_{1}(v, y)w_{(1)},
x_{1})\cdot\nno\\
&&\quad\quad\cdot
{\cal Y}_{2}(w_{(2)}, x_{2})w_{(3)}\rangle_{W_{4}}
\lbar_{x_{1}
=z_{1}, x_{2}= z_{2}}
\nno\\
&&\quad +
\langle w'_{(4)}, {\cal Y}_{1}(w_{(1)}, x_{1})
{\cal Y}_{2}(Y_{2}(v, x)w_{(2)}, x_{2})w_{(3)}\rangle_{W_{4}}
\lbar_{x_{1}
=z_{1}, x_{2}= z_{2}}\nno\\
&&=\res_{y}(x_{1}-x_{2})^{-1}\delta\left(\frac{x-y}{x_{1}-x_{2}}
\right)
\langle w'_{(4)}, {\cal Y}_{1}(Y_{1}(v, y)w_{(1)},
x_{1})\cdot\nno\\
&&\quad\quad\cdot
{\cal Y}_{2}(w_{(2)}, x_{2})w_{(3)}\rangle_{W_{4}}
\lbar_{x_{1}
=z_{1}, x_{2}=z_{2}}
\nno\\
&&\quad +
\langle w'_{(4)}, {\cal Y}_{1}(w_{(1)}, x_{1})
{\cal Y}_{2}(Y_{2}(v, x)w_{(2)}, x_{2})w_{(3)}\rangle_{W_{4}}
\lbar_{x_{1}
=z_{1}, x_{2}=z_{2}}.
\end{eqnarray}
Using the Jacobi identity for the intertwining operators ${\cal Y}_{1}$ and
${\cal Y}_{2}$,  the right-hand side of (16.20) is equal to
\begin{eqnarray}
\lefteqn{
\langle w'_{(4)}, Y_{1}(v, x-x_{2}){\cal Y}_{1}(w_{(1)},
x_{1})
{\cal Y}_{2}(w_{(2)}, x_{2})w_{(3)}\rangle_{W_{4}}
\lbar_{x_{1}
=z_{1}, x_{2}=z_{2}}}
\nno\\
&&\quad -\langle w'_{(4)}, {\cal Y}_{1}(w_{(1)},
x_{1})Y_{5}(v, x-x_{2})
{\cal Y}_{2}(w_{(2)}, x_{2})w_{(3)}\rangle_{W_{4}}
\lbar_{x_{1}
=z_{1}, x_{2}=z_{2}}
\nno\\
&&\quad +\res_{y}x^{-1}\delta\left(\frac{y-x_{2}}{x}\right)
\langle w'_{(4)}, {\cal Y}_{1}(w_{(1)}, x_{1})Y_{5}(v, y)\cdot\nno\\
&&\quad\quad\cdot
{\cal Y}_{2}(w_{(2)}, x_{2})w_{(3)}\rangle_{W_{4}}
\lbar_{x_{1}
=z_{1}, x_{2}=z_{2}}\nno\\
&&- \res_{y}x^{-1}\delta\left(\frac{x_{2}-y}{x}\right)
\langle w'_{(4)}, {\cal Y}_{1}(w_{(1)}, x_{1})\cdot\nno\\
&&\quad\quad\cdot
{\cal Y}_{2}(w_{(2)}, x_{2})Y_{3}(v, y)w_{(3)}\rangle_{W_{4}}
\lbar_{x_{1}
=z_{1}, x_{2}=z_{2}}\nno\\
&&=\res_{y}x^{-1}\delta\left(\frac{y-x_{2}}{x}\right)
\langle w'_{(4)}, Y_{4}(v, y){\cal Y}_{1}(w_{(1)}, x_{1})\cdot\nno\\
&&\quad\quad\cdot
{\cal Y}_{2}(w_{(2)}, x_{2})w_{(3)}\rangle_{W_{4}}
\lbar_{x_{1}
=z_{1}, x_{2}=z_{2}}\nno\\
&&- \res_{y}x^{-1}\delta\left(\frac{x_{2}-y}{x}\right)
\res_{x_{2}}x_{2}^{n}
\langle w'_{(4)}, {\cal Y}_{1}(w_{(1)}, x_{1})\cdot\nno\\
&&\quad\quad\cdot
{\cal Y}_{2}(w_{(2)}, x_{2})Y_{2}(v, y)w_{(3)}\rangle_{W_{4}}
\lbar_{x_{1}
=z_{1}, x_{2}=z_{2}}.
\end{eqnarray}
Thus
\begin{eqnarray}
\lefteqn{\res_{y}(z_{1}-z_{2})^{-1}\delta\left(\frac{x-y}{z_{1}-z_{2}}
\right)\cdot}\nno\\
&&\quad \quad\cdot \sum_{r\in {\Bbb C}}a_{r}(w'_{(4)}, Y_{1}(v, y)
w_{(1)}, w_{(2)}, w_{(3)})
 e^{(\Delta+r+1)\log (z_{1}-z_{2})}e^{(-r-1)\log z_{2}}
\nno\\
&&+\sum_{r\in {\Bbb C}}a_{r}(w'_{(4)},
w_{(1)}, Y_{1}(v, y)w_{(2)}, w_{(3)})\cdot\nno\\
&&\quad \quad\cdot e^{(\Delta+r+1)\log (z_{1}-z_{2})}
e^{(-r-1)\log z_{2}}\nno\\
&&=\res_{y}x^{-1}\delta\left(\frac{y-x_{2}}{x}\right)
\langle w'_{(4)}, Y(v, y){\cal Y}_{1}(w_{(1)}, x_{1})\cdot\nno\\
&&\quad\cdot
{\cal Y}_{2}(w_{(2)}, x_{2})w_{(3)}\rangle_{W_{4}}
\lbar_{x_{1}
=z_{1}, x_{2}=z_{2}}\nno\\
&&- \res_{y}x^{-1}\delta\left(\frac{x_{2}-y}{x}\right)
\res_{x_{2}}x_{2}^{n}
\langle w'_{(4)}, {\cal Y}_{1}(w_{(1)}, x_{1})\cdot\nno\\
&&\quad\cdot
{\cal Y}_{2}(w_{(2)}, x_{2})Y(v, y)w_{(3)}\rangle_{W_{4}}
\lbar_{x_{1}
=z_{1}, x_{2}=z_{2}}\nno\\
&&=\sum_{m\in {\Bbb Z}}\sum_{l\ge 0}
(-1)^{l}{{m}\choose {l}}x^{-m-1}x_{2}^{l}\langle w_{(4)}, v_{m-l}
{\cal Y}_{1}(w_{(1)}, x_{1})\cdot\nno\\
&&\quad\quad\cdot
{\cal Y}_{2}(w_{(2)}, x_{2})w_{(3)}\rangle_{W_{4}}
\lbar_{x_{1}
=z_{1}, x_{2}=z_{2}}\nno\\
&&\quad -\sum_{m\in {\Bbb Z}}\sum_{ l\ge 0}
(-1)^{l+m}{{m}\choose {l}}x^{-m-1}x_{2}^{m-l}\langle w_{(4)},
{\cal Y}_{1}(w_{(1)}, x_{1})\cdot\nno\\
&&\quad\quad\cdot
{\cal Y}_{2}(w_{(2)}, x_{2})v_{l}
w_{(3)}\rangle_{W_{4}}\lbar_{x_{1}
=z_{1}, x_{2}=z_{2}}\nno\\
&&=\sum_{m\in {\Bbb Z}}\sum_{l\ge 0}
(-1)^{l}{{m}\choose {l}}x^{-m-1}\sum_{r\in {\Bbb C}}a_{r}(v^{*}_{m-l}w'_{(4)},
w_{(1)}, w_{(2)}, w_{(3)})\cdot\nno\\
&&\quad \quad\cdot e^{(\Delta+r+1)
\log (z_{1}-z_{2})}e^{(l-r-1)\log z_{2}}\nno\\
&&\quad -\sum_{m\in {\Bbb Z}}\sum_{ l\ge 0}
(-1)^{l+m}{{m}\choose {l}}x^{-m-1}
\sum_{r\in {\Bbb C}}a_{r}(v^{*}_{m-l}w'_{(4)},
w_{(1)}, w_{(2)}, w_{(3)})\cdot\nno\\
&&\quad \quad\cdot e^{(\Delta+r+1)\log (z_{1}-z_{2})}e^{(-m+l-r-1)\log z_{2}}.
\end{eqnarray}
The intermediate steps in the equality (16.20) holds only when
(14.7) holds. But since
both sides of (16.22) can be analytically extended to
$|z_{2}|>|z_{1}-z_{2}|>0$,
they are equal when
$|z_{2}|>|z_{1}-z_{2}|>0$. By Lemma 14.5, the coefficients of both sides of
(16.22) in powers of $e^{\log z_{2}}$ are equal, that is,
\begin{eqnarray}
\lefteqn{\res_{y}(z_{1}-z_{2})^{-1}\delta\left(\frac{x-y}{z_{1}-z_{2}}
\right)\cdot a_{n}(w'_{(4)}, Y_{1}(v, y)
w_{(1)}, w_{(2)}, w_{(3)})}
\nno\\
&&+a_{n}(w'_{(4)},
w_{(1)}, Y_{1}(v, y)w_{(2)}, w_{(3)})\nno\\
&&=\sum_{m\in {\Bbb Z}}\sum_{l\ge 0}
(-1)^{l}{{m}\choose {l}}x^{-m-1}a_{n+l}(v^{*}_{m-l}w'_{(4)},
w_{(1)}, w_{(2)}, w_{(3)})\nno\\
&&\quad -\sum_{m\in {\Bbb Z}}\sum_{ l\ge 0}
(-1)^{l+m}{{m}\choose {l}}x^{-m-1}
a_{n+m-l}(v^{*}_{m-l}w'_{(4)},
w_{(1)}, w_{(2)}, w_{(3)}).\nno\\
&&
\end{eqnarray}
By (16.17) and (16.23), we obtain
\begin{eqnarray}
\lefteqn{(v^{*}_{m}\beta_{n}(w'_{(4)}, w_{(3)}))(w_{(1)}
\otimes w_{(2)})}\nno\\
&&=\sum_{l\ge 0}
(-1)^{l}{{m}\choose {l}}(\beta_{n+l}(v^{*}_{m-l}w'_{(4)}, w_{(3)}))
(w_{(1)}, w_{(2)})\nno\\
&&\quad -\sum_{ l\ge 0}
(-1)^{l+m}{{m}\choose {l}}(\beta_{n+m-l}(w'_{(4)}, v_{l}w_{(3)}))(w_{(1)},
w_{(2)}).
\end{eqnarray}
By induction,
\begin{eqnarray}
\lefteqn{((v^{*}_{1})_{m_{1}}\cdots (v^{*}_{r})_{m_{r}}\beta_{n}(w'_{(4)},
w_{(3)})
(w_{1)}\otimes w_{(2)})=}\nno\\
&&=\sum_{i\in {\Bbb N}}\sum_{\begin{array}{c}
\mbox{\scriptsize $j_{1}>\cdots >j_{i}$}\\
\mbox{\scriptsize $j_{i+1}>\cdots >j_{r}$}
\\
\mbox{\scriptsize $\{j_{1}, \dots, j_{r}\}=\{1, \dots, r\}$}\end{array}}
\sum_{l_{1},  \dots, l_{r}\ge 0}(-1)^{l_{1}+\cdots +l_{r}+(m_{j_{i+1}}+1)+
\cdots
+ (m_{j_{r}}+1)}\cdot\nno\\
&&\hspace{2em}\cdot
{{m_{j_{1}}}\choose {l_{1}}}\cdots {{m_{j_{r}}}\choose {l_{r}}}
(\beta_{n+m_{j_{i+1}}\cdots +m_{j_{r}}+l_{1}+\cdots +l_{i}-l_{i+1}-
\cdots -l_{r}}\nno\\
&&\hspace{2em}
((v^{*}_{j_{1}})_{m_{j_{1}}-l_{1}}\cdots (v^{*}_{j_{i}})_{m_{j_{i}}-l_{i}}
w'_{(4)}, v_{l_{i+1}}\cdots v_{l_{r}}w_{(3)}))(w_{(1)}\otimes w_{(2)})
\;\;\;\;\;\;\;\;\;\;\;\;
\end{eqnarray}
for any $m_{1}, \dots, m_{r}\in {\Bbb Z}$ and any $v_{1}, \dots,
v_{r}\in V$.
Taking $v=\omega$ and $m=1$ in (16.24), we have
\begin{eqnarray}
\lefteqn{(L'_{P(z_{1}-z_{2})}(0)\beta_{n}(w'_{(4)},
w_{(3)}))(w_{(1)}
\otimes w_{(2)})}\nno\\
&&=(\beta_{n}(L'(0)w'_{(4)}, w_{(3)}))(w_{(1)}
\otimes w_{(2)})\nno\\
&&\quad -(\beta_{n+1}(L'(1)w'_{(4)}, w_{(3)}))(w_{(1)}
\otimes w_{(2)})\nno\\
&&\quad -(\beta_{n}(w'_{(4)}, L(0)w_{(3)}))(w_{(1)}
\otimes w_{(2)})\nno\\
&&\quad +(\beta_{n+1}(w'_{(4)}, L(-1)w_{(3)}))(w_{(1)}
\otimes w_{(2)}).
\end{eqnarray}
When $|z_{1}|>|z_{2}|>|z_{0}|>0$ where $z_{0}=z_{1}-z_{2}$,
\begin{eqnarray}
\lefteqn{-\sum_{n\in {\Bbb C}}a_{n}(L'(1)w'_{(4)},
w_{(1)}, w_{(2)}, w_{(3)})
e^{(\Delta+n+1)\log (z_{1}-z_{2})}e^{(-n-1)\log z_{2}}}
\nno\\
&&\quad +\sum_{n\in {\Bbb C}}a_{n}(w'_{(4)},
w_{(1)}, w_{(2)}, L(-1)w_{(3)}) e^{(\Delta+n+1)\log (z_{1}-z_{2})}
e^{(-n-1)\log z_{2}}\nno\\
&&= -
\langle L'(1)w'_{(4)},
{\cal Y}_{1}(
w_{(1)}, x_{1})
{\cal Y}_{2}(w_{(2)}, x_{2})w_{(3)}\rangle_{W_{4}}
\lbar_{x_{1}
=z_{2}+z_{0}, \; x_{2}=z_{2}}\nno\\
&&\quad +
\langle w'_{(4)}, {\cal Y}_{1}(
w_{(1)}, x_{1})
{\cal Y}_{2}(w_{(2)}, x_{2})L(-1)w_{(3)}\rangle_{W_{4}}
\lbar_{x_{1}
=z_{2}+z_{0}, \; x_{2}=z_{2}}.
\end{eqnarray}
By the commutator formula for $L(-1)$ and intertwining operators and the
$L(-1)$-derivative property for intertwining operators, the right-hand side of
(16.27) is equal to
\begin{eqnarray}
\lefteqn{-
\langle w'_{(4)},
\frac{d}{dx_{1}}({\cal Y}_{1}(
w_{(1)}, x_{1}))
{\cal Y}_{2}(w_{(2)}, x_{2})w_{(3)})\rangle_{W_{4}}
\lbar_{x_{1}
=z_{2}+z_{0}, \; x_{2}=z_{2}}}\nno\\
&&\quad -
\langle w'_{(4)}, {\cal Y}_{1}(
w_{(1)}, x_{1})
\frac{d}{dx_{2}}({\cal Y}_{2}(w_{(2)}, x_{2})w_{(3)})\rangle_{W_{4}}
\lbar_{x_{1}
=z_{2}+z_{0}, \; x_{2}=z_{2}}\nno\\
&&=-\frac{\p}{\p z_{2}}\biggl(
\langle w'_{(4)},
({\cal Y}_{1}(
w_{(1)}, x_{1})
{\cal Y}_{2}(w_{(2)}, x_{2})w_{(3)})\rangle_{W_{4}}
\lbar_{x_{1}
=z_{2}+z_{0}, \; x_{2}=z_{2}}\biggr)\nno\\
&&=-\frac{\p}{\p z_{2}}\Biggl(\sum_{n\in {\Bbb C}}a_{n}(w'_{(4)},
w_{(1)}, w_{(2)}, w_{(3)}) e^{(\Delta+n+1)\log z_{0}}e^{(-n-1)\log z_{2}}
\Biggr)\nno\\
&&=-\sum_{n\in {\Bbb C}}(-n-1)a_{n}(w'_{(4)},
w_{(1)}, w_{(2)}, w_{(3)}) e^{(\Delta+n+1)\log z_{0}}e^{(-n-2)\log z_{2}},
\nno\\
&&
\end{eqnarray}
where $\frac{\p}{\p z_{2}}$ is the partial derivation with respect to $z_{2}$
acting on functions of $z_{0}$ and $z_{2}$.
By the calculations {}from (16.27) to (16.28), the left-hand side of (16.27)
and the right-hand side of (16.28) are
equal when (14.7) holds.
Since they both can be
extended to $|z_{2}|>|z_{1}-z_{2}|>0$, they are equal
 when $|z_{2}|>|z_{1}-z_{2}|>0$. By Lemma 14.5, their coefficients  in powers
of $e^{\log z_{2}}$ are equal.
Thus
\begin{eqnarray}
\lefteqn{ -(\beta_{n+1}(L'(-1)w'_{(4)}, w_{(3)}))(w_{(1)}
\otimes w_{(2)})}\nno\\
&&\quad +(\beta_{n+1}(w'_{(4)}, L(-1)w_{(3)}))(w_{(1)}
\otimes w_{(2)})\nno\\
&&=- a_{n+1}(L'(-1)w'_{(4)},
w_{(1)}, w_{(2)}, w_{(3)})
 + a_{n+1}(w'_{(4)},
w_{(1)}, w_{(2)}, L(-1)w_{(3)})\nno\\
&&=(n+1)a_{n}(w'_{(4)}, w_{(1)}, w_{(2)}, w_{(3)})\nno\\
&&=(n+1)\beta_{n}(w'_{(4)}, w_{(3)})(w_{(1)}\otimes w_{(2)}).
\end{eqnarray}
Substituting (16.29) into (16.26),
we have
\begin{eqnarray}
\lefteqn{(L'_{P(z_{1}-z_{2})}(0)\beta_{n}(w'_{(4)},
w_{(3)}))(w_{(1)}
\otimes w_{(2)})=}\nno\\
&&=(\wt w'_{(4)}+n+1-\wt w_{(3)})(\beta_{n}(w'_{(4)},
w_{(3)}))(w_{(1)}\otimes w_{(2)})\;\;\;\;\;
\end{eqnarray}
for any $w_{(1)}\in W_{1}$ and $w_{(2)}\in W_{2}$.
Thus $\beta_{n}(w'_{(4)}, w_{(3)})$, $n\in {\Bbb C}$,
are weight vectors when
$w_{(3)}$ and $w'_{(4)}$ are homogeneous and their
weights are $\mbox{\rm wt}\ w'_{(4)}+n+1-\mbox{\rm wt}\ w_{(3)}$.
By (16.23) we see that $(v^{*}_{1})_{m_{1}}\cdots
(v^{*}_{r})_{m_{r}}\beta_{n}(w'_{(4)},
w_{(3)})$,
$n\in {\Bbb C}$, are weight vectors when
$v_{1}, \dots, v_{r}$, $w_{(3)}$ and $w'_{(4)}$ are homogeneous and their
weights are
\begin{equation}
-(\mbox{\rm wt}\ v_{1}-m_{1}-1)-\cdots -(\mbox{\rm wt}\ v_{r}-m_{r}-1)
+\mbox{\rm wt}\ w'_{(4)}+n+1-\mbox{\rm wt}\ w_{(3)}.
\end{equation}
By the definition of $\beta_{n}(w'_{(4)}, w_{(3)})$, $\beta_{n}(w'_{(4)},
w_{(3)})=0$ when $$\wt \beta_{n}(w'_{(4)}, w_{(3)})
=\mbox{\rm wt}\ w'_{(4)}+n+1-\mbox{\rm wt}\ w_{(3)}\le N.$$
Thus by (16.25)
\begin{eqnarray}
\lefteqn{(v^{*}_{1})_{m_{1}}\cdots (v^{*}_{r})_{m_{r}}\beta_{n}(w'_{(4)},
w_{(3)})=0}\nno\\
&&\quad\quad \mbox{\rm when}\;\; \wt
(v^{*}_{1})_{m_{1}}\cdots (v^{*}_{r})_{m_{r}}\beta_{n}(w'_{(4)},
w_{(3)})\le N.
\end{eqnarray}
For fixed $n\in {\Bbb C}$, $w_{(3)}\in W_{3}$ and $w'_{(4)}\in W'_{4}$,
let $W_{\beta_{n}(w'_{(4)}, w_{(3)})}$ be the smallest
graded space  containing
$\beta_{n}(w'_{(4)}, w_{(3)})$ and stable under the action of
$Y^{\prime *}_{P(z_{1}-z_{2})}$ or equivalently under $Y'_{P(z_{1}-z_{2})}$.
Then (16.32) shows that
the homogeneous subspace $(W_{\beta_{n}(w'_{(4)}, w_{(3)})})_{(l)}$ of a
fixed wieight $l\in {\Bbb C}$ of $W_{\beta_{n}(w'_{(4)}, w_{(3)})}$
 is $0$ when $l$ is sufficiently small.
Since $\gamma(F_{1}; I, F_{2})'(w'_{(4)})$ satisfies the
$P(z_{1}, z_{2})$-compatibility condition, by the proof of (14.51),
$\beta_{n}$ satisfies the $P(z_{1}-z_{2})$-compatibility condition.
Thus $W_{\beta_{n}(w'_{(4)}, w_{(3)})}$ is a
finitely-generated lower-truncated generalized module.
By assumption, it is in fact a module.
 This proves that for any
$n\in {\Bbb C}$,
$\beta_{n}(w'_{(4)}, w_{(3)})$ satisfies the
$P(z_{1}-z_{2})$-local
grading-restriction
condition. So $\beta_{n}(w'_{(4)}, w_{(3)})$ is an element of
$W_{1}\hboxtr_{P(z_{1}-z_{2})}W_{2}$. Since $\beta_{n}(w'_{(4)}, w_{(3)})$,
$n\in {\Bbb C}$, are all elements of $W_{1}\hboxtr_{P(z_{1}-z_{2})}W_{2}$,
it follows that $\gamma(F_{1}; I, F_{2})'(w'_{(4)})$ satisfies the
$P(z_{1}-z_{2})$-local
grading-restriction condition. \epfv

Combining Theorems 14.10 and 16.2, we obtain:

\begin{theo}
Assume that $V$ is a rational vertex operator algebra and
that all irreducible $V$-modules are ${\Bbb R}$-graded.  Also assume that $V$
satisfies the following conditions:
\begin{enumerate}

\item Every finitely-generated lower-truncated generalized $V$-module
is a $V$-module.

\item The products or the iterates of the intertwining operators
for $V$ have the convergence and extension property.
\end{enumerate}
Then for any $V$-module $W_{1}$, $W_{2}$ and $W_{3}$
and any complex numbers $z_{1}$ and $z_{2}$ satisfying (14.7),
there is a unique isomorphism
${\cal A}_{P(z_{1}), P(z_{2})}^{P(z_{1}-z_{2}), P(z_{2})}$ {}from
$W_1\boxtimes_{P(z_1)}(W_2 \boxtimes_{P(z_2)}W_3)$
to $(W_1 \boxtimes_{P(z_1-z_2)}W_2)\boxtimes_{P(z_2)}W_3$ such that
(14.76) holds. \epf
\end{theo}

We now would like to know when every finitely-generated
lower-truncated
generalized $V$-module
is a module. We consider the following conditions for $V$:

\begin{description}

\item[Condition A] The vertex operator algebra $V$ satisfies the following
conditions:
\begin{enumerate}
\item $V$ is finitely generated by
$v_{1}, \dots$, $v_{k}\in V$.

\item For any $i_{1}, \dots, i_{m}\in {\Bbb Z}$ satisfying
$1\le i_{1}, \dots, i_{m}\le k$ and any $n_{1}, \dots, n_{m}
\in {\Bbb C}$, the operator
\begin{equation}
(v_{i_{1}})_{n_{1}}\cdots (v_{i_{m}})_{n_{m}}
\end{equation}
on $V$ can be written  as a linear combination of the operators
\begin{equation}
(v_{j_{1}})_{l_{1}}\cdots (v_{l_{p}})_{l_{p}},\;\;\;\; p>0, \;1\le j_{1},
\dots, j_{p}\le k, \;l_{1}, \dots, l_{p}\in {\Bbb C},
\end{equation}
having the property that there exist
integers $p_{1}$ and $p_{2}$ satisying $0\le p_{1}\le p_{2}\le p$
such that
\begin{eqnarray}
&&\wt v_{j_{q}}-l_{q}-1<0, \;\;\; 1\le j\le p_{1}\nno\\
&&\wt v_{j_{q}}-l_{q}-1> 0, \;\;\; p_{1}< j\le p_{2}\nno\\
&&\wt v_{j_{q}}-l_{q}-1=0, \;\;\; p_{2}< j\le p.
\end{eqnarray}

\item For any lower-truncated generalized
$V$-module $W$ generated by $w\in W$,
the subspace of $W$ spanned by the elements of the form
\begin{equation}
(v_{i_{1}})_{n_{1}}\cdots (v_{i_{m}})_{n_{m}}w, \;\;
\wt v_{i_{1}}-n_{1}-1=\cdots=\wt v_{i_{1}}-n_{1}-1= 0
\end{equation}
is finite-dimensional.
\end{enumerate}
\end{description}

\begin{propo}
Let $V$ be a vertex operator algebra satisfying Condition A. Then
every finitely-generated
lower-truncated  generalized $V$-module $W$ is a module.
\end{propo}
\pf
We can assume that $W$ is generated by one element $w$.
It is clear {}from the Condition A that given any $n\in {\Bbb C}$,
there are only finitely many elements in $W_{(n)}$ of the form (16.34)
satisfying (16.35). So $W$ is a module.
\epfv

Combining this proposition with Theorem 16.3, we obtain:

\begin{theo}
Let $V$ be a rational vertex operator algebra.
Assume that  $V$ satisfies the Condition A,
that all irreducible $V$-modules are ${\Bbb R}$-graded and
that the products or the iterates of the
intertwining operators for $V$ have the convergence and extension property.
Then for any $V$-module $W_{1}$, $W_{2}$ and $W_{3}$
and any complex numbers $z_{1}$ and $z_{2}$ satisfying (14.7),
there is a unique isomorphism
${\cal A}_{P(z_{1}), P(z_{2})}^{P(z_{1}-z_{2}), P(z_{2})}$ {}from
$W_1\boxtimes_{P(z_1)}(W_2 \boxtimes_{P(z_2)}W_3)$
to $(W_1 \boxtimes_{P(z_1-z_2)}W_2)\boxtimes_{P(z_2)}W_3$ such that
(14.76) holds for any $w_{(1)}\in W_{1}$, $w_{(2)}\in W_{2}$ and
$w_{(3)}\in W_{3}$. \epf
\end{theo}

\noindent {\small \sc Department of Mathematics, University of Pennsylvania,
Philadelphia, PA 19104}

\noindent{\it and}

\noindent {\small \sc Department of Mathematics, Rutgers
University, New Brunswick, NJ 08903} ({\it current address})

\noindent{\em E-mail address}: yzhuang@math.rutgers.edu

\end{document}